\documentclass[12pt,oneside,reqno]{amsart}

\usepackage{appendix}
\usepackage{amsfonts}
\usepackage{amsmath}
\usepackage{amssymb}
\usepackage{amsthm}
\usepackage{enumerate}
\usepackage{graphicx}
\usepackage{hyperref}
\usepackage{mathrsfs}
\usepackage{stmaryrd}
\usepackage{color}
\usepackage{enumitem}
\usepackage{amsmath, amsthm, amssymb, bm, graphicx, hyperref, mathrsfs,geometry}
\newcommand{\be}{\begin{eqnarray}}
\newcommand{\ee}{\end{eqnarray}}
\newcommand{\ce}{\begin{eqnarray*}}
\newcommand{\de}{\end{eqnarray*}}

\newtheorem{theorem}{Theorem}[section]
\newtheorem{proposition}[theorem]{Proposition}
\newtheorem{lemma}[theorem]{Lemma}
\newtheorem{corollary}[theorem]{Corollary}
\newtheorem{remark}[theorem]{Remark}
\newtheorem{definition}[theorem]{Definition}

\newtheorem{assumption}[theorem]{Assumption}

\def\bt{\begin{theorem}}
\def\et{\end{theorem}}
\def\bp{\begin{proposition}}
\def\ep{\end{proposition}}
\def\bl{\begin{lemma}}
\def\el{\end{lemma}}
\def\bc{\begin{corollary}}
\def\ec{\end{corollary}}
\def\bd{\begin{definition}}
\def\ed{\end{definition}}
\def\br{\begin{remark}}
\def\er{\end{remark}}
\def\bx{\begin{Examples}}
\def\ex{\end{Examples}}
\def\ba{\begin{assumption}}
\def\ea{\end{assumption}}

\def\s{\sigma}

\def\[{{\Big[}}
\def\]{{\Big]}}
\def\<{{\langle}}
\def\>{{\rangle}}
\def\({{\Big(}}
\def\){{\Big)}}

\def\geq{\geqslant}
\def\leq{\leqslant}

\def\no{\nonumber}

\def\min{{\mathord{{\rm min}}}}
\def\max{{\mathord{{\rm max}}}}

\def\cB{{\mathcal B}}
\def\cC{{\mathcal C}}
\def\cD{{\mathcal D}}

\def\cF{{\mathcal F}}
\def\cG{{\mathcal G}}
\def\cH{{\mathcal H}}

\def\cK{{\mathcal K}}
\def\cL{{\mathcal L}}

\def\cO{{\mathcal O}}
\def\cP{{\mathcal P}}

\def\cS{{\mathcal S}}

\def\cX{{\mathcal X}}

\def\mD{{\mathbb D}}
\def\mE{{\mathbb E}}

\def\mK{{\mathbb K}}

\def\mN{{\mathbb N}}

\def\mP{{\mathbb P}}
\def\mQ{{\mathbb Q}}
\def\mR{{\mathbb R}}

\def\sF{{\mathscr F}}

\def\Om{\Omega}

\allowdisplaybreaks

\begin{document}

	\title{Integration by Parts Formulas of Mckean-Vlasov SDEs with Jumps and Some Applications}
	\author{Yao Chen$^{1}$, Jiagang Ren$^{1}$, Hua Zhang$^{2,*}$}

\dedicatory{$^{1}$School of Mathematics, Sun Yat-Sen University,\\
Guangzhou, Guangdong 510275, P.R.China\\
$^{2}$School of Statistics and Data Science, Jiangxi University of Finance and Economics,\\
Nanchang, Jiangxi 330013, P.R.China\\
Emails: Y. Chen: cheny2378@mail2.sysu.edu.cn\\
J. Ren: renjg@mail.sysu.edu.cn\\
H. Zhang: zh860801@163.com}

\thanks{This work is supported by National Natural Science Foundation of China (Grant Nos. 12261038 and 12371152), and Natural Science Foundation of Jiangxi Province (Grant Nos. 20232BAB201004 and 20242BAB23003).}

\thanks{$^*$Corresponding Author.}

	
\date{}


	\begin{abstract}
		In this article, we establish integration by parts formulas for the solutions of McKean-Vlasov stochastic differential equations with jumps under elliptic coefficients. The derived formulas accommodate both derivatives with respect to real-valued variables and measure-valued variables, interpreted through the Lions' derivative. As applications,
		we obtain estimates  for the derivatives of the density functions of the McKean-Vlasov SDEs, and relying on the  integration by parts formulas, we subsequently prove the existence and uniqueness of classical solutions to the associated PDEs with irregular terminal conditions.
	\end{abstract}
	
	\maketitle

	\section{Introduction}
	Stochastic differential equations whose coefficients depend on the distribution of the solution, commonly known as McKean-Vlasov stochastic differential equations (MVSDEs in short), have become a fundamental framework in probability theory for modeling nonlinear and nonlocal interactions. Their origin can be traced back to the seminal works of Kac \cite{MK,MK1}, where mean-field limits and propagation of chaos were introduced to derive kinetic-type evolution equations. Since then, MVSDEs have played an essential role in the analysis of interacting particle systems, statistical physics, nonlinear Fokker-Planck equations, and mean-field games. Important contributions in this direction include the works of Lasry and Lions on mean-field games \cite{LL}, Lions' lectures at Coll\`ege de France  \cite{LP}, and the lecture notes of Cardaliaguet \cite{C}. These MVSDEs provide a probabilistic representation to the solutions of a class of nonlinear PDEs, a special case of which was first studied by McKean \cite{HP}.
	
	A central development in the modern theory of MVSDEs is the introduction of a differential calculus on the space of probability measures, initiated by Lions \cite{LP}. This measure-valued calculus, often referred to as the Lions' derivative, has enabled precise analysis of differentiability for functionals of probability measures and has facilitated the derivation of PDEs posed on the Wasserstein space. Building on this framework, \cite{BR,HL} established  second-order It\^o's formulas for mean-field It\^o's processes and identified the fully nonlinear, nonlocal PDEs associated with McKean-Vlasov dynamics. These works demonstrate that the differentiability of the flow with respect to both spatial and measure variables is central for understanding the analytic properties of the associated PDEs.

	Another important aspect concerns regularity and smoothing properties of McKean-Vlasov dynamics. Even when the terminal condition is irregular, the interaction between diffusion and distribution dependence can create smoothing effects, ensuring the existence of density functions with regularity properties. A notable contribution in this direction is due to Crisan and McMurray \cite{CM}, who established an integration by parts formula for MVSDEs with uniformly elliptic diffusion coefficients using Malliavin calculus, thereby obtaining gradient estimates of densities and proving existence and uniqueness of classical solutions to the corresponding PDEs with non-differentiable terminal conditions.
	
	The motivation to study MVSDEs with jumps arises  on the one hand from previous study of MVSDEs driven by Brownian motion \cite{CM}, but on the other hand also from the fact that jump processes can describe some financial phenomenon.  For instance, many works on jumps in finance build upon the seminal jump-diffusion model of Merton \cite{MR} and the more general framework of L\'evy processes \cite{CP}. The need to incorporate such discontinuous processes into McKean-Vlasov models has led to the development of MVSDEs with jumps \cite{DL}, which are capable of capturing systemic risk and large market moves in a multi-agent setting. The presence of jumps destroys classical smoothing properties and introduces nonlocal behaviors. Existing works on MVSDEs with jumps include the study of kinetic limits \cite{DJ}, quasilinear SPDEs of McKean-Vlasov type \cite{KK}. More recently, \cite{HL} investigated  SDEs with jumps and their relations with nonlocal PDEs, highlighting new challenges specific to jump-driven dynamics. It is also worth noting that Li in \cite{L} also studied mean-field backward stochastic differential equations with jumps and established a probabilistic representation for the solutions of associated nonlocal quasi-linear integral-PDEs.
	
	The research on the regularity of the laws of SDEs with jumps has been carried
	out for a long time. Especially, after Malliavin announced the breakthrough article \cite{MP} which established stochastic calculus of variations, Bichteler, Gravereaux and Jocod
	in \cite{KB} rapidly proposed an extension of the Malliavin calculus to the case of stochastic differential equations with jumps. Many subsequent studies develop local operators that act on either the size of jumps (cf. \cite{KB}, etc.) or the instants of jumps (cf. \cite{CPE}, etc.). Other approaches introduce finite difference operators together with Fock space representations (cf. \cite{PJ}, etc.).
	
	Later, in a series of works, Bouleau and Denis present systematically the lent particle method which
	simplifies the previous approaches, (see \cite{BL} and references therein). The simplicity
	and efficiency of this method are also shown there with many illustrative examples
	connected to stochastic calculus with L\'evy processes. \cite{RZ2019, RZ2023} make use of this method to obtain the regularity of the laws of SDEs with jumps in under some kind of H\"ormander's conditions. In this work, by the lent partial method, we derive the regularity of the laws of MVSDEs with jumps and some density estimates  under an   ellipticity assumption.
	
    As an application of this regularity we prove the existences and uniqueness of the classical solutions of a class of nonlocal nonlinear PDEs on Wasserstein space (see \cite{BL} below) and obtain a stochastic representation of the solutions even when the initial condition is not smooth. This complements \cite{HL} which concerns the initial condition is smooth.

    Let us give more details. Let $(\Omega,\cF, \{\cF_r\}, \mP)$ be a filtered  probability space satisfying the usual condition (the precise space we will work with will be introduced later), and  suppose the coefficients  $ c:  \mR^N\times \mR^m \times \cP_2(\mR^N) \rightarrow \mR^N$ and $b:\mR^N\times \cP_2(\mR^N) \rightarrow \mR^N $ are given.  We consider the following  MVSDE with jumps:
	\ce
	X_{r,t}^{\theta}&=&\theta+\int_r^t b(X_{s}^{\theta},[X_{s}^{\theta}])ds+\int_r^t \int_{\Xi} c(X_{s-}^{\theta},u,[X_{s-}^{\theta}])\widetilde{N}(ds,du),
	\de
where  $\widetilde{N}$ is the compensated  process of some $\{\cF_r\}$-Poisson point process $N$ on
a measure space $(\Xi,\cG, \lambda)$ with compensator $\nu(dt,du)=\lambda(du)dt$ .
	Let the initial time $r\in [0,T]$ and the initial condition $\theta \in L^2(\cF_r;\mR^N)(:=L^2(\Omega,\cF_r,\mP;\mR^N))$, which is a square-integrable random variable.  (Hence, in particular, $\theta$ is independent of the martingale ${\widetilde{N}((r, t] \times A)}$ for all $A \in \Xi$.)  Throughout, we denote by $[\xi]$ the law of a random variable $\xi$ and $\cP_2(\mR^N)$ the set of probability measures on ~$\mR^N$ with finite second moment.   Without any loss of generality we may and will suppose that
$\cF_0\subset \cF$ satisfies the following conditions:
	\begin{enumerate}
		\item the Poisson random measure $N$ is independent of $\cF_0$;
		\item $\cF_0$ is rich enough such that $\cP_2(\mR^N)=\{ P_{\vartheta}, \vartheta \in L^2(\cF_0;\mR^N)\}, N\geq 1$.
	\end{enumerate}
	If the initial time $r=0$, $X_{r,t}^{\theta}$ will be abbreviated as $X_t^{\theta}$ and the above MVSDE with jumps takes the following form:
	\be\label{SDE1}
	X_{t}^{\theta}&=&\theta+\int_0^t b(X_{s}^{\theta},[X_{s}^{\theta}])ds+\int_0^t \int_{\Xi} c(X_{s-}^{\theta},u,[X_{s-}^{\theta}])\widetilde{N}(ds,du).
	\ee

	 In \cite{CM}, Crisan and McMurray developed integration by parts formulas for parameters $x$ and $\mu$, 
	  and investigated the regularity of the solutions of MVSDEs driven by Brownian motion using Malliavin calculus. In this work, we first derive  integration by parts formulas using lent particle method, and subsequently establish the regularity properties and associated estimates for the density function.
	
	 In turn, these integration by parts formulas enable us to use MVSDEs to define the solution of a class of nonlocal integral-PDEs that has the form
	 \be\label{pde1}
\begin{cases}
	 (\partial_t-\cL)U(t,x,[\theta])=0 &\quad \text{for} \quad (t,x,[\theta])\in (0,T]\times\mR^N\times\cP_2(\mR^N ),\\
	 U(0,x,[\theta])=g(x,[\theta]) &\quad\text{for} \quad (x,[\theta])\in \mR^N\times\cP_2(\mR^N ),
\end{cases}
	 \ee
	 where $g:\mR^N\times\cP_2(\mR^N )\to \mR$ (not necessarily smooth) and the nonlocal operator ~$\cL$ acts on sufficiently enough functions ~$F:\mR^N\times\cP_2(\mR^N )\to \mR^N$ and is defined by
	 \ce
	 \cL F(x,[\theta])&=&
	 \mE[\partial_{x} F(x,[\theta])b(x,[\theta])] + \mE[\partial_{\mu} F(x,[\theta]) b(x,[\theta])] \\
	 &&+\int_{\Xi} \mE[F(x+c(x,u,[\theta]))-F(x,[\theta])-\partial_xF(x,[\theta]) c(x,u,[\theta]) ]\lambda(du)\\
	 &&+ \int_\Xi\mE[  F(x,[\theta+ c(x,u,[\theta])])-F(x,[\theta])-\partial_\mu F(x,[\theta])  c(x,u,[\theta]) ] \lambda(du).
	 \de
	The last two terms in the description of $\cL F(x,[\theta])$ involve the derivative with respect to the measure variable which was first introduced by Lions \cite{C} (the detailed definition will be presented in Section 2).
	
	It is a classical result established in \cite{F,H} that linear parabolic PDEs on $[0,T]\times \mR^N$
	admit classical solutions under uniform ellipticity or H\"{o}rmander conditions, even with non-differentiable initial data.   In this paper, we investigate whether analogous results hold for PDE (\ref{pde1}) under an   ellipticity condition. Specifically, we prove the existence of classical solutions to PDE (\ref{pde1}) when the initial function
	$g$ is not differentiable. Our approach employs a probabilistic representation for the classical solution of PDE (\ref{pde1}), which is expressed in terms of a functional of $X_t^{\theta}$
	and the solution of the following decoupled equation:
	\ce
	X_{r,t}^{x,[\theta]}=x+ \int_r^t b(X_{s}^{x,[\theta]},[X_{s}^{\theta}]) ds+\int_r^t \int_{\Xi} c(X_{s-}^{x,[\theta]},u,[X_{s-}^{\theta}])\widetilde{N}(ds,du).
	\de
	If $r=0$, we abbreviate it as $X_t^{x,[\theta]}$, and thus
	\be \label{SDE}
	X_{t}^{x,[\theta]}=x+ \int_0^t b(X_{s}^{x,[\theta]},[X_{s}^{\theta}]) ds+\int_0^t \int_{\Xi} c(X_{s-}^{x,[\theta]},u,[X_{s-}^{\theta}])\widetilde{N}(ds,du).
	\ee
Here, the  decoupled equation means the coefficients $b$ and $c$ do not depend on the distribution of $X^{x,\theta}_{t}$, but rather depend
on that of $X_t^{\theta}$. In the following, we consider a certain class of functions $g:\mR^N\times \cP_2(\mR^N)\to \mR$
	\ce
	U(t,x,[\theta]):=\mE[g(X_t^{x,[\theta]},[X_t^{\theta}])],
	\de
	where $(t,x,[\theta])\in [0,T]\times\mR^N\times \cP_2(\mR^N)$. It can be seen from Theorem \ref{pdesolution} that $U(t,x,[\theta])$ solves  PDE (\ref{pde1}).

In \cite{HL}, Hao and Li established that the map $(x,[\theta])\mapsto \mE[g(X_t^{x,[\theta]},[X_t^{\theta}])]$ is differentiable, under the assumption that the initial value $g$ is differentiable. Given the sufficient smoothness of $g$, it is unnecessary to impose non-degeneracy conditions on the coefficients. In contrast to their work, we prove the differentiability of $(x,[\theta])\mapsto \mE[g(X_t^{x,[\theta]},[X_t^{\theta}])]$ under the  ellipticity condition, even when the initial value $g$ is assumed to be non-differentiable. This result is achieved through an extension of the integration by parts formulas established in Section 4.

	The organization of this paper is as follows. The notations, the lent particle method and the basic results related to MVSDEs with jumps are described in Section 2. Section 3 is devoted to the differentiability of solutions to MVSDEs. We present results concerning their derivatives with respect to the parameter $(x,\mu)$, as well as their Malliavin derivatives.  Section 4 presents the integration by parts formula with jumps. The key step involves using the lent particle method to derive a fundamental identity relating the carr\'e du champ operator and the Malliavin derivative,
which then lays the groundwork for the formula. The smoothness of densities of the solutions to MVSDEs is discussed in Section 5. In Section 6, we turn our attention to the associated PDEs when the initial condition $g$ is not differentiable. Finally, the proofs of some technical results are given in Appendix.
	
	Throughout the paper, unimportant constants will be denoted by $C$ with or without index, whose values may vary from place to place.
	\section{Preliminary}
	\subsection{Notation}

	Let $\cP_2(\mR^N)$ be the set of measures on $(\mR^N, \cB(\mR^N))$ with finite second moments. We equip $\cP_2(\mR^N)$ with the $2$-Wasserstein metric, i.e., for $\mu, \nu\in \cP_2(\mR^N)$, we set
	\ce
\begin{gathered}
	W_{2}(\mu,\nu):=\operatorname*{inf}\left\{\left(\int_{\mathbb{R}^{2N}}|x-y|^{2} m(dx,dy)\right)^{\frac{1}{2}}, m\in\mathcal{P}_{2}(\mathbb{R}^{2N}), \mathrm{such~that}\right. \\
	m(A\times\mathbb{R}^N)=\mu(A),~A\in\mathcal{B}(\mathbb{R}^N),~m(\mathbb{R}^N\times B)=\nu(B),~B\in\mathcal{B}(\mathbb{R}^N)\Biggr\}.
\end{gathered}
	\de
	
	 We denote the $L^p$ norm on $(\Omega,\sF,\mP)$ by $\|\cdot\|_p$ for $p\geq1$ and we also introduce the space $\cS_t^p$ of c\`{a}dl\`{a}g  $\{\cF_s\}$-adapted process $\varphi$ on $[0,t]$, satisfying
	\ce
	\|\varphi\|_{\cS_t^p}=(\mE[\sup_{s\in[0,t]}|\varphi_s|^p])^{ \frac 1 p}<\infty.
	\de
	\subsection{Set-up}

	Let us first specify the general set-up in which we will work. We follow \cite{BL}, to which we refer the reader for more details.
	\subsubsection{Dirichlet Structure on the bottom space}

	We assume that ~$\Xi$ is a separable Hausdorff space, $\cG$ its Borel ~$\sigma$-algebra such that ~$\{x\}$ belongs to $\cG$ for any ~$x\in \Xi$, $\nu$ a ~$\sigma$-finite diffuse measure on $\cG$. Let $(\mathbf{d},e)$ be a local symmetric Dirichlet form on $L^2(\nu)$ which admits a carr\'{e} du champ operator $\gamma$. The structure $(\Xi,\cG,\nu,\mathbf{d},\gamma )$ is called the bottom structure.
	
	The main property of local Dirichlet forms we shall use is the Energy Image Density (EID) Property, which is introduced in \cite[ Chapter VI, section 1]{BN}.

    We denote by $\mN^*$ the set of all positive integer.
	\bd
	We say that the Dirichlet form $(\mathbf{d},e) $ satisfies (EID) if for any $N\in \mN^*$ and any $\mR^N$-valued function $U$ whose components are in $\mathbf{d}$, we have
	\ce
	U_{*} [(\det \gamma[U,U^*])\cdot \nu]  \ll \lambda^N,
	\de
	where $\det$ denotes the determinant, $U_{*}\nu$ the image measure by $U$ of the measure $\nu$, and $\lambda^N $ is the $N$-dimensional Lebesgue measure on $(\mR^N,\cB(\mR^N))$.
	\ed
	 We assume that the bottom space satisfies the following hypotheses.
	 \ba
	 \begin{enumerate}
	 	\item The generator $a$ of the Dirichlet form $(\mathbf{d},e) $ with domain $\cD(a)\subset \mathbf{d}$ satisfies the following property: there exists a subspace $L$ of $\cD(a) \cap L^1(\nu)$  such that for any $f\in L$, $\gamma[f]\in L^2(\nu)$.
	 	\item The structure $(\Xi,\cG,\nu,\mathbf{d},\gamma )$ satisfies (EID).
	 	\item $\Xi$ admits a partition of the form $\Xi:=B \cup ( \cup_{k=1}^\infty A_k)$, where for all $k$, $A_k \in \cG$ with $\nu(A_k) <\infty$ and $\nu(B)=0$, in such a way that for any $k \in \mN^*$ there is a local Dirichlet structure with carr\'{e} du champ operator
	 	\ce
	 	\cS_k=(A_k, \cG|_{A_k}, \nu|_{A_k}, \mathbf{d}_k, \gamma_k),
	 	\de
	 	such that for all $\psi \in\mathbf{d}$, $\psi|_{A_k} \in \mathbf{d}_k$ and $\gamma[\psi]|_{A_k}=\gamma[\psi|_{A_k}]$.
	 	\item Any finite product of structures $\cS_k$ satisfies (EID).
	 	\item $\mathbf{d} $ is a separable Hilbert space equipped with the  scalar product
	 	\ce
	 	(\psi,\phi)_\mathbf{d}=(\psi,\phi)_{L^2(\nu)}+e(\psi,\phi).
	 	\de
	 \end{enumerate}
	 \ea
	
	Since $\mathbf{d}$ is separable, the bottom Dirichlet structure admits a gradient operator, i.e., there exists a
	separable Hilbert space $H$ and a linear map $\nabla$ from d into $L^2(\nu;H)$ such that
	\ce
	\gamma[\phi]=\|\nabla \phi \|_H^2,\quad\forall u\in\mathbf{d}.
	\de
	More concretely, let $R=[0,1]$, $\mathcal{R}$ the Lebesgue $\s$-algebra, and $\rho$ the Lebesgue measure. We can take
	\ce
	H=L_0^{2}(R,\mathcal{R},\rho)=\{g\in L^{2}(R,\mathcal{R},\rho);\int_{R}g(r)\rho(dr)=0\}.
	\de
	 The corresponding gradient will be denoted by $\flat$, which we assume without any loss of generality that constants belong to $\mathbf{d}_{loc}$ (see \cite[page 15]{BL})and so that l$^\flat=0.$

	\subsubsection{Dirichlet structure on the upper space}
	From now on we set $X=\mathbb{R}_+\times\Xi, \cX=\mathcal{B}(\mathbb{R}_+)\times\mathcal{G}$ and $\mu=dt\times\nu.$ Define the Dirichlet structure on $(X,\mathcal{X},\mu)$ to be the product of the trivial one on $(L^2(\mathbb{R}_+,dt),0)$ and $(\mathbf{d},e)$ and we keep the same notations $\mathbf{d},e,b,\gamma$ and $\mathbf{a}$, etc, which act only on the second variable for operators corresponding to this new Dirichlet form.
	
	Denote by $\Omega$ the configuration space of $X$, i.e.,
	\ce
	 \Omega&:=&\{\omega=\sum_{i=1}^{\infty} \varepsilon_{x_{i}};  x_{i} \in X,  \forall i, \text{and}\quad \sharp \{K \cap\{x_{i}, i=1,2, \ldots\}\}<\infty  \\
 &&\text{for any compact} \quad  K \subset X \} .
	\de
	Here $\varepsilon_x$ stands for the Dirac measure at $x$ and $\sharp(A)$ is the cardinal of set $A$. Then, due to \cite[page 51]{BL}, it is possible to totally order $\Omega$ and we shall fix and denote by $ \prec$ such a total order relation. Then $\Omega$ can be rewritten in a unique way as
	
	\ce
	\Omega:=\{\omega=\sum_{i=1}^\infty\varepsilon_{x_i};\quad x_i\in X,\quad\forall i,\quad\text{and}\quad x_1\prec x_2\prec\cdots\prec x_n\prec\cdots\}.
	\de
	Let $N$ be the Poisson random measure with intensity $\nu$ defined on $(\Omega,\mathcal{F},\mathbb{P})$ where $N(\omega)=\omega,\mathcal{F}$
	the $\sigma$-algebra generated by $N$ and $\mathbb{P}$ the law of $N.$
	
	We now introduce the creation and annihilation operators $\varepsilon^{+}$ and $\varepsilon^{-}$:
	\ce
	\forall (t,y)\in X,\quad\forall\omega\in\Omega,\\
	\varepsilon^{+}_{(t,y)}(\omega)=\omega 1_{\{(t,y)\in supp\omega\}}+(\omega+\varepsilon_{(t,y)})1_{\{(t,y)\notin supp\omega\}},\\
	\forall (t,y)\in X,\quad\forall\omega\in\Omega,\\
	\varepsilon^{-}_{(t,y)}(\omega)=\omega 1_{\{(t,y)\notin supp\omega\}}+(\omega_1-\varepsilon_{(t,y)})1_{\{(t,y)\in supp\omega\}}.
	\de
	Denote $\mathbb{P}_{N}:=\mathbb{P}(d\omega)N_{\omega}(dt,dy)$. Then it is well known (see \cite[Lemma 4.2]{BL}) that the map $(\omega,(t,y))\mapsto(\varepsilon^{+}_{(t,y)}(\omega))$ sends $\mathbb{P}_{N}$-negligible sets to $\mathbb{P}\times\mu$-negligible ones, and the map $(\omega,(t,y))\mapsto(\varepsilon^{-}_{(t,y)}(\omega),(t,y))$ sends $\mathbb{P}\times\mu$-negligible sets to $\mathbb{P}_{N}$-negligible ones.
	
	If $N(\omega)=\sum_{i=1}^\infty\varepsilon_{x_i}$, then define
	\ce
	N\odot\rho(\omega,\hat{\omega}):=\sum_{i=1}^\infty\varepsilon_{(x_i,r_i(\hat{\omega}))},
	\de
	where $(r_i)$ is a sequence of i.i.d random variables independent of $N$ whose common law is $\rho$ and which are defined on some probability space $(\widehat{\Omega},\widehat{\mathcal{F}},\widehat{\mathbb{P}}).$ Hence $N\odot\rho$ is defined on the product probability space $(\Omega,\mathcal{F},\mathbb{P})\times(\widehat{\Omega},\widehat{\mathcal{F}},\widehat{\mathbb{P}}).$ It is a Poisson measure on $X\times R$ with compensator $\mu\times\rho$ which is called the marked Poisson measure.
	
Starting from the Dirichlet structure on the bottom space, one can construct a Dirichlet structure on the upper space. We denote
the associated semigroup by $\{T_t\}$, its generator by $A$ and the domain of $A$ by $\cD(A)$. For details we refer to \cite[Section 4.3]{BL}.

	\subsubsection{Sobolev Spaces on the Bottom Space}
	Let $E$ be a separable Hilbert space. We denote by $S(E)$ the set of $E$-valued functions defined on $\Xi$ such that there exist $n\in\mathbb{N}^*,{e}_1,\ldots,{e}_n$ in $E$ and $\varphi_1,\ldots,\varphi_n$ in $\mathbf{d}$ with
	\ce
	u=\sum_{{i=1}}^n\varphi_i {e}_i.
	\de
	If $u=\sum_{i=1}^n\varphi_i{e}_i$ belongs to $S(E)$, we can define its gradient $u^\flat=\sum_{i=1}^n\varphi_i^\flat {e}_i$ as an element of $L^2(\nu;H\otimes E).$ We denote by $\mathbf{d}(E)$ the completion of $S(E)$ with respect to the norm:
	\ce
\|u\|_{\mathbf{d}(E)}^2:=\|u\|_{L^2(\nu)}+\|u^\flat\|_{L^2(\nu;H\otimes E)}^2.
\de
	We will use the following assumption.
	\begin{assumption}\label{2.4}
	 There exists a dense vector subspace $\mathbf{d}_0 \subset \mathbf{d}$ such that each element $u$ in $\mathbf{d}_0$ is such that
		\begin{enumerate}
			\item $u \in \bigcap_{p \geqslant 2} L^p(\nu)$;
			\item  $u$ is infinitely differentiable in the sense that $u^\flat \in \mathbf{d}(H)$, $u^{2\flat} = (u^\flat)^\flat \in \mathbf{d}(H^{\hat{\otimes}2})$, $\ldots$, $u^{(k+1)\flat} = (u^{k\flat})^\flat \in \mathbf{d}(H^{\hat{\otimes}(k+1)})$, $\ldots$;
			\item for all $k \in \mathbb{N}^*$, $u^{k\flat} \in \bigcap_{p \geqslant 2} L^p(\nu; H^{\hat{\otimes}k})$.
		\end{enumerate}
		\end{assumption}
		We introduce
		\ce
		\mathbf{d}_0(E) := \{ u = \sum_{i=1}^n \varphi_i {e}_i \in S(E) \mid {e}_i\in E, \varphi_i \in \mathbf{d}_0, i = 1, \ldots, n \}.
		\de
		Now we can define the Hilbert-valued Sobolev spaces on the bottom space.

	\bd
     Let $k \in \mathbb{N}^*$, $p \geqslant 2$. We denote by $\mathbf{d}^{k,p}(E)$ the completion of $\mathbf{d}_0(E)$ w.r.t the norm
	\ce
	\|u\|_{k,p} := \|u\|_{L^p(\nu;E)} + \|u^\flat\|_{L^p(\nu;H\hat{\otimes}E)} + \cdots + \|u^{k\flat}\|_{L^p(\nu;H^{\hat{\otimes}k}\hat{\otimes}E)}.
	\de
	We also set
	\ce
	\mathbf{d}^{k,\infty}(E) := \bigcap_{p \geqslant 2} \mathbf{d}^{k,p}(E) \quad \text{and} \quad \mathbf{d}^\infty(E) := \bigcap_{k \in \mathbb{N}^*} \mathbf{d}^{k,\infty}(E).
	\de
	For simplification of notations, we denote $\mathbf{d}^{k,p}$ for $\mathbf{d}^{k,p}(\mathbb{R})$.
	\ed
	\subsubsection{Sobolev Spaces on Poisson Space}
	Set
	\ce
	\mD_0:=\{  \varphi ( \widetilde{N}(f_1),\cdots , \widetilde{N}(f_k));\varphi \in \cC_c^\infty(\mR^k),f_1,\cdots,f_k \in  \mathbf{d}^{\infty},k\in \mN^*\}
	\de
	and
	\ce
	\mD_0(E):=\{F\mid F=\sum_{i=1}^kF_i{e}_i,k\in\mathbb{N}^*,F_i\in\mD_0,{e}_i\in E,i=1,\ldots,k\}.
	\de
	For all $k\in \mN^*$, $p\geq 2$, we denote by $\mD^{n,p}(E)$ the closure of $\mD_0(E)$ w.r.t. the norm
	\ce
	\|F\|_{\mD^{k,p}}&:=&\|F\|_{L^{p}(\mP;E)}+\sum_{i=1}^{k}\|D^{i}F\|_{L^p(\mP;L^2(\hat{\mP}^k;E))},
	\de
    where $D$ is the gradient operator on Poisson space, and $D^i$ is the $i$-th gradient operator.

	In the same way, for all $k\in \mN^*$, $p\geq 2$, we denote by $\overline{\mD}^{\infty}$, the vector subspace of elements $F$ in $\mD^{\infty} \bigcap \cD(A)$ such that $A[F]\in \mD^{\infty}$.  Then define $\overline{\mD}_0(E)$ by
	\ce
		\overline{\mD}_0(E):=\{F\mid F=\sum_{i=1}^kF_i{e}_i,k\in\mathbb{N}^*,F_i\in\overline{\mD}^{\infty},{e}_i\in E,i=1,\ldots,k\}.
	\de
	we denote by $\overline{\mD}^{k,p}(E)$ the closure of $\overline{\mD}_0(E)$ w.r.t. the norm
	\ce
	\|F\|_{\overline{\mD}^{k,p}}&:=&\|F\|_{\mD^{k,p}}+\|A[F]\|_{\mD^{k,p}} .
	\de
$D$ admits a dual operator $\delta$ which is  densely defined, mapping from $ L^2(\mP\times \widehat{\mP})$ to $L^2(\mP)$(see \cite[ Section 5.1]{BL}). $\delta$ is called the divergence operator on the Poisson space. Then $\delta$ is an $ L^2(\mP;E)$ valued operator densely defined on $L^2(\mP;H \otimes E)$, and for any $F\in \mD^{1,2}(E)$, we have
	\ce
	\mathbb{E}[\langle F,\delta Z\rangle_E]=\mathbb{E}[\langle DF,Z\rangle_{H\otimes E}].
	\de
	For any $F_1,F_2\in \mD^{1,2}$, we introduce the following notation
	\ce
	\Gamma[F_1,F_2]:=\langle DF_1,DF_2\rangle_{H}
	\de
	and
	\ce
	\Gamma[F_1]:=\Gamma[F_1,F_1].
	\de

	Let
	\ce
	\mathbb{D}^{k,\infty}(E):=\bigcap_{p\geqslant2}\mathbb{D}^{k,p}(E)\quad\mathrm{and}\quad\mathbb{D}^{\infty}(E):=\bigcap_{k\in\mathbb{N}^{*},p\geqslant2}\mathbb{D}^{k,p}(E),\\
	\overline{\mathbb{D}}^{k,\infty}(E):=\bigcap_{p\geqslant2}\overline{\mathbb{D}}^{k,p}(E)\quad\mathrm{and}\quad\overline{\mathbb{D}}^{\infty}(E):=\bigcap_{k\in\mathbb{N}^{*},p\geqslant2}\overline{\mathbb{D}}^{k,p}(E).
	\de
	For simplification of notations, we denote $\mD^{k,p}$ for $\mD^{k,p}(\mR)$.

	We define the following sets of process:
	\begin{enumerate}[label=\textbullet]
		
		\item $\mathcal{H}_{\mathbb{D}^{k,p}}$ the space of real-valued processes which belong to $L^2([0,T];\mathbb{D}^{k,p})$, i.e.,
\ce
		\|H\|_{\cH_{\mathbb{D}^{k,p}} }:=\mE[\int_0^T \|H(t)\|_{\mD^{k,p}}^2dt]^\frac{1}{2}<\infty.
		\de
		\item $\mathcal{H}_{\overline{\mathbb{D}}^{k,p}}$ the space of real-valued processes which belong to $L^2([0,T];\overline{\mathbb{D}}^{k,p})$, i.e.,
\ce
		\|H\|_{\cH_{\overline{\mathbb{D}}^{k,p}} }:=\mE[\int_0^T \|H(t)\|_{\overline{\mD}^{k,p}}^2dt]^\frac{1}{2}<\infty.
		\de
	\end{enumerate}

	In a natural way, we set
	\ce
	\mathcal{H}_{\mathbb{D}^\infty}:=\bigcap_{k\in\mathbb{N}^*,p\geq 2}\mathcal{H}_{\mathbb{D}^{k,p}},\quad \mathcal{H}_{\overline{\mathbb{D}}^\infty}:=\bigcap_{k\in\mathbb{N}^*,p \geq 2}\mathcal{H}_{\overline{\mathbb{D}}^{k,p}},
	\de
	and
	\ce
	\mathcal{H}_{\mathbb{D}^{k,\infty}}:=\bigcap_{p\geqslant2}\mathcal{H}_{\mathbb{D}^{k,p}},\quad \mathcal{H}_{\overline{\mathbb{D}}^{k,\infty}}:=\bigcap_{p\geqslant2}\mathcal{H}_{\overline{\mathbb{D}}^{k,p}}.
	\de

	We also denote by $\mathcal{H}_{\mathbb{D}^{k,p}}^N$, $\mathcal{H}_{\mathbb{D}^{k,\infty}}^N$ and  $\mathcal{H}_{\mathbb{D}^{\infty}}^N$  the space of $\mathbb{R}^N$-valued processes such that each coordinate belongs, respectively, to   $\mathcal{H}_{\mathbb{D}^{k,p}},\mathcal{H}_{\mathbb{D}^{k,\infty}}$ and $\mathcal{H}_{\mD^{\infty}}$, and we equip them with the standard norms of the product topology. Similarly,  $\mathcal{H}_{\overline{\mathbb{D}}^{k,p}}^N$, $\mathcal{H}_{\overline{\mathbb{D}}^{k,\infty}}^N$ and  $\mathcal{H}_{\overline{\mathbb{D}}^{\infty}}^N$  the space of $\mathbb{R}^N$-valued processes such that each coordinate belongs, respectively, to $\mathcal{H}_{\overline{\mathbb{D}}^{k,p}},\mathcal{H}_{\overline{\mathbb{D}}^{k,\infty}}$ and $\mathcal{H}_{\overline{\mD}^{\infty}}$.

\subsection{ Differential in $\cP_2(\mR^N)$}
	Now we recall the notion of Lions derivatives with respect to measure,
for details see \cite{C, CD, CM}.
	
Let $\Om'=[0,1]$, $\cF'$ the Lebesgue $\s$-algebra and $P'$ the Lebesgue measure. It is well known that any distribution can be realized
on the probability space $(\Om', \cF', P')$.
For a function $U:\cP_2(\mR^N)\rightarrow\mR$, define its lift $U': L^2(\Omega';\mR^N)\mapsto \mR$ by
$$
U'(X):=U([X]).
$$
 If $U'$ is Fr\'echet differentiable at $X$, then we say $U$ is differentiable at $[X]$.
 Then, identifying the dual of $L^2(\Omega';\mR^N)$ as itself, the Fr\'echet derivative
 $DU'(X')\in L^2(\Omega';\mR^N)\mapsto \mR$ is such that
 \be\label{direction}
	DU'(X)(\gamma)=\<DU'(X),\gamma'\>=\mE'[DU'(X)\gamma]\quad \forall \gamma\in L^2(\Omega';\mR^N),
	\ee
 where $\mE'$ is the expectation under $\mP'$.
 By \cite[Theorem 6.2]{C}, in this case there exists a $[X]$-almost surely uniquely
 defined function $\partial_\mu U([X], \cdot): \mR^N\mapsto \mR^N$ such that
 \ce
\partial_{\mu}U([X],X)=DU'(X).
\de
	Furthermore, we will focus on functions where for each $\mu\in\cP_2(\mR^N)$, there exists a unique version of such function $\partial_{\mu}U(\mu,\cdot)$ which is assumed to be a prior continuous as a function (see the discussion in \cite{CM})
	\ce
	\cP_2(\mR^N) \times \mR^N \ni (\mu,v)  \mapsto \partial_{\mu}U(\mu,v).
	\de
	
	To get a more general result, we extend the derivatives to higher order. For a function $f:\cP_2(\mR^N)\rightarrow\mR^N$, we can apply the above discussion straightforwardly to each component $f=(f^1,\cdots,f^N)$. Then the derivatives $\partial_{\mu}f^i$, $1\leq i\leq N$, takes values in $\mR^N$, and we denote $(\partial_{\mu}f^i)_j:\cP_2(\mR^N)\times\mR^N\rightarrow\mR$ for $j=1,\cdots,N$.
	
	For a fixed $v\in\mR^N$, if we are able to  differentiate $\cP_2\ni\mu\mapsto(\partial_{\mu}f^i)_j(\mu,v)\in\mR$ , we can get the second order derivatives. If the derivative of the mapping $\cP_2\ni\mu\mapsto(\partial_{\mu}f^i)_j(\mu,v)$ exists and there is a continuous version of
	\ce
	\cP_2(\mR^N)\times\mR^N\times\mR^N\ni(\mu,v_1,v_2)\mapsto\partial_{\mu}(\partial_{\mu}f^i)_j(\mu,v_1,v_2)\in\mR^N,
	\de
	then it is unique. It is natural to have a multi-index notation $\partial_{\mu}^{(j,k)}f^i:=(\partial_{\mu}(\partial_{\mu}f^i)_j)_k$ to ease the notation.
	
	We can recursively define higher order derivative of $f$ with respect to the measure which is a mapping
	\ce
	\cP_2(\mR^N)\times(\mR^N)^n\ni v_j\mapsto\partial_{\mu}^{\alpha}f^{i_0}(\mu,v_1,\cdots,v_n)\in\mR^N,
	\de
	 if for $\alpha=(i_1,\cdots,i_n)$ and each $(i_0,\cdots,i_n)\in\{1,\cdots,N\}^{n+1}$, the derivative
	\ce
	\partial_{\mu}(\partial_{\mu}\cdots(\partial_{\mu}f^{i_0})_{i_1}\cdots)_{i_n}
	\de
	exists.

	If, for some $j\in\{1,\cdots,N\}$ and all $(\mu,v_1,\cdots,v_{j-1},v_{j+1},\cdots,v_n)\in\cP_2(\mR^N)\times(\mR^N)^{n-1}$,
	\ce
	\mR^N\ni v_j\mapsto\partial_{\mu}^{\alpha}f^{i_0}(\mu,v_1,\cdots,v_n)\in\mR
	\de
	is $l$-times continuously differentiable, we denote the derivatives $\partial_{v_j}^{\beta_j}\partial_{\mu}^{\alpha}f^{i_0}$, for
a multi-index $\beta_j=(x_1,x_2,
\cdots,x_n)$  on $\{1,\cdots,N\}$  with $\#\beta_j:=\dim(\beta_j)\leq l$. Similar to the above, we will denote by $\boldsymbol{\beta}$ the $n$-tuple of multi-indices $(\beta_1,\cdots,\beta_n)$. We also associate a length to $\boldsymbol{\beta}$ by
	\ce
	|\boldsymbol{\beta}|:=|\beta_1|+\cdots+|\beta_n|,
	\de
	and denote $\#\boldsymbol{\beta}:=n$. Then we denote by $\cB_n$ the collection of all such $\mathbf{\beta}$ with $\#\mathbf{\beta}=n$, and $\cB:=\cup_{n\geq1}\cB_n$. Again, to lighten the notation, we use
	\ce \partial_{\boldsymbol{v}}^{\boldsymbol{\beta}}\partial_{\mu}^{\alpha}f^i(\mu,\boldsymbol{v}):=\partial_{v_n}^{\beta_n}\cdots\partial_{v_1}^{\beta_1}\partial_{\mu}^{\alpha}f^i(\mu,v_1,\cdots,v_n).
	\de
	We will need the following result from \cite[Lemma 5.1]{BR}.

	\bl\label{exchange}
	Let ~$g:\mR\times \cP_2(\mR)\to \mR$ and suppose that the derivative functions
	\ce
	(x,\mu,v)\in\mathbb{R}\times\mathcal{P}_2(\mathbb{R})\times\mathbb{R}\rightarrow\left(\partial_x\partial_\mu g(x,\mu,v),\partial_\mu\partial_xg(x,\mu,v)\right)\in\mathbb{R}\times\mathbb{R}
  	\de
	both exist and are Lipschitz continuous, i.e., there exists a constant ~$C > 0$ such that
	\ce
	|\partial_x\partial_\mu g(x,\mu,v),-\partial_x\partial_\mu g(x^{\prime},\mu^{\prime},v^{\prime})|
		\leq C\left(|x-x^{\prime}|+W_2(\mu,\mu^{\prime})+|v-v^{\prime}|\right),\\
    |\partial_\mu\partial_x g(x,\mu,v),-\partial_\mu\partial_x g(x^{\prime},\mu^{\prime},v^{\prime})|
		\leq C\left(|x-x^{\prime}|+W_2(\mu,\mu^{\prime})+|v-v^{\prime}|\right).
	\de
	Then,
	\ce
	\partial_x\partial_\mu g(x,\mu,v)=\partial_\mu \partial_x g(x,\mu,v).
	\de
	\el
	
	Now we introduce the following definition which will be used in the sequel.
	
	\bd
	Let $f:\mR^N\times\Xi\times\cP_2(\mR^N)\rightarrow\mR^N$.
	\begin{enumerate}
		\item
		We say that $f\in\cC_{b,Lip}^{1,1}(\mR^N\times\Xi\times\cP_2(\mR^N);\mR^N)$ if the following is true: $\partial_{\mu}f$ and $\partial_x f$ exist and there exists $H\in L^2(\Xi, \lambda)$ such that
		\ce
		|\partial_xf(x,u,\mu)|+|\partial_{\mu}f(x,u,\mu,{v})|\leq C H(u).
		\de
		 Moreover, suppose that $\partial_{\mu}f$ and $\partial_xf$ are all Lipschitz in the sense that for all $(x,\mu,{v})$, $(x^{\prime},\mu^{\prime},v^{\prime})\in\mR^N\times\cP_2(\mR^N)\times\mR^N$,
		\ce
		|\partial_{\mu}f(x,u,\mu,{v})-\partial_{\mu}f(x',u,\mu',{v}')|&\leq&C H(u)(|x-x'|+|{v}-{v}'|+W_2(\mu,\mu')),\\
		|\partial_xf(x,u,\mu)-\partial_xf(x',u,\mu')|&\leq&C H(u)(|x-x'|+W_2(\mu,\mu')).
		\de
		\item
		We say that $f\in\cC_{b,Lip}^{k,k}( \mR^N\times\Xi\times\cP_2(\mR^N)\rightarrow\mR^N)$ if the following holds true: for all multi-indices $\alpha$ and $\gamma$ on $\{1,\cdots,N\}$ and all $\boldsymbol{\beta}\in\cB$   satisfying $\#\alpha+\#\boldsymbol{\beta}+\# \gamma\leq k$ the derivatives
		\ce
		\partial_x^{\gamma}\partial_{\boldsymbol{v}}^{\boldsymbol{\beta}}\partial_{\mu}^{\alpha}f(x,u,\mu,\boldsymbol{v}),\partial_{\boldsymbol{v}}^{\boldsymbol{\beta}}\partial_{\mu}^{\alpha}\partial_x^{\gamma}f(x,u,\mu,\boldsymbol{v}),\partial_{\boldsymbol{v}}^{\boldsymbol{\beta}}\partial_x^{\gamma}\partial_{\mu}^{\alpha}f(x,u,\mu,\boldsymbol{v})
		\de
		exist and are bounded, Lipschitz continuous.
		\item We say that $f\in \cC_{b,Lip}^{k,k}(\mR^N\times \cP_2(\mR^N);\mR^N)$  if $f$ does not depend on $\Xi$ but otherwise satisfy the conditions in part~$(1)$ and $(2)$.
	\end{enumerate}
	\ed
	
	\br
	\begin{enumerate}
\item If $\gamma=(\gamma_1,\cdots,\gamma_k,\cdots,\gamma_i),\gamma_k\in \{ 1,\cdots,N\}$, then we define the high order derivative of $\partial^\gamma_x f$ as follows
        \ce
        \partial^\gamma_x f= \partial_{x_{\gamma_i}} f(\cdots(\partial_{x_{\gamma_2}} (\partial_{x_{\gamma_1}} f))).
        \de
		\item Due to Lemma \ref{exchange}, we have
		\ce
		\partial_x^{\gamma}\partial_{\boldsymbol{v}}^{\boldsymbol{\beta}}\partial_{\mu}^{\alpha}f(x,u,\mu,\boldsymbol{v})=\partial_{\boldsymbol{v}}^{\boldsymbol{\beta}}\partial_{\mu}^{\alpha}\partial_x^{\gamma}f(x,u,\mu,\boldsymbol{v})=\partial_{\boldsymbol{v}}^{\boldsymbol{\beta}}\partial_x^{\gamma}\partial_{\mu}^{\alpha}f(x,u,\mu,\boldsymbol{v}).
		\de

    \end{enumerate}
	 \er
	\subsection{Assumptions}
	 We will  need the following assumptions.
	\ba\label{R}
	For all multi-indices $\alpha$ and $\gamma$ on $\{1,\cdots,N\}$ and all $\boldsymbol{\beta}\in\cB$   satisfying $\#\alpha+\#\boldsymbol{\beta}+\# \gamma\leq k$,
	\begin{enumerate}
		\item \begin{enumerate}
			\item for $u\in {\Xi}$,   $c(\cdot,u,\cdot)\in \cC_{b,Lip}^{k,k}(\mR^N\times  \cP_2(\mR^N);\mR^N)$, and
			\ce
			\sup_{x\in \mR^N}|\partial_x^\gamma \partial_{\boldsymbol{v}}^{\boldsymbol{\beta}} \partial_{\mu}^{\alpha} c(x,\cdot,\mu,\boldsymbol{v})|\in \bigcap_{p\geq 2}L^p(\Xi,\nu),~~~\forall \mu\in \cP_2(\mR^N), \boldsymbol{v}\in (\mR^N)^{\# \boldsymbol{\beta}};
			\de
			\item $ |c(0,\cdot,\mu)|\in \bigcap_{p\geq 2}L^p(\Xi,\nu), ~~~\forall  \mu\in \cP_2(\mR^N)$;
			\item for $x\in \mR^N$, $\mu\in \cP_2(\mR^N)$, $\boldsymbol{v}\in (\mR^N)^{\# \boldsymbol{\beta}}$, $ \partial_x^\gamma \partial_{\boldsymbol{v}}^{\boldsymbol{\beta}} \partial_{\mu}^{\alpha} c(x,\cdot,\mu,\boldsymbol{v}) \in \bar{\mathbf{d}}^{m,p}$ and
			\ce
			\forall m\in \mN^*, q\geq 2,\quad \sup_{x\in \mR^N} \|\partial_x^\gamma \partial_{\boldsymbol{v}}^{\boldsymbol{\beta}} \partial_{\mu}^{\alpha} c(x,\cdot,\mu,\boldsymbol{v}) \|_{\bar{\mathbf{d}}^{m,q}}\leq C;
			\de
			\item
			$\forall u\in\Xi,$ $\mu\in \cP_2(\mR^N)$, $\boldsymbol{v}\in (\mR^N)^{\# \boldsymbol{\beta}}$
			\ce
			\sup_{x\in \mathbb{R}^N}\left|(I+\partial_x c(x,\cdot,\mu))^{-1}\times c(x,\cdot,\mu)\right|\in \bigcap_{p\geq 2}L^p(\Xi,\nu);
			\de
		\end{enumerate}
		\item for all  $\mu\in \cP_2(\mR^N)$, $\boldsymbol{v}\in (\mR^N)^{\# \boldsymbol{\beta}}$, $b(x,\mu)\in \cC_{b,Lip}^{k,k}(\mR^N\times \cP_2(\mR^N);\mR^N)$
		\ce
		\sup_{x\in\mathbb{R}^N}|\partial_x^\gamma \partial_{\boldsymbol{v}}^{\boldsymbol{\beta}} \partial_{\mu}^{\alpha} b(x,\mu,\boldsymbol{v})|\leqslant C.
		\de
	\end{enumerate}
	\ea

		\ba\label{H_solution}
	Let the coefficients  $b:\mR^N\times \cP_2(\mR^N) \to \mR^N$ and $c:\mR^N\times \Xi\times \cP_2(\mR^N) \to \mR^N$ satisfy:
	\begin{enumerate}
		\item the coefficient $b$ is Lipschitz;
		\item the coefficient $c$ is Borel measurable. Moreover, there exists $H\in L^2(\Xi,\lambda)$ such that, for all $u\in \Xi$, $x,x^{\prime}\in\mR^N$, $\mu,\mu^{\prime} \in \cP_2(\mR^N)$,
		\ce
		|c(x,u,\mu)|&\leq& C H(u),\\
		|c(x,u,\mu)-c(x^{\prime},u,\mu^{\prime}) |&\leq& C H(u)(|x-x^{\prime}|+W_2(\mu,\mu^{\prime})).
		\de
	\end{enumerate}
	\ea
	
	\ba
	$c^{\flat}(x,u,\mu,r)$ satisfies  Lipschitz continuous condition and there exists $H\in L^2(\Xi,\lambda)$ such that, for $x,x^{\prime}\in \mR^N$, $u\in \Xi$, $r\in R$ and $\mu ,\mu^\prime \in \cP_2(\mR^N)$,
	\ce
	|c^{\flat}(x,u,\mu,r)-c^{\flat}(x^{\prime},u,\mu^{\prime},r)|	\leq C H(u)|r|(|x-x^\prime|+W_2(\mu,\mu^\prime)).
	\de
	\ea
	
	\ba[Ellipticity]\label{ellipticity}Assume

(i) Assumption \ref{R} holds,
	 and  the intensity measure $dt\times \nu$ of $N$ is such that $\nu$ has  an infinite mass near some $0$;

(ii) there is an open set $\cO\subset \Xi$ such that $\gamma[c(s,x,\mu)](u)$ is continuous and invertible on
$\mR_+\times \mR^N\times\cP_2(\mR^N)\times \cO$ and
	\ce	\gamma[c(s^{\prime},x,\mu)](u)\geqslant\frac{1}{1+|x|^\delta}\psi(u)I_N,~~~
\forall(s^{\prime},x,u,\mu)\in[0,s]\times\mathbb{R}^N\times\mathcal{O}
\times\cP_2(\mR^N),
\de
where $\geq $ denotes the order relation in the set of symmetric and positives matrices, $\delta$ is a constant, $I_N$ is the identity matrix in $\mR^{N\times N}$, and $\psi$ is an $\mathbb{R}_+\setminus\{0\} $-valued measurable function on $\cO$ such that
	\ce
	\left(\int_0^t\int_{\mathcal{O}}\psi(u)N(ds,du)\right)^{-1}\in\bigcap_{p\geqslant2}L^p(\mathbb{P}).
	\de
	\ea
	
\subsection{ Basic results on MVSDEs driven by L\'{e}vy processes}

	\bp\label{eau}
Under Assumption \ref{H_solution}, Eq.~(\ref{SDE1}) and Eq.~(\ref{SDE}) admit unique strong solutions $X^\theta$ and
$X^{x,\theta}$ respectively, and for $p\geq 2$, there exists a constant $C$ which depends on the Lipschitz constants of  $b$ and $c$, such that for all \ce
	(x,\theta,t),(x^{\prime},\theta^{\prime},t^{\prime}) \in \mR^N\times L^p(\cF_0;\mR^N)\times [0,T],
	\de
	we have
	\be\label{estX}
	\|X^{\theta}\|_{S_t^p}\leq C(1+\|\theta\|_p),
	\ee
	\be\label{estX2}
	\|X^{x,\theta}\|_{S_t^p}\leq C(1+|x|+\|\theta\|_2 )
	\ee
	and
	\be\label{est3}
	 (\mE[ \sup_{s\in[0,h]}(|X_{s}^{\theta}-\theta|^p+ |X_{s}^{x,[\theta]}-x|^p)])^{\frac{1}{p}} \leq C h  (1+\|\theta\|_p).
	\ee
	Moreover,
	\be\label{estX3}
	\left\|X^{x,[\theta]}-X^{x^{\prime},[\theta^{\prime}]}\right\|_{\mathcal{S}_t^p}\leq C\left(|x-x^{\prime}|+\|\theta-\theta^{\prime}\|_2\right).
	\ee
	Finally, we have the flow property, i.e.,
	\ce
	\big(  X_{t,r}^{X_{s,t}^{x,[\theta]},[X_{s,t}^{\theta}]} , X_{t,r}^{X_{s,t}^{\theta}}\big)=( X_{s,r}^{x,[\theta]}, X_{s,r}^\theta),\quad r\in [t,T]
	\de
	for any $0\leq s \leq t \leq T$, $x\in \mR^N$ and $\theta \in L^2(\cF_s;\mR^N)$.

	\ep
	\begin{proof}
	The existence and uniqueness of the solutions and the flow property are proven in \cite[Theorem 3.1]{HL}.
	
For the estimates, we prove only (\ref{estX3}), as the others can be obtained analogously. Using Kunita's second inequality,   H\"{o}lder's inequality and the Lipschitz continuous on coefficients, we have
    \ce
    &&\left\|X^{x,[\theta]}-X^{x^{\prime},[\theta^{\prime}]}\right\|_{\mathcal{S}_t^p}^p\\
    &=& \mE[ \sup_{s\in [0,t]} |(x-x^{\prime})+\int_0^s
    b(X_v^{x,[\theta]},[X_v^\theta])-b(X_v^{x^{\prime},[\theta^{\prime}]},[X_v^{\theta^{\prime}}])dv\\
    &&+ \int_0^s \int_{\Xi} c(X_{v-}^{x,[\theta]},u,[X_{v-}^{\theta}]) -c(X_{v-}^{x^{\prime},[\theta^{\prime}]},u,[X_{v-}^{\theta^{\prime}}])\widetilde{N}(dv,du) |^p]\\
    &\leq & C(|x-x^{\prime}|^p +\mE[\int_0^t |b(X_v^{x,[\theta]},[X_v^\theta])-b(X_v^{x^{\prime},[\theta^{\prime}]},[X_v^{\theta^{\prime}}])|^pdv]\\
    &&+\mE[\int_0^t \big|\int_{\Xi} |c(X_{v-}^{x,[\theta]},u,[X_{v-}^{\theta}]) -c(X_{v-}^{x^{\prime},[\theta^{\prime}]},u,[X_{v-}^{\theta^{\prime}}])|^2 \lambda(du)\big|^{\frac p 2}dv ]\\
    &&+ \mE[(\int_0^t \int_{\Xi} |c(X_{v-}^{x,[\theta]},u,[X_{v-}^{\theta}]) -c(X_{v-}^{x^{\prime},[\theta^{\prime}]},u,[X_{v-}^{\theta^{\prime}}])|^p \lambda(du)dv )] )\\
    &\leq & C(|x-x^{\prime}|^p +\mE[\int_0^t \big||X_v^{x,[\theta]}-X_v^{x^{\prime},[\theta^{\prime}]}|+W_2(X_v^\theta,X_v^{\theta^{\prime}})\big|^pdv]\\
    &&+\mE[\int_0^t \big|\int_{\Xi} CH(u)^2\big(|X_{v-}^{x,[\theta]}-X_{v-}^{x^{\prime},[\theta^{\prime}]}|+W_2(X_{v-}^{\theta},X_{v-}^{\theta^{\prime}} )\big)^2 \lambda(du)\big|^{\frac p 2}dv ]\\
    &&+ \mE[\int_0^t \int_{\Xi}  CH(u)^2\big(|X_{v-}^{x,[\theta]}-X_{v-}^{x^{\prime},[\theta^{\prime}]}|^p+W_2(X_{v-}^{\theta},X_{v-}^{\theta^{\prime}} )^p\big)
    \lambda(du)dv] )\\
    &\leq &C(|x-x^{\prime}|^p  +\int_0^t W_2(X_v^{\theta},X_v^{\theta^{\prime}})^p dv+\mE[\int_0^t \sup_{s\in [0,v]} |X_s^{x,[\theta]}-X_s^{x^{\prime},[\theta^{\prime}]}|^p dv]).
    \de
    By Gr\"onwall's inequality,
    \ce
    \left\|X^{x,[\theta]}-X^{x^{\prime},[\theta^{\prime}]}\right\|_{\mathcal{S}_t^p}^p
    &\leq &C(|x-x^{\prime}|^p  +\int_0^t W_2(X_v^{\theta},X_v^{\theta^{\prime}})^p dv).
    \de
Since
	\ce
	\|X_v^{\theta}-X_v^{\theta^{\prime}} \|_2^2
	 &\leq& C(\mE[| \theta-\theta^{\prime}|^2]+\mE[\int_0^v |b(X_r^{\theta},[X_r^\theta])-b(X_r^{\theta^{\prime}},[X_r^{\theta^{\prime}}])|^2dr]\\
    &&+ \mE[(\int_0^v \int_{\Xi} |c(X_{r-}^{\theta},u,[X_{r-}^{\theta}]) -c(X_{r-}^{\theta^{\prime}},u,[X_{r-}^{\theta^{\prime}}])|^2 \lambda(du)dr )] )\\
	&\leq&  C (\mE[| \theta-\theta^{\prime}|^2]+  \int_0^v \mE[ |X_r^{\theta}-X_r^{\theta^{\prime}}|^2]dr+\int_0^v W_2(X_v^\theta, X_v^{\theta^{\prime}})^2dr)\\
	&\leq& C (\mE[| \theta-\theta^{\prime}|^2]+  \int_0^v \mE[ |X_r^{\theta}-X_r^{\theta^{\prime}}|^2]dr),
	\de
	by Gr\"ownall's inequality again,
	\ce
	\|X_v^\theta-X_v^{\theta^{\prime}} \|_2 &\leq& C\| \theta-\theta^{\prime}\|_2.
	\de
But $W_2(X_v^{\theta},X_v^{\theta^{\prime}})^2
	\leq\|X_v^{\theta}-X_v^{\theta^{\prime}} \|_2^2$,	
	hence
	\ce
	\left\|X^{x,[\theta]}-X^{x^{\prime},[\theta^{\prime}]}\right\|_{\mathcal{S}_t^p}^p&\leq& C(| x-x^{\prime}|^p + \| \theta-\theta^{\prime}\|_2^p )\leq C(| x-x^{\prime}| + \| \theta-\theta^{\prime}\|_2 )^p.
	\de

	\end{proof}
	
	\section{Regularities of Solutions}
	This section studies the differentiability of the solution $X_t^{x,[\theta]}$, investigating both its classical differentiability with respect to $(x,\mu)$ and its Malliavin differentiability. These regularity results are a prerequisite for the integration by parts formulas derived in Section 4.
	
	\subsection{First Order Derivatives}
	\bt\label{dfo}
	Suppose that $b \in\cC_{b,Lip}^{1,1}(\mR^N\times\cP_2(\mR^N);\mR^N)$ and  $c\in\cC_{b,Lip}^{1,1}(\mR^N\times\Xi \times\cP_2(\mR^N);\mR^N)$. Then the following statements hold:
	\begin{enumerate}
		\item
		For all ~$t\in[0,T]$,  the map $\mR^N \ni x \mapsto X_t^{x,[\theta]}\in L^2(\cF_t;\mR^N)$ is $L^2$-differentiable.
More precisely, there exists a unique process $\partial_x X_t^{x,[\theta]}$ such that
		\ce
		\mE[\sup_{s\in[0,t]} |X_s^{x+h,[\theta]}-X_s^{x,[\theta]}-\partial_x X_s^{x,[\theta]}h|]=o(|h^2|)\quad \mR^N \ni h \to 0.
		\de
		 Moreover, $\partial_x X_t^{x,\theta}$ is the unique solution of the following SDE with jumps:
		\be\label{p_x}
		\partial_xX_t^{x,[\theta]}&=&I+\int_0^t \partial b(X_{s}^{x,[\theta]},[X_{s}^{\theta}]) \partial_x X_{s}^{x,[\theta]} ds\no\\
		&&+\int_0^t \int_{\Xi}\partial c(X_{s-}^{x,[\theta]},u,[X_{s-}^{\theta}])\partial_xX_{s-}^{x,[\theta]}\widetilde{N}(ds,du),
		\ee
		where $\partial b$ and $\partial c $ are the derivatives with respect to the first variable respectively.
		\item
		For all $t\in[0,T]$, the mapping  $\theta\mapsto X_t^{x,[\theta]}:L^2(\cF_0;\mR^N)\to L^2(\cF_t;\mR^N)$ is Fr\'{e}chet differentiable.
	Moreover, for all $(t,x,\theta)\in [0,T]\times \mR^N \times L^2(\cF_0;\mR^N)$, and $p\geq1$, the map $\cP_2(\mR^N)\ni[\theta]\mapsto X_t^{x,[\theta]}\in L^p(\Omega)$ is differentiable. Hence, $\partial_{\mu}X_t^{x,[\theta]}(v)$ exists and it satisfies the following SDE with jumps:
		\be\label{p_mu}
		\partial_{\mu}X_t^{x,[\theta]}(\boldsymbol{v})&=& \int_{0}^{t}\left\{\partial b\left(X_{s}^{x,[\theta]},\left[X_{s}^{\theta}\right]\right) \partial_{\mu}X_s^{x,[\theta]}(\boldsymbol{v})\right. \no\\
		&&+{\mathbb{E}^{\prime}}\left[\partial_{\mu} b\left(X_{s}^{x,[\theta]},\left[X_{s}^{\theta}\right], ({X}_{s}^{{v},[\theta]})^{\prime}\right) \partial_{x} ({{X}_{s}^{{v},[\theta]}})^{\prime}\right] \no\\
		&&+\left.\mathbb{E}^{\prime}\left[\partial_{\mu} b\left(X_{s}^{x,[\theta]},\left[X_{s}^{\theta}\right], ({X}_{s}^{{\theta}^{\prime}})^{\prime}\right) (\partial_{\mu}X_s^{\theta^{\prime},[\theta]}(\boldsymbol{v}))^{\prime}\right]\right\} ds\no\\
		&&+\int_{0}^{t}\int_{\Xi}\left\{\partial c\left(X_{s-}^{x,[\theta]},u,\left[X_{s-}^{\theta}\right]\right) \partial_{\mu}X_{s-}^{x,[\theta]}(\boldsymbol{v})\right. \no\\
		&&+{\mathbb{E}^{\prime}}\left[\partial_{\mu} c\left(X_{s-}^{x,[\theta]},u,\left[X_{s-}^{\theta}\right], ({X}_{s-}^{{v},[\theta]})^{\prime}\right) \partial_{x} ({{X}_{s-}^{{v},[\theta]}})^{\prime}\right] \no\\
		&&+\left.\mathbb{E}^{\prime}\left[\partial_{\mu} c\left(X_{s-}^{x,[\theta]},u,\left[X_{s-}^{\theta}\right], ({X}_{s-}^{{\theta}^{\prime}})^{\prime}\right) (\partial_{\mu}X_{s-}^{\theta^{\prime},[\theta]}(\boldsymbol{v}))^{\prime}\right]\right\} \widetilde{N}(ds,du),
		\ee
		where $({X}_s^{{\theta}^{\prime}})^{\prime}$ is a copy of $X_s^{\theta}$ on a different  probability space $({\Omega}^{\prime},{\sF}^{\prime},{\mP}^{\prime})$, $\partial_x({X}_s^{v,[\theta]})^{\prime}$ is a copy of $\partial_x X_s^{v,[\theta]}$,   $\partial_{\mu}({X}_s^{{\theta}^{\prime},[\theta]})^{\prime}(\boldsymbol{v})=\partial_{\mu}({X}_s^{x,[\theta]})^{\prime}(\boldsymbol{v})|_{x={\theta}^{\prime}}$ with $\theta^{\prime}\in L^2(\cF_0;\mR^N)$ and $[\theta^{\prime}]=[\theta]$. Moreover, the Fr\'{e}chet derivative $\cD X_t^{x,[\theta]}$ holds for all $\gamma, {\gamma}^{\prime}\in L^2(\cF_0;\mR^N)$:
		\be\label{cD}
		\cD X_t^{x,[\theta]}(\gamma)={\mE}^{\prime}[\partial_\mu X_t^{x,[\theta]}(\theta^{\prime}){\gamma}^{\prime}].
		\ee
		\item Under Assumption \ref{R},
		for all $t\in[0,T]$,  $X_t^{x,[\theta]},X_t^{\theta}\in \cH_{{\mD}^{1,\infty}}^N$.
		Moreover,
		\be\label{Md}
		DX_t^{\theta}&=&\int_0^t \partial b(X_s^{\theta},[X_s^\theta]) DX_s^{\theta} ds+\int_0^t \int_{\Xi} \partial c( X_{s-}^{\theta},u,[ X_{s-}^{\theta}]) DX_{s-}^{\theta} \widetilde{ N}(ds,du)\no\\
		&&+\int_0^t \int_{\Xi\times R } c^{\flat}( X_{s-}^{\theta},u,[ X_{s-}^{\theta}],r) N\odot \rho (ds,du,dr).
		\ee

	\end{enumerate}
	\et
	
	\begin{proof}
	
(1) The  $ X^{x,[\theta]}$ satisfies a classical  time-inhomogeneous SDE with jumps, it follows from \cite[Lemma  6.5.1]{ HK} that the first derivative satisfies Eq. (\ref{p_x}).

(2) It is shown that the map $L^2(\cF_0)\ni\theta \mapsto  X_t^{x,[\theta]}$ is Fr\'echet  differentiable in \cite[Theorem 4.3]{HL}, so $\partial_{\mu} X_t^{x,[\theta]}(v)$ exists and satisfies (\ref{p_mu}). Finally, we can use the idea in \cite{HL} to prove (\ref{cD}). The calculations are  similar, hence its proof is omitted here.

(3) Let~$ X^{\theta,n}$ denote the Picard approximation of the solution to the MVSDE (\ref{SDE1}), given by
		\ce\label{ab}
		X_t^{\theta,0}&=&\theta, \quad t\in [0,T],\\
		X_t^{\theta,n}&=&\theta+\int_0^t b(X_s^{\theta,n},[X_s^{\theta,n-1}])ds+\int_0^t \int_{\Xi} c( X_{s-}^{\theta,n},u,[ X_{s-}^{\theta,n-1}]) \widetilde{ N}(ds,du).
		\de	
		For each $n\geq 1$, $ X^{\theta,n}$ solves a  classical SDE with jumps.
		
	 By  Kunita's first inequality, H\"older's inequality and the assumption on  coefficients, we have
		\be
		&&\mE[\sup_{s\in[0,t]}|X_{s}^{\theta}-X_{s}^{\theta,n}|^{2}]\no\\
		&\leq& C \mE[\sup_{s\in[0,t]}|\int_0^s [b(X_r^\theta,[X_r^\theta])-b(X_r^{\theta,n},[X_r^{\theta,n-1}])]dr|^{2}]\no\\
		&&+ C \mE[\sup_{s\in[0,t]}| \int_0^s \int_{\Xi}[c(X_{r-}^\theta,u,[X_{r-}^\theta] )-c(X_{r-}^{\theta,n},u,[X_{r-}^{\theta,n-1}])] \widetilde{N}(dr,du)|^2 ]\no\\
		&\leq& C \mE[\int_0^t [b(X_r^\theta,[X_s^\theta])-b(X_r^{\theta,n},[X_r^{\theta,n-1}])]^2 dr]\no\\
		&&+ C \mE[ \int_0^t \int_{\Xi}[c(X_{r-}^\theta,u,[X_{r-}^\theta] )-c(X_{r-}^{\theta,n},u,[X_{r-}^{\theta,n-1}])]^2 \lambda(du)dr ]\no\\
		&\leq& C \mE[\int_0^t |X_r^\theta -X_r^{\theta,n}|^2 +W_2^2( [X_r^\theta], [X_r^{\theta,n-1}] ) dr]\no\\
		&&+ C \mE[ \int_0^t \int_{\Xi}H(u)^2(|X_{r}^\theta -X_{r}^{\theta,n}|^2 +W_2^2( [X_{r}^\theta], [X_{r}^{\theta,n-1}] ) ) \lambda(du)dr]\no\\
		&\leq & C \int_{0}^{t} \mathbb{E}\left[\left|X_{r}^{\theta}-X_{r}^{\theta, n}\right|^{2}+\mathbb{E}\left|X_{r}^{\theta}-X_{r}^{\theta, n-1}\right|^{2}\right] d r.
		\ee
		Set
		\ce
		\phi_n(t)=	\mE[\sup_{s\in[0,t]}|X_{s}^{\theta}-X_{s}^{\theta,n}|^{2}].
		\de
Then the inequality (\ref{ab}) can be written as
		\ce
		\phi_n(t) \leq C \int_0^t ( \phi_n(r)+\phi_{n-1}(r) ) dr.
		\de
	    We now prove the following inequality by  induction
		\ce
		\phi_n(t)\leq C_0 \frac{(Ct e^{Ct})^n}{n!}.
		\de
		When $n=0$,
		\ce
		\phi_0(t)&=&\mE[ \sup_{s\in[0,t]} | X_s^{\theta}-\theta|^2]\leq C(1+\int_0^t\mE[|X_r^\theta|^2] dr)
		\leq C(1+t\|\theta\|_2):=C_0.
		\de
		 Suppose that $\phi_{n-1}(t)\leq C_0 \frac{(Ct e^{Ct})^{n-1}}{(n-1)!}$.  Then by  induction
		\ce
		\phi_n(t) &\leq& C \int_0^t ( \phi_n(r)+\phi_{n-1}(r) ) dr\\
		&\leq& C \int_0^t ( \phi_n(r)+C_0 \frac{(Cr e^{Cr})^{n-1}}{(n-1)!} ) dr\\
		&\leq& C \int_0^t  \phi_n(r)dr+ C_0 \frac{(Ct)^{n} e^{Ct(n-1)}}{(n)!},
		\de
which gives in virtue of Gr\"onwall inequality,
		\ce
		\phi_n(t)\leq C_0 \frac{(Ct)^{n} e^{Ct(n-1)}}{(n)!} exp( \int_0^t C dt)\leq C_0 \frac{(Ct e^{Ct})^n}{n!}.
		\de
Hence
$$\lim_{n\to \infty}\phi_n(t)\to 0.
$$
Thus
		 \be\label{8.7}
		 \lim_{n\to+\infty}\mE[\sup_{s\in[0,t]}|X_{s}^{\theta}-X_{s}^{\theta,n}|^{2}]=0.
		 \ee
		As a consequence of \cite[Proposition 8.19 and 8.23 and Lemma 8.24]{BL}, it is clear that for all ~$n,X^{\theta,n}$ belongs
to ~$ \cH_{\mD^{1,p}}^N$ and we have ~$\forall t\in [0,T]$
		\ce
		DX_t^{\theta,n}&=&\int_0^t \partial b(X_s^{\theta,n}, [X_s^{\theta,n-1}])  D X_s^{\theta,n}ds +\int_0^t \int_{\Xi} \partial c( X_{s-}^{\theta,n},u,[ X_{s-}^{\theta,n-1}]) D X_{s-}^{\theta,n} \widetilde{ N}(ds,du)\\
		&&+\int_0^t \int_{\Xi\times R} c^{\flat}( X_{s-}^{\theta,n},u,[ X_{s-}^{\theta,n-1}],r) N\odot \rho (ds,du,dr).
		\de
		Then thanks to  Kunita's second inequality and the  hypotheses we made on coefficients, we have
		\ce
		&&\mathbb{E}\widehat{\mathbb{E}}[|DX_t^{\theta,n}|^p]\\
 &\leq& C( \mathbb{E}\widehat{\mathbb{E}}[|\int_0^t \partial b(X_s^{\theta,n}, [X_s^{\theta,n-1}])  D X_s^{\theta,n}ds|^p ]\\
		&&+ \mathbb{E}\widehat{\mathbb{E}}[|\int_0^t \int_{\Xi} \partial c( X_{s-}^{\theta,n},u,[ X_{s-}^{\theta,n-1}]) D X_{s-}^{\theta,n} \widetilde{ N}(ds,du)|^p]\\
		&&+\mathbb{E}\widehat{\mathbb{E}}[|\int_0^t \int_{\Xi\times R} c^{\flat}( X_{s-}^{\theta,n},u,[ X_{s-}^{\theta,n-1}],r) N\odot \rho (ds,du,dr)|^p ] )\\
		& \leq &  \mathbb{E}\widehat{\mathbb{E}}[|\int_0^t \partial b(X_s^{\theta,n}, [X_s^{\theta,n-1}])  D X_s^{\theta,n}ds|^p ]\\
		&&+\mathbb{E}\widehat{\mathbb{E}}[\int_0^t \int_{\Xi} |\partial c( X_{s-}^{\theta,n},u,[ X_{s-}^{\theta,n-1}]) D X_{s-}^{\theta,n}|^p ]\lambda(du)ds\\
		&&+ \mathbb{E}\widehat{\mathbb{E}}[\int_0^t | \int_{\Xi} (\partial c( X_{s-}^{\theta,n},u,[ X_{s-}^{\theta,n-1}]) D X_{s-}^{\theta,n})^2 \lambda(du)|^{\frac p 2}ds]\\
		&&+ C \mathbb{E}[|\int_0^t \int_{\Xi\times R} |c^{\flat}( X_{s-}^{\theta,n},u,[ X_{s-}^{\theta,n-1}],r)|^p \lambda(du) \rho (dr)|ds ]\\
		&&+ C \mathbb{E}[|\int_0^t |\int_{\Xi\times R} |c^{\flat}( X_{s-}^{\theta,n},u,[ X_{s-}^{\theta,n-1}],r)|^2 \lambda(du) \rho (dr)|^{\frac p 2} ds ] \\
		& \leq &  C\mathbb{E}\widehat{\mathbb{E}}[|\int_0^t ( D X_s^{\theta,n})^p ds|  ]\\
		&&+C \mathbb{E}[|\int_0^t \int_{\Xi\times R} H^2(u) r^p(1+|X_{s-}^{\theta,n}|^p+\mE[|X_{s-}^{\theta,n-1}|^p]) \lambda(du) \rho (dr)| ds ] \\
		&&+ C \mathbb{E}[\int_0^t |\int_{\Xi\times R} H^2(u) r^2(1+|X_{s-}^{\theta,n}|^2+\mE[|X_{s-}^{\theta,n-1}|^2] ) \lambda(du) \rho (dr)|^{\frac p 2} ds ] \\
		&\leq& C(1+|x|+\|\theta\|_2)^p+C\int_0^t  \mathbb{E}\widehat{\mathbb{E}}[( D X_s^{\theta,n})^p ]ds).
		\de
		By Gr\"onwall's ineqaulity,  for all ~$n\in \mN^{*}$ and all $t\in [0,T]$,
		\ce
		\forall n\in\mathbb{N},\forall t\in[0,T],\mathbb{E}\widehat{\mathbb{E}}\left[|DX_t^{\theta,n}|^p\right]\leqslant C e^{C t}.
		\de
		Hence the sequence ~$(X^{\theta,n})$ is bounded in ~$\cH_{\mD^{1,p}}^N$ which is a Banach space.  So, there is a sequence of convex combinations of ~$X^{\theta,n}$ which converges to a process ~$Y\in \cH_{\mD^{1,p}}^N$. But by (\ref{8.7}) we know that ~$X^{\theta,n}$ tends to ~$X^{\theta}$ in ~$ L^2([0,T];\mR^N)$, hence $Y=X^{\theta} \in \cH_{\mD^{1,p}}^{N}$.
		
		Moreover, by \cite[Proposition 8.19 and 8.22]{BL}, we know that $DX_t^{\theta}$ satisfies
		\ce
		DX_t^{\theta}&=&\int_0^t \partial b(X_s^{\theta},[X_s^\theta]) DX_s^{\theta} ds+\int_0^t \int_{\Xi} \partial c( X_{s-}^{\theta},u,[ X_{s-}^{\theta}])  DX_{s-}^{\theta} \widetilde{ N}(ds,du)\no\\
		&&+\int_0^t \int_{\Xi\times R} c^{\flat}( X_{s-}^{\theta},u,[ X_{s-}^{\theta}],r) N\odot \rho (ds,du,dr).
		\de
Consequently,
\ce
&&\mE \widehat{\mE}[ \sup_{s\in[0,t]}|DX_s^{\theta}-DX_s^{\theta,n}|^p]\\
		&\leq&\mE \widehat{\mE}[ |\int_0^t (\partial b(X_v^{\theta},[X_v^\theta]) DX_v^{\theta}-\partial b(X_v^{\theta,n},[X_v^{\theta,n-1}]) DX_v^{\theta,n}) dv|^p] \\
		&&+\mE \widehat{\mE}[ \sup_{s\in[0,t]}|\int_0^s \int_{\Xi} (\partial c( X_{v-}^{\theta},u,[ X_{v-}^{\theta}])  DX_{v-}^{\theta}-\partial c( X_{v-}^{\theta,n},u,[ X_{v-}^{\theta,n-1}])  DX_{v-}^{\theta,n}) \widetilde{ N}(dv,du)|^p]\\
		&&+ \mE \widehat{\mE}[ \sup_{s\in[0,t]}|\int_0^s \int_{\Xi\times R} [c^{\flat}( X_{v-}^{\theta},u,[ X_{v-}^{\theta}],r)-c^{\flat}( X_{v-}^{\theta,n},u,[ X_{v-}^{\theta,n-1}],r)] N\odot \rho (dv,du,dr)|^p].
\de
Here we calculate only the second term; the others are the same. For all $p\geq 2$, by Kunita's second inequality, H\"older's inequality, the assumptions on coefficients, the boundedness of $\mE \widehat{\mE}[D X_s^{\theta,n}]$ and $W_2([X_{s-}^{\theta}],[X_{s-}^{\theta,n}])\leq (\mE[|X_{s-}^{\theta}-X_{s-}^{\theta,n}|^2])^{\frac 1 2} $, and noticing that the integrand can be written as
\ce
&&\partial c( X_{v-}^{\theta}, u, [ X_{v-}^{\theta}]) \, DX_{v-}^{\theta} - \partial c( X_{v-}^{\theta,n}, u, [ X_{v-}^{\theta,n-1}]) \, DX_{v-}^{\theta,n} \\
&=& \Bigl( \partial c( X_{v-}^{\theta}, u, [ X_{v-}^{\theta}]) - \partial c( X_{v-}^{\theta,n}, u, [ X_{v-}^{\theta}]) \Bigr) \, DX_{v-}^{\theta} \\
&&+ \Bigl( \partial c( X_{v-}^{\theta,n}, u, [ X_{v-}^{\theta}]) - \partial c( X_{v-}^{\theta,n}, u, [ X_{v-}^{\theta,n-1}]) \Bigr) \, DX_{v-}^{\theta} \\
&&+ \partial c( X_{v-}^{\theta,n}, u, [ X_{v-}^{\theta,n-1}]) \, \Bigl( DX_{v-}^{\theta} - DX_{v-}^{\theta,n} \Bigr),
\de
 we have
\ce
&&\mE \widehat{\mE}[ \sup_{s\in[0,t]}|\int_0^s \int_{\Xi} (\partial c( X_{v-}^{\theta},u,[ X_{v-}^{\theta}])  DX_{v-}^{\theta}-\partial c( X_{v-}^{\theta,n},u,[ X_{v-}^{\theta,n-1}])  DX_{v-}^{\theta,n}) \widetilde{ N}(dv,du)|^p]\\
&\leq & C(\|X_s^{\theta}-X_s^{\theta,n}\|_2^p+\|X_s^{\theta,n}-X_s^{\theta,n-1}\|_2^p +\mE \widehat{\mE}[\sup_{s\in[0,t]}|\int_0^s |DX_v^{\theta}-DX_v^{\theta,n} |^p  dv  ]).
\de
		Then it follows that
		\ce
		\mE \widehat{\mE}[ \sup_{s\in[0,t]}|DX_s^{\theta}-DX_s^{\theta,n}|^p]&\leq& C(\|X_s^{\theta}-X_s^{\theta,n}\|_2^p+\|X_s^{\theta,n}-X_s^{\theta,n-1}\|_2^p\\
    &&+\mE \widehat{\mE}[\sup_{s\in[0,t]}\int_0^s |DX_v^{\theta}-DX_v^{\theta,n} |^p  dv  ]).
		\de
So by Gr\"onwall's inequality and (\ref{8.7})
		\ce
		\lim_{n\to +\infty}\mE \widehat{\mE}[\sup_{t\in [0,T]} |DX_t^{\theta}-DX_t^{\theta,n}|^p]=0.
		\de
		Therefore, we have proved that for all $p\geq 2$, $X^{\theta,n}$ converges to $X^{\theta}$ in $\cH_{\mD^{1,p}}^N$.
 		
		The measure term in the coefficients of the equation for $X_t^{x,[\theta]}$ being deterministic,  $DX_t^{x,[\theta]}$ satisfies the standard equation  which can be straightly derived from \cite[Proposition 8.26]{BL}.
	\end{proof}
	
	\subsection{Higher Order Derivatives}
	In order to derive the higher order differentiability and a derivative formula of the density function,
some special classes of processes are used in \cite{K,KS} in the context of the classic Malliavin calculus on the Wiener space.
Here we need an analogous notion on the Poisson space, which will play an important role in the proof of integration by parts formulas.
	
	\bd
	Let $E$ be a separable Hilbert space and let $r\in\mR$, $q,M\in\mN$. We denote by $\mK_r^q(E,M)$ the set of adapted processes $\Psi:[0,T]\times\mR^N\times\Xi \times\cP_2(\mR^N) \rightarrow \cH_{{\mD}^{M,\infty}}(E)$ satisfying the following:
	\begin{enumerate}
		\item
		For any multi-indices $\alpha,\boldsymbol{\beta},\gamma$ satisfying $\#\alpha+\#\boldsymbol{\beta}+\# \gamma\leq M$, the function 
\ce
\partial_x^{\gamma}\partial_{\boldsymbol{v}}^{\boldsymbol{\beta}}\partial_{\mu}^{\alpha}\Psi(t,x,u,[\theta],\boldsymbol{v}) 
\de
exists, the mapping
		\ce
		\mR^N\times\cP_2(\mR^N)\times (\mR^N)^{\# \boldsymbol{\beta}}  \ni(x,[\theta],\boldsymbol{v})\mapsto\partial_x^{\gamma}\partial_{\boldsymbol{v}}^{\boldsymbol{\beta}}\partial_{\mu}^{\alpha}\Psi(t,x,u,[\theta],\boldsymbol{v})\in L^p(\Omega)
		\de
		  is continuous  and the mapping
        \ce
        [0,T]\ni t\mapsto\partial_x^{\gamma}\partial_{\boldsymbol{v}}^{\boldsymbol{\beta}}\partial_{\mu}^{\alpha}\Psi(t,x,u,[\theta],\boldsymbol{v})\in L^p(\Omega)
        \de
        is c\`{a}dl\`{a}g  for all $p\geq 1$.
		\item
		For any $p\geq 1$ and $m\in\mN$ with $\#\alpha+\#\boldsymbol{\beta}+\# \gamma+m\leq M$, we have
		\ce
		\sup_{\boldsymbol{v}\in(\mR^N)^{\#\boldsymbol{\beta}}}\sup_{t\in[0,T]}t^{-r/2}\| \partial_x^{\gamma}\partial_{\boldsymbol{v}}^{\boldsymbol{\beta}}\partial_{\mu}^{\alpha}\Psi(t,x,u,[\theta],\boldsymbol{v})\|_{{\mD}^{m,p}(E)}\leq C(1+|x|+\|\theta\|_2)^q.
		\de
	\end{enumerate}
	\ed
	The higher derivatives of $X_t^{x,[\theta]}$ is summarized in the following  theorem.
	\bt\label{KS}
	Suppose $ b \in\cC_{b,Lip}^{k,k}(\mR^N\times\cP_2(\mR^N);\mR^N)$ and  $c\in\cC_{b,Lip}^{k,k}(\mR^N\times \Xi \times\cP_2(\mR^N);\mR^N)$. Then $(t,x,[\theta])\mapsto X_t^{x,[\theta]}\in\mK_0^1(\mR^N,k)$. If, in addition, $b$ and $c$ are uniformly bounded then $(t,x,[\theta])\mapsto X_t^{x,[\theta]}\in\mK_0^0(\mR^N,k)$.
	\et
	The proof is quite mechanical and  given in the Appendix \ref{appA}.

	\section{Integration by Parts Formulas}
		
We now introduce some operators acting on $\mK_r^q(E,M)$. These operators will appear later in the integration by parts formulas.
For $\Psi\in \mK_r^q(\mR;n),$ $\gamma=(i)$, set
	\ce
	Z_{(i)}^1(\Psi)(t,x,[\theta])&:=&\delta\big( [(DX_t^{x,[\theta]})^{*} (\Gamma[X_t^{x,[\theta]}])^{-1}  \partial_x X_t^{x,[\theta]}]_i  \Psi(t,x,[\theta])  \big),\\
	 Z_{(i)}^2(\Psi)(t,x,[\theta])&:=&\delta\big( [(DX_t^{x,[\theta]})^{*} (\Gamma[X_t^{x,[\theta]}])^{-1} ]_i  \Psi(t,x,[\theta])  \big),\\
	Z_{\mu,(i)}^1(\Psi)(t,x,[\theta],\boldsymbol{v})&:=&\delta\big( [(DX_t^{x,[\theta]})^{*} (\Gamma[X_t^{x,[\theta]}])^{-1}  \partial_{\mu} X_t^{x,[\theta]}]_i  \Psi(t,x,[\theta]) \big),
	\de
where $(DX_t^{x,[\theta]})^{*}$ denotes the transpose of the matrix ~$DX_t^{x,[\theta]}$, and $[\cdot]_i,i\in \{1,\cdots,N\}$ denotes
the $i$-th column of the matrix.
	
 For multi-indexes $\gamma=(\gamma_1,\cdots,\gamma_k,\cdots,\gamma_i)$ and $\alpha=(\alpha_1,\cdots,\alpha_k,\cdots,\alpha_{j})$, $\gamma_k,\alpha_k\in \{ 1,\cdots,N\}$, we define by induction
    \ce
    Z_{\gamma}^1(\Psi)(t,x,[\theta])&:=& Z_{(\gamma_i)}^1\circ\cdots \circ Z_{(\gamma_2)}^1 \circ Z_{(\gamma_1)}^1(\Psi)(t,x,[\theta]),\\
    Z_{\gamma}^2(\Psi)(t,x,[\theta])&:=& Z_{(\gamma_i)}^2\circ\cdots \circ Z_{(\gamma_2)}^2 \circ Z_{(\gamma_1)}^2(\Psi)(t,x,[\theta]),\\
    Z_{\gamma}^3(\Psi)(t,x,[\theta])&:=&Z_{\gamma}^1(\Psi)(t,x,[\theta]) + \partial_x^\gamma \Psi(t,x,[\theta]),\\
    Z_{\gamma,\alpha}^4(\Psi)(t,x,[\theta])&:=&(Z_{\alpha}^2( Z_{\gamma}^1(\Psi))(t,x,[\theta]))+ Z_{\alpha}^2( \partial_x^\gamma \Psi)(t,x,[\theta]),\\
    Z_{\mu,\alpha}^1(\Psi)(t,x,[\theta],\boldsymbol{v})&:=& Z_{\mu,(\alpha_j)}^1\circ\cdots \circ Z_{\mu,(\alpha_2)}^1 \circ Z_{\mu,(\alpha_1)}^1(\Psi)(t,x,[\theta],\boldsymbol{v}),\\
    Z_{\mu,\alpha}^2(\Psi)(t,x,[\theta],\boldsymbol{v})&:=& Z_{\mu,\alpha}^1(\Psi)(t,x,[\theta],\boldsymbol{v})+ \partial_{\mu}^{\alpha} \Psi(t,x,[\theta]),\boldsymbol{v}),\\
	Z_{\mu,\gamma,\alpha}^3(\Psi)(t,x,[\theta],\boldsymbol{v})&:=&Z_{\gamma}^2(Z_{\mu,\alpha}^1(\Psi) ) (t,x,[\theta],\boldsymbol{v} )+Z_{\gamma}^2(\partial_{\mu}^{\alpha} \Psi)(t,x,[\theta],\boldsymbol{v}).
    \de
	
\ba\label{a}
Assume that there exists $a\in (0,1)$ such that the limit
	 \ce
	 r_1=\lim_{\lambda\to \infty}\frac{1}{\lambda^{a} }\int_{\Xi} ( e^{-\lambda \psi(u)}-1)\lambda(du)
	 \de
	 exists and belongs to $(-\infty,0)$.
\ea	
     \bl\label{G}
	 Under Assumption \ref{ellipticity} and Assumption \ref{a},  for $p\geq 2$ we have
	\ce
	\mE[ |\Gamma[X_t^{x,[\theta]}]|^{-p}] \leq Ct^{-\frac{p}{a}},
	\de
where $C=C_p>0$ is a constant.  Therefore, in particular, the solution $X_t^{x,[\theta]}$ of Eq.~(\ref{SDE}) admits a density
in $C^{\infty}(\mR^N)$ for all $t \geq s$.
	 \el
	 	\begin{proof}
	 		From \cite[page 216]{BL}, there exists a random variable $V$ such that $ V^{-1} \in \bigcap_{p\geq 2} L^p(\mP)$ with
	 		\ce
	 		\Gamma[X_t^{x,[\theta]}]\geq V \int_0^t \int_{\cO} \psi(u)N(ds,du)I_d.
	 		\de
	 		Hence
	 		\ce
	 		\mE[|(\Gamma[X_t^{x,[\theta]}])^{-1}|^p]&\leq& \mE[|V \int_0^t \int_{\cO} \psi(u)N(ds,du)|^{-p} ]\\
	 		&\leq &  \mE[ \big((\int_0^t \int_{\cO} \psi(u)N(ds,du))^{-1}\big)^{p} (V^{-1})^{p} ]\\
	 		&\leq &( \mE[ \big((\int_0^t \int_{\cO} \psi(u)N(ds,du))^{-1}\big)^{2p}])^{\frac{1}{2}} (\mE[(V^{-1})^{2p}])^{\frac{1}{2}}\\
	 		&\leq &C( \mE[ \big((\int_0^t \int_{\cO} \psi(u)N(ds,du))^{-1}\big)^{2p}])^{\frac{1}{2}} .
	 		\de
	 		Set $U=(\int_0^t \int_{\cO} \psi(u)N(ds,du))^{-1}$,  we have by \cite[Lemma 8.29]{BL},
	 		\ce
	 		\mP( U\geq x) \underset{x \rightarrow+\infty} {\sim} e^{-\frac{1-a}{a} (-a r_1)^{\frac{a}{1-a}}  t^{\frac{1}{1-a}} x^{b}}:=e^{r_2 x^{b}},
	 		\de
	 		where $b$ and $r_2$ satisfy:
	 		\ce
	 		\frac{1}{a}=\frac{1}{b}+1 \quad \text{and} \quad |atr_1|^{\frac{1}{a}}=|btr_2|^{\frac{1}{b}}.
	 		\de
	 		Therefore
	 		\ce
	 		\mE[ |U|^p]
	 		&=& p\int_0^\infty |x|^{p-1}\mP(U\geq x) dx\\
	 		&\leq& p\int_0^\infty |x|^{p-1}e^{-\frac{1-a}{a} (-a r_1)^{\frac{a}{1-a}}  t^{\frac{1}{1-a}} x^{b}} dx\\
	 		&=& p\int_0^\infty |x|^{p-1}e^{ -c x^{b}} dx\\
	 		&=&p\int_{0}^{\infty} c^{-\frac{p-1}{b}} u^{\frac{p-1}{b}} e^{-u}  \frac{1}{b} c^{-\frac{1} {b}} u^{\frac{1}{b}-1} d u\\
	 		&=&p\frac{1}{b} c^{-\frac{p}{b}} \int_{0}^{\infty} u^{\frac{p}{b}-1} e^{-u} d u\\
	 		&=&p\frac{1}{b} c^{-\frac{p}{b}} \Gamma( \frac{p}{b}),
	 		\de
	 		where $c=\frac{1-a}{a} (-a r_1)^{\frac{a}{1-a}}  t^{\frac{1}{1-a}}$.

        Then from \cite[Theorem 7.22, Proposition 8.28]{BL}, it follows that the solution $X_t^{x,[\theta]}$ of Eq. (\ref{SDE}) admits a density in $C^{\infty}(\mR^N)$.
	 	\end{proof}

We will need still some more lemmas.

	\bl\label{guocheng}
	Suppose that  $ b \in\cC_{b,Lip}^{k,k}(\mR^N\times\cP_2(\mR^N);\mR^N)$ and $c\in\cC_{b,Lip}^{k,k}(\mR^N\times \Xi\times\cP_2(\mR^N);\mR^N)$, $\Psi\in\mK_r^q(\mR,n)$, and  Assumption \ref{ellipticity} and Assumption \ref{a} are satisfied. Then
	\ce
	Z_{\gamma}^1(\Psi)(t,x,[\theta])&\in&\mK_{r-\frac{2}{a}\#\gamma}^{ q+2\#\gamma}(\mR, [(k-2) \wedge (n-1)]-\#\gamma ),\\
	Z_{\gamma}^2(\Psi)(t,x,[\theta])&\in&\mK_{r-\frac{2}{a}\#\gamma}^{ q+\#\gamma}(\mR, [(k-2) \wedge (n-1)]-\#\gamma ),\\
	Z_{\gamma}^3(\Psi)(t,x,[\theta]) &\in& \mK_{r-\frac{2}{a}\#\gamma}^{q+2\#\gamma}(\mR,[(k-2)\wedge (n-1)]-\#\gamma),\\
	Z^4_{\gamma,\alpha}(\Psi)(t,x,[\theta]) &\in& \mK_{ r-\frac{2}{a}(\#\gamma+\#\alpha)}^{ q+2\#\gamma+\#\alpha}(\mR, [(k-2) \wedge (n-1)]-\#\gamma-\#\alpha),\\
	Z_{\mu,\alpha}^1(\Psi)(t,x,[\theta],\boldsymbol{v})&\in&\mK_{r-\frac{2}{a}\#\alpha}^{ q+2\#\alpha}(\mR, [(k-2) \wedge (n-1)]-\#\alpha ),\\
	Z_{\mu,\alpha}^2(\Psi)(t,x,[\theta],\boldsymbol{v})&\in&\mK_{r-\frac{2}{a}\#\alpha}^{ q+2\#\alpha}(\mR, [(k-2) \wedge (n-1)]-\#\alpha ),\\
	Z_{\mu,\gamma,\alpha}^3(\Psi) (t,x,[\theta],\boldsymbol{v})&\in& \mK_{r-\frac{2}{a}(\#\gamma+\#\alpha)}^{ q+\#\gamma+2\#\alpha}(\mR, [(k-2) \wedge (n-1)]-\#\gamma-\#\alpha).
	\de
	\el
	The proof is provided in the Appendix \ref{appB}.

	\bl\label{bar}
	Under Assumption \ref{R}, for $m\in \mN$, we have $X_t^{x,[\theta]}\in \cH_{\overline{\mD}^{m,\infty}}^N$.
	\el
	\begin{proof}
		Since  $X_t^{x,[\theta]}$ satisfies a  standard  classical SDE with jumps, this is a straight result form \cite[Proposition 8.26]{BL}.
	\end{proof}

	\bl
		Under Assumption \ref{R} and Assumption \ref{ellipticity}, let $\Psi(t,x,\theta) \in \mK_r^q(\mR;n)$. Then
	\begin{enumerate}
		\item $DX_t^{x,[\theta]}\in (\overline{\mD}^{m,\infty})^{N}$;
	    \item $DX_t^{x,\theta}(\Gamma[X_t^{x,[\theta]}])^{-1}   \Psi(t,x,[\theta])\in Dom(\delta)$;
		\item $DX_t^{x,\theta}(\Gamma[X_t^{x,[\theta]}])^{-1}  \partial_x X_t^{x,[\theta]} \Psi(t,x,[\theta])\in Dom(\delta)$;
		\item $DX_t^{x,\theta}(\Gamma[X_t^{x,[\theta]}])^{-1}  \partial_\mu X_t^{x,[\theta]}(\boldsymbol{v})  \Psi(t,x,[\theta])\in Dom(\delta)$.
	\end{enumerate}
	\el
	\begin{proof}
	From Theorem \ref{KS}, Lemma \ref{bar} and $\Psi \in \mK_r^q(\mR;M)$,  we know that $ DX_t^{x,[\theta]} \in (\overline{\mD}^{m,\infty})^{N}$, $ \partial_x X_t^{x,[\theta]} $, $\partial_\mu X_t^{x,[\theta]}(\boldsymbol{v}) $ and $\Psi(t,x,[\theta])$ are in $(\mD^{m,\infty})^{N\times N}$. By \cite[Lemma 7.7, Corollary 7.9, Lemma 7.18]{BL} and  Assumption \ref{ellipticity}, we have
	\ce
	(\Gamma[X_t^{x,[\theta]}])^{-1}   \Psi(t,x,[\theta])\in (\mD^{m,\infty})^{N\times N},\\
	(\Gamma[X_t^{x,[\theta]}])^{-1}  \partial_x X_t^{x,[\theta]}  \Psi(t,x,[\theta])\in (\mD^{m,\infty})^{N\times N},\\
	(\Gamma[X_t^{x,[\theta]}])^{-1}  \partial_\mu X_t^{x,[\theta]}(\boldsymbol{v}) \Psi(t,x,[\theta])\in(\mD^{m,\infty})^{N\times N}.
	\de
	Finally, by \cite[lemma 7.20]{BL}, the proof is completed.
	\end{proof}

		Now, we first prove the integration by parts formulas with respect to parameter $x$.

	\bt\label{IPF variable} Under Assumption \ref{ellipticity} and Assumption \ref{a}, suppose that $f\in\cC_b^{\infty}(\mR^N,\mR)$ and $\Psi\in\mK_r^q(\mR,n)$. We have:
\begin{enumerate}
		\item
		If $\#\gamma\leq (n-1)\wedge (k-2)$, then
		\ce
		\mE[\partial_x^{\gamma}(f(X_t^{x,[\theta]}))\Psi(t,x,[\theta])]=\mE[f(X_t^{x,[\theta]})Z_{\gamma}^1(\Psi)(t,x,[\theta])].
		\de
		\item
		If $\#\gamma\leq (n-1)\wedge(k-2)$, then
		\ce
		\mE[(\partial^{\gamma}f)(X_t^{x,[\theta]})\Psi(t,x,[\theta])]=\mE[f(X_t^{x,[\theta]})Z_{\gamma}^2(\Psi)(t,x,[\theta])].
		\de
		\item
		If $\#\gamma\leq (n-1)\wedge (k-2)$, then
		\ce
		\partial_x^{\gamma}\mE[f(X_t^{x,[\theta]})\Psi(t,x,[\theta])]=\mE[f(X_t^{x,[\theta]})Z_{\gamma}^3(\Psi)(t,x,[\theta])].
		\de
		\item
		If $\#\gamma+\#\alpha\leq (n-1)\wedge(k-2)$, then
		\ce
		\partial_x^{\gamma}\mE[(\partial^{\alpha}f)(X_t^{x,[\theta]})\Psi(t,x,[\theta])]=\mE[f(X_t^{x,[\theta]})Z_{\gamma,\alpha}^4(\Psi)(t,x,[\theta])].
		\de
	\end{enumerate}
	\et
	\begin{proof}
(1)
 According to  \cite[Theorem 5.1]{BL}, we have
		\ce
		\widehat{\mE}[ DX_t^{x,[\theta]} (DX_t^{x,[\theta]})^{*}   ]=\Gamma[X_t^{x,[\theta]} ].
		\de
		This formula allows us to make the following computation for $f\in \cC_b^{\infty}(\mR^N;\mR)$, $i\in \{ 1,\cdots,N\}$
		\ce
		 &&\mE[ \partial_{x} f(X_t^{x,[\theta]})\Psi(t,x,[\theta])]\\
		 &=&\mE[\big(\partial f(X_t^{x,[\theta]})  \partial_x X_t^{x,[\theta]}\big)_i \Psi(t,x,[\theta])]\\
		 &=& \mE[\big(\partial f(X_t^{x,[\theta]}) \widehat{\mE}[DX_t^{x,[\theta]} (DX_t^{x,[\theta]})^{*} ](\Gamma[X_t^{x,[\theta]}])^{-1}  \partial_x X_t^{x,[\theta]} \big)_i  \Psi(t,x,[\theta])]\\
		 &=& \mE\widehat{\mE}[\big(\partial f(X_t^{x,[\theta]}) (DX_t^{x,[\theta]}) (DX_t^{x,[\theta]})^{*} (\Gamma[X_t^{x,[\theta]}])^{-1}  \partial_x X_t^{x,[\theta]} \big)_i \Psi(t,x,[\theta])]\\
		 &=& \mE\widehat{\mE}[Df(X_t^{x,[\theta]})\big((DX_t^{x,[\theta]})^{*} (\Gamma[X_t^{x,[\theta]}])^{-1}  \partial_x X_t^{x,[\theta]} \big)_i \Psi(t,x,[\theta])]\\
		 &=& \mE[(f(X_t^{x,[\theta]})) \delta(\big((DX_t^{x,[\theta]})^{*} (\Gamma[X_t^{x,[\theta]}])^{-1}  \partial_x X_t^{x,[\theta]} \big)_i \Psi(t,x,[\theta]))].
		 \de
		This proves the result for the case $\#\gamma=1$. By Lemma \ref{guocheng}, when $\#\gamma=1$,  $Z^1_{\gamma}(\Psi)\in \mK_{r-\frac{2}{a}}^{q+2}(\mR, [(k-2)\wedge (n-1)]-1) $. Repeating the same computation  $\#\gamma-1$ times, we conclude that the result holds for all multi-indices $\gamma$ such that $\#\gamma \leq (k-2)\wedge (n-1)$.

(2) With standard notation, we consider the column vector:
		\ce
		\Gamma[f(X_t^{x,[\theta]}), X_t^{x,[\theta]}]=\left(\Gamma\left[f(X_t^{x,[\theta]}), (X_t^{x,[\theta]})_{i}\right]\right)_{1 \leqslant i \leqslant N} .
		\de
		As a consequence  of the functional calculus related to the local Dirichlet forms \cite[Section 1.6]{BN}, we have for all ~$j\in\{ 1,\cdots,N\}$:
		\ce
		\Gamma[f(X_t^{x,[\theta]}), (X_t^{x,[\theta]})_i]=\sum_{i=1}^N \partial_j f(X_t^{x,[\theta]})\left(\Gamma\left[(X_t^{x,[\theta]})_j, (X_t^{x,[\theta]})_{i}\right]\right)_{1 \leqslant i \leqslant N},
		\de
namely,
		\ce
		 \Gamma[f(X_t^{x,[\theta]}),X_t^{x,[\theta]}]=\Gamma[X_t^{x,[\theta]},X_t^{x,[\theta]}] \partial f(X_t^{x,[\theta]}).
		 \de
Hence
		 \ce
		 \partial f(X_t^{x,[\theta]})=( \Gamma[X_t^{x,[\theta]},X_t^{x,[\theta]}])^{-1} \Gamma[f(X_t^{x,[\theta]}), X_t^{x,[\theta]}].
		 \de
        Then
		\ce
		&&\mE[(\partial^{i})f(X_t^{x,[\theta]}) \Psi(t,x,[\theta])]\\
		&=&\mE[ \big(\Gamma[f(X_t^{x,[\theta]}), (X_t^{x,[\theta]})^* ] (\Gamma[ X_t^{x,[\theta]}])^{-1}\big)_i  \Psi(t,x,[\theta]) ]\\
		&=&\mE\widehat{\mE}[\big(\langle Df(X_t^{x,[\theta]}), (X_t^{x,[\theta]})^{*}\rangle (\Gamma[ X_t^{x,[\theta]}])^{-1}\big)_i  \Psi(t,x,[\theta]) ]\\
		&=&\mE[f(X_t^{x,[\theta]})\delta(\big((D X_t^{x,[\theta]})^{*} (\Gamma[ X_t^{x,[\theta]}])^{-1}\big)_i  \Psi(t,x,[\theta]) ) ].
		\de
		By Lemma \ref{guocheng}, when $\#\gamma=1$, $ Z_{\gamma}^2 (\Psi) \in \mK_{r-\frac{2}{a}}^{q+2}(\mR,[(k-2)\wedge (n-1)]-1)$.
By iterating this procedure $\#\gamma-1$ times the result extends to the case for all $\#\gamma\leq (k-2)\wedge (n-1).$
		
(3) Here
		\ce
		&&\partial_x^\gamma\mE[ f(X_t^{x,[\theta]})  \Psi(t,x,[\theta])]\\
		&=&\mE\big[ \partial_x^\gamma[(f(X_t^{x,[\theta]})  \Psi(t,x,[\theta])]\big]\\
		&=& \mE\big[ \partial_x^\gamma f(X_t^{x,[\theta]})  \Psi(t,x,[\theta])+ f(X_t^{x,[\theta]})  \partial_x^\gamma \Psi(t,x,[\theta]) \big]\\
		&=& \mE\big[  f(X_t^{x,[\theta]})  Z_{\gamma}^1(\Psi)(t,x,[\theta]) + f(X_t^{x,[\theta]})  \partial_x^\gamma \Psi(t,x,[\theta])\big]\\
		&=& \mE\big[  f(X_t^{x,[\theta]})  (Z_{\gamma}^1(\Psi)(t,x,[\theta]) + \partial_x^\gamma \Psi(t,x,[\theta]))\big]\\
        &=& \mE\big[  f(X_t^{x,[\theta]}) Z_{\gamma}^3(\Psi)(t,x,[\theta]) \big].
		\de
		By Lemma \ref{guocheng}, when $\#\gamma=1$,  $Z^3_{\gamma}(\Psi)\in \mK_{r-\frac{2}{a}}^{q+2}(\mR, [(k-1)\wedge (n-1)]-1) $. We can therefore iterate this argument another $\#\gamma-1$ times to obtain the results for all $\gamma$ satisfying $\#\gamma \leq (k-2)\wedge (n-1)$.

		(4)
 From parts (2) and (3),
		\ce
		&&\partial_x^{\gamma}\mE[(\partial^{\alpha}f)(X_t^{x,[\theta]})\Psi(t,x,[\theta])]\\
		&=& \mE[\partial_x^{\gamma}[(\partial^{\alpha}f)(X_t^{x,[\theta]}) \Psi(t,x,[\theta])]]\\
		&=& \mE[\partial_x^{\gamma}(\partial^{\alpha}f)(X_t^{x,[\theta]})\Psi(t,x,[\theta])+ (\partial^{\beta}f)(X_t^{x,[\theta]}) \partial_x^{\gamma} \Psi(t,x,[\theta])]\\
		&=& \mE[\partial_x^{\gamma}(\partial^{\alpha}f)(X_t^{x,[\theta]}) \Psi(t,x,[\theta])\big]+\mE\big[ (\partial^{\beta}f)(X_t^{x,[\theta]}) \partial_x^{\gamma} \Psi(t,x,[\theta])]\\
		&=&\mE[\partial_x^{\gamma}(\partial^{\alpha}f)(X_t^{x,[\theta]})\Psi(t,x,[\theta])\big]+\mE[ f(X_t^{x,[\theta]}) Z_{\alpha}^2( \partial_x^\gamma \Psi)(t,x,[\theta])]\\
		&=& \mE[(\partial^{\alpha}f)(X_t^{x,[\theta]})Z_{\gamma}^1(\Psi)(t,x,[\theta]))\big]+\mE[ f(X_t^{x,[\theta]}) Z_{\alpha}^2( \partial_x^\gamma \Psi)(t,x,[\theta])]\\
		&=& \mE[f(X_t^{x,[\theta]}) Z_{\alpha}^2( Z_{\gamma}^1(\Psi))(t,x,[\theta]))\big]+\mE[ f(X_t^{x,[\theta]}) Z_{\alpha}^2( \partial_x^\gamma \Psi)(t,x,[\theta])]\\
		&=&\mE[f(X_t^{x,[\theta]}) (Z_{\alpha}^2( Z_{\gamma}^1(\Psi))(t,x,[\theta]))+ Z_{\alpha}^2( \partial_x^\gamma \Psi)(t,x,[\theta]))]\\
        &=&\mE[f(X_t^{x,[\theta]}) Z_{\gamma,\alpha}^4(\Psi)(t,x,[\theta])].
		\de
	\end{proof}
We now consider the derivatives of the function ~$[\theta] \mapsto \mE[f( X_t^{x,[\theta]})] $.

	\bt\label{IPF measure}
	Under Assumption \ref{ellipticity} and Assumption \ref{a}, suppose $f\in\cC_b^{\infty}(\mR^N,\mR)$ and $\Psi\in\mK_r^q(\mR,n)$. We have:
	\begin{enumerate}
		\item
		If $\#\alpha\leq (n-1)\wedge(k-2)$, then
		\ce
		\mE[\partial_{\mu}^{\alpha}(f(X_t^{x,[\theta]}))(\boldsymbol{v})\Psi(t,x,[\theta])]=\mE[f(X_t^{x,[\theta]})Z_{\mu,\alpha}^1(\Psi)(t,x,[\theta],\boldsymbol{v})].
		\de
		\item
		If $\#\alpha\leq (n-1)\wedge(k-2)$, then
		\ce
		\partial_{\mu}^{\alpha}\mE[f(X_t^{x,[\theta]})\Psi(t,x,[\theta])](\boldsymbol{v})=\mE[f(X_t^{x,[\theta]})Z_{\mu,\alpha}^2(\Psi)(t,x,[\theta],\boldsymbol{v})].
		\de
		\item
		If $\#\alpha+\#\gamma\leq (n-1)\wedge(k-2)$, then
		\ce
		\partial_{\mu}^{\alpha}\mE[(\partial^{\gamma}f)(X_t^{x,[\theta]})\Psi(t,x,[\theta])](\boldsymbol{v})=\mE[f(X_t^{x,[\theta]})Z_{\mu,\gamma,\alpha}^3(\Psi)(t,x,[\theta],\boldsymbol{v})].
		\de
	\end{enumerate}
	\et
	\begin{proof}

 (1)  Here we use again the relationship
        \ce
		\widehat{\mE}[ DX_t^{x,[\theta]}\cdot (DX_t^{x,[\theta]})^{*}   ]=\Gamma[X_t^{x,[\theta]} ].
		\de
       For $\alpha=(i)$,
 		\ce
 		&&\mE[\partial_{\mu}(f(X_t^{x,[\theta]}))(v)\Psi(t,x,[\theta])]\\
 		&=&\mE[\big(\partial f(X_t^{x,[\theta]})  \partial_{\mu}(X_t^{x,[\theta]})(v)\big)_{i}  \Psi(t,x,[\theta]) ]\\
 		&=&\mE[\big(\partial f(X_t^{x,[\theta]}) \widehat{\mathbb{E}}[DX_{t}^{x,[\theta]}(DX_{t}^{x,[\theta]})^{*}](\Gamma[X_{t}^{x,[\theta]}])^{-1}  \partial_{\mu}X_t^{x,[\theta]}(v)\big)_{i}  \Psi(t,x,[\theta]) ]\\
 		&=&\mE\widehat{\mathbb{E}}[\big(\partial f(X_t^{x,[\theta]})  DX_{t}^{x,[\theta]}(DX_{t}^{x,[\theta]})^{*}(\Gamma[X_{t}^{x,[\theta]}])^{-1}   \partial_{\mu}X_t^{x,[\theta]}(v)\big)_{i}  \Psi(t,x,[\theta]) ]\\
 		&=&\mE\widehat{\mathbb{E}}[ Df(X_t^{x,[\theta]}) \big((DX_{t}^{x,[\theta]})^{*}(\Gamma[X_{t}^{x,[\theta]}])^{-1}  \partial_{\mu}X_t^{x,[\theta]}(v)\big)_{i}  \Psi(t,x,[\theta]) ]\\
 		&=& \mE[f(X_t^{x,[\theta]}) \delta(\big((DX_{t}^{x,[\theta]})^{*}(\Gamma[X_{t}^{x,[\theta]}])^{-1}  \partial_{\mu}X_t^{x,[\theta]}(v)\big)_{i}  \Psi(t,x,[\theta]) )].
 		\de
 		This proves for case $\#\alpha=1$, we can iterate this argument another ~$\#\alpha-1$ times to obtain the result for all ~$\#\alpha$.

 (2) By the chain rule.
 		\ce
 		&&\partial_{\mu}^{\alpha}\mE[f(X_t^{x,[\theta]})\Psi(t,x,[\theta])](\boldsymbol{v})\\
 		&=& \mE[\partial_{\mu}^{\alpha}f(X_t^{x,[\theta]})\Psi(t,x,[\theta])](\boldsymbol{v})]+\mE[f(X_t^{x,[\theta]}) \partial_{\mu}^{\alpha} \Psi(t,x,[\theta])(\boldsymbol{v})]\\
 		&=&\mE[ f(X_t^{x,[\theta]}) Z_{\mu,\alpha}^1(\Psi)(t,x,[\theta],\boldsymbol{v}) ]+\mE[f(X_t^{x,[\theta]}) \partial_{\mu}^{\alpha} \Psi(t,x,[\theta])(\boldsymbol{v})]\\
 		&=&\mE[f(X_t^{x,[\theta]})  \big(Z_{\mu,\alpha}^1(\Psi)(t,x,[\theta],\boldsymbol{v})+ \partial_{\mu}^{\alpha} \Psi(t,x,[\theta])(\boldsymbol{v})\big) ]\\
 		&=& \mE[f(X_t^{x,[\theta]})  Z_{\mu,\alpha}^2(\Psi)(t,x,[\theta],\boldsymbol{v})],
 		\de
 		where
 		\ce
 		Z_{\mu,\alpha}^2(\Psi)(t,x,[\theta],\boldsymbol{v})= Z_{\mu,\alpha}^1(\Psi)(t,x,[\theta],\boldsymbol{v})+ \partial_{\mu}^{\alpha} (\Psi)(t,x,[\theta],\boldsymbol{v}).
 		\de

 (3) From part (1) and (2), and Theorem \ref{IPF variable} part (2),
 		\ce
 		&&\partial_{\mu}^{\alpha}\mE[(\partial^{\gamma}f)(X_t^{x,[\theta]})\Psi(t,x,[\theta])](\boldsymbol{v})\\
 		&=& \mE[(\partial_{\mu}^{\alpha}\partial^{\gamma}f)(X_t^{x,[\theta]})(\boldsymbol{v}) \Psi(t,x,[\theta])]+\mE[(\partial^{\gamma}f)(X_t^{x,[\theta]})  \partial_{\mu}^{\alpha} \Psi(t,x,[\theta]) (\boldsymbol{v})]\\
 		&=& \mE[\partial_{\mu}^{\alpha}(\partial^{\gamma}f(X_t^{x,[\theta]}))(\boldsymbol{v}) \Psi(t,x,[\theta])]+\mE[(\partial^{\gamma}f)(X_t^{x,[\theta]}) \partial_{\mu}^{\alpha} \Psi(t,x,[\theta]) (\boldsymbol{v})]\\
 		&=&\mE[(\partial^{\gamma})f(X_t^{x,[\theta]}) Z_{\mu,\alpha}^1(\Psi)(t,x,[\theta],\boldsymbol{v} ) ]+\mE[f(X_t^{x,[\theta]}) Z_{\gamma}^2(\partial_{\mu}^{\alpha} \Psi)(t,x,[\theta],\boldsymbol{v}) ]\\
 		&=&\mE[f(X_t^{x,[\theta]}) Z_{\gamma}^2(Z_{\mu,\alpha}^1(\Psi) ) (t,x,[\theta],\boldsymbol{v} ) ]+\mE[f(X_t^{x,[\theta]}) Z_{\gamma}^2(\partial_{\mu}^{\alpha} \Psi)(t,x,[\theta],\boldsymbol{v}) ]\\
 		&=& \mE[f(X_t^{x,[\theta]}) Z_{\mu,\gamma,\alpha}^3(\Psi)(t,x,[\theta],\boldsymbol{v}) ],
 		\de
 		where
 		\ce
 		Z_{\mu,\gamma,\alpha}^3(\Psi)(t,x,[\theta],\boldsymbol{v})=Z_{\gamma}^2(Z_{\mu,\alpha}^1(\Psi) ) (t,x,[\theta],\boldsymbol{v} )+Z_{\gamma}^2(\partial_{\mu}^{\gamma} \Psi)(t,x,[\theta],\boldsymbol{v}).
 		\de

	\end{proof}

	We now derive the integration by parts formulas for the derivatives of the following mapping
	\ce
	x\mapsto\mE[f(X_t^{x,\delta_x})].
	\de

	\bt\label{IPF}
	Under Assumption \ref{ellipticity} and Assumption \ref{a}, let $f\in\cC_b^{\infty}(\mR^N,\mR)$. Then for all multi-indices $\gamma$ on $\{1,\cdots,N\}$ with $ \#\gamma \leq (k-2)\wedge (n-1)$,
	\ce
	\partial_x^{\gamma}\mE[f(X_t^{x,\delta_x})]=\mE[f(X_t^{x,\delta_x}) (Z_{\gamma}^1(1)(t,x,\delta_x) +Z_{\mu,\gamma}^1(1)(t,x,\delta_x) ) ].
	\de
	In particular, the following estimate is satisfied:
	\ce
	|\partial_x^{\gamma}\mE[f(X_t^{x,\delta_x})]|\leq C\|f\|_{\infty}t^{-\frac{\#\gamma}{a} }(1+|x|)^{2\#\gamma}.
	\de
	\et
	\begin{proof}
		Since $\delta_x$ depends on $x$, we have
		\ce
		\partial_x^i\mE[f(X_t^{x,\delta_x})]=\partial_z^i\mE[f(X_t^{z,\delta_x})]|_{z=x}+\partial_{\mu}^i\mE[f(X_t^{x,[\theta]})(v)]|_{[\theta]=\delta_x,v=x}.
		\de
		Then for $\#\gamma=1$, the  following equalities  by Theorem \ref{IPF variable} and Theorem \ref{IPF measure}.
		\ce
		\partial_z^i\mE[f(X_t^{z,\delta_x})]|_{z=x}
		&=&\mE[ f(X_t^{x,\delta_x}) Z^1_{(i)}(1)(t,x,\delta_x) ],\\
		\partial_{\mu}^i\mE[f(X_t^{x,[\theta]})(v)]|_{[\theta]=\delta_x,v=x}
		&=&\mE[ f(X_t^{x,[\theta]}) Z_{\mu,(i)}^1(1)(t,x,\delta_x)].
		\de
		The proof is completed by repeating this procedure for another $\#\gamma-1$ times.
	\end{proof}
	
	As a direct consequence of Theorem \ref{IPF} and Theorem \ref{IPF variable}, we have the following result.
	
	\bc\label{IPF general}
	Under Assumption \ref{ellipticity} and Assumption  \ref{a}, let $f\in\cC_b^{\infty}(\mR^N,\mR)$, $\gamma$ and $\beta$ are multi-indices on $\{1,\cdots,N\}$ with $\#\gamma+\#\alpha \leq (k-2)\wedge (n-1)$. Then
	\ce
	\partial_x^{\gamma}\mE[(\partial^{\alpha}f)(X_t^{x,\delta_x})]=\mE[f(X_t^{x,\delta_x}) Z_{\alpha}^2(Z_{\gamma}^1(1) +Z_{\mu,\gamma}^1(1))(t,x,\delta_x)].
	\de
	and
	\be\label{ss}
	Z_{\alpha}^2(Z_{\gamma}^1(1) +Z_{\mu,\gamma}^1(1))(t,x,\delta_x)\in\cK_{-\frac{2}{a}(\#\gamma+\#\alpha )}^{2\#\gamma+\#\alpha}(\mR,(n-1) \wedge (k-2)-\#\gamma-\#\alpha).
	\ee
	\ec
	
	\begin{proof}
		Theorem \ref{IPF} gives
		 \ce
		 \partial_x^{\gamma}\mE[(\partial^{\alpha}f)(X_t^{x,\delta_x})]=\mE[(\partial^{\alpha}f)(X_t^{x,\delta_x})  (Z_{\gamma}^1(1)(t,x,\delta_x) +Z_{\mu,\gamma}^1(1)(t,x,\delta_x) ) ].
		 \de
		 Then using part~$(2)$ of  Theorem \ref{IPF variable}, we get
		 \ce
		 &&\mE[(\partial^{\alpha}f)(X_t^{x,\delta_x}) (Z_{\gamma}^1(1)(t,x,\delta_x) +Z_{\mu,\gamma}^1(1)(t,x,\delta_x) ) ]\\
		 &=&\mE[f(X_t^{x,\delta_x}) Z_{\alpha}^2(Z_{\gamma}^1(1) +Z_{\mu,\gamma}^1(1))(t,x,\delta_x)].
		 \de
        Similar to Lemma \ref{guocheng}, we can easily prove (\ref{ss}).
	\end{proof}	
	
	\br
	Denote by $\cK_r^q(\mR,M)$  the set of all elements of  $\mK_r^q(\mR,M)$ that do not depend
on the measure variable $\mu\in \cP_2(\mR^N)$. From \cite[Lemma 2.11]{CM}, if $\Psi\in \mK_r^q(\mR,M)$, then $ \Psi(t,x):=\Psi(t,x,\delta_x) \in \cK_r^q(\mR,M) $.
	\er
	\section{Smoothness of Densities}
	
	We are now in a position to prove the existence and smoothness of the density for the solution to Eq. (\ref{SDE}).	The following lemma will be used in the proof of Theorem \ref{smoothness}.
	
	\bl\label{est4}
	Under  Assumption \ref{H_solution}, for any $R>0$, $p\geq 2$, $\theta\in L^2(\cF_0;\mR^N)$ and $\vartheta \geq 1$, there exists a constant $C>0$ such that for any $t\geq 0$,
	\ce
		\mathbb{P}(\sup_{s\in(0,t]}|X_s^{x,[\theta]}-x|\geq R)\leq C\frac{t^\vartheta(1+\|\theta\|_p^\vartheta)}{R^\vartheta}.
	\de
	\el
	
	\begin{proof}
		Using (\ref{est3}) and Chebyshev's inequality, for any $\vartheta \geq 1$, there exists a constant $C>0$ such that
		\ce
		\begin{aligned}
			\mathbb{P}(\sup_{s\in(0,t]}|X_s^{x,[\theta]}-x|\geq R)\leq\frac{\mathbb{E}[\sup_{s\in(0,t]}|X_s^{x,[\theta]}-x|^\vartheta]}{R^\vartheta}\leq C\frac{t^\vartheta(1+\|\theta\|_p^\vartheta)}{R^\vartheta}.
		\end{aligned}
		\de
	\end{proof}

	\bt\label{smoothness}
	Under  Assumption \ref{ellipticity} and Assumption \ref{a}, suppose that $b \in\cC_{b,Lip}^{k,k}(\mR^N\times\cP_2(\mR^N);\mR^N)$ and  $c\in\cC_{b,Lip}^{k,k}(\mR^N\times\Xi\times\cP_2(\mR^N);\mR^N)$. Let $\gamma,\beta$ be multi-indices on $\{1,\cdots,N\}$, and let $\# \gamma+\#\beta+N+1 \leq k$. Then
	\begin{enumerate}
		\item
		for all $t\in[0,T]$, $X_t^{x,\delta_x}$ has density $p(t,x,\cdot)$ such that
		\ce
		(x,y)\mapsto\partial_x^{{\gamma}}\partial_y^{\beta}p(t,x,y)
		\de
		is continuous;
		\item
		moreover, for  any $\vartheta>0$ , there exist $C>0$ and  $\tilde{\vartheta}_1 >\frac{2\#\gamma+\#\beta+1}{a} $,
such that
			\ce
			\max_{0\leq2\#\gamma\leq k}|\partial_x^\gamma p(t,x,x)|&\leq& Ct^{-\frac{2\#\gamma+1}{a}}(1+|x|)^{2\#\gamma+1}
\de
and
\ce
    \max_{0\leq \#\gamma+\#\beta\leq k}|(1+|y-x|^2)^{k/2}\partial_x^\gamma\partial_y^\beta p(t,x,y)|&\leq& C\frac{t^{\tilde{\vartheta}_1 -\frac{2\#\gamma+\#\beta+1}{a}}(1+|x|)^{2\#\gamma+\#\beta+1}}{|y|^{\tilde{\vartheta}_1}}
		\de
		for all $(t,x,y) \in\{(t,x,y) \in (0,T] \times \Xi \times \mR^N; y\neq x \}$.
	\end{enumerate}
	\et
	
	\begin{proof}
		Let ~$ \eta=(1,2,\cdots,N)$ and define the multi-dimensional indicator function
		\ce
		{1}_{\left\{y_{0}>y\right\}}:=\prod_{i=1}^{N} {1}_{\left\{y_{0}^{i}>y^{i}\right\}}.
		\de
		For any ~$g\in \cC_0^{\infty}$, define the function ~$f$ by
		\be\label{6.3}
		f(y_0):=\int_{\mathbb{R}^N}g(y){1}_{\{y_0>y\}}dy.
		\ee
Then $f\in  \cC_p^{\infty}(\mR^N;\mR)$ and $ \partial^{\eta} f(y_0)=g(y_0)$.
		
		(1)  We have
			\be\label{density}
			&&\partial_x^{\gamma}\mathbb{E} \left[ \left( \partial^\beta g \right)\left(X_t^{x, \delta_x} \right) \right]\nonumber\\
			&=& \partial_x^{\gamma} \mathbb{E} \left[ \left( \partial^{\beta}( \partial^\eta f) \right)\left(X_t^{x, \delta_x} \right)\right]\nonumber\\
			&=& \mathbb{E} \left[ f\left(X_t^{x, \delta_x} \right) \cdot Z^2_{\beta+\eta} \left( {Z}^1_{\gamma}(1) + {Z}^1_{\mu,\gamma}(1) \right)(t,x,\delta_x) \right]\nonumber\\
			&=& \mE[\left(\int_{\mathbb{R}^N}g(y)\mathbf{1}_{\{X_t^{x,\delta x}>y\}}dy\right)  Z^2_{\beta+\eta} \left( {Z}^1_{\gamma}(1) + {Z}^1_{\mu,{\gamma}}(1) \right)(t,x,\delta_x)]\nonumber\\
			&=& \int_{\mR^N} g(y) \mE[\mathbf{1}_{\{X_t^{x,\delta x}>y\}} Z^2_{\beta+\eta} \left( {Z}^1_{\gamma}(1) + {Z}^1_{\mu,{\gamma}}(1) \right)(t,x,\delta_x) ]dy,
			\ee
			where each step separately uses the:~$ \partial^{\eta} f=g$; Corollary \ref{IPF general};  (\ref{6.3}), and Fubini's theorem. Then it is easy to check that, for any $R>0$ and $t\in [0,T]$, there exists $C=C(R,t)>0$ such that
            \ce
            \sup_{|x|\leq R}|\partial_x^{\alpha} \mE[(\partial^{\beta} g)(X_t^{x,\delta_x})]|\leq C\|g\|_{\infty}.
            \de
			From \cite[Lemma 3.1]{TS}, we know that $X_t^{x,\delta_x}$ has a density function $ p(t,x,\cdot)$ and $\partial_x^{\alpha} \partial_y^{\beta}p(t,x,y) $ exists. Once we know that a smooth density exists, it follows from (\ref{density}) that we can identify ~$ \partial_x^{{\gamma}}\partial_y^{\beta}p(t,x,y)$ as
		\ce
		\partial_x^{\gamma}\partial_y^{\beta}p(t,x,y)=(-1)^{|\beta|}\mE[\mathbf{1}_{\{X_t^{x,\delta x}>y\}} Z^2_{\beta+\eta} \left( {Z}^1_{\gamma}(1) + {Z}^1_{\mu,{\gamma}}(1) \right)(t,x,\delta_x) ].
		\de
		
(2)		By the chain rule and  part $(3)$ of Theorem \ref{IPF variable}, we have
		\ce
		\partial^{\gamma}_x\mE[\varphi(X_T^{x,\delta_x})]
		&=&\mE[\varphi(X_T^{x,\delta_x}) Z_{\gamma}^3(1)(t,x,\delta_x) ].
		\de
		Let
		\ce
		\mQ_{\gamma}^x (T,dy)=Z_{\gamma}^3(1)(t,x,\delta_x) \mP( (X_T^{x,\delta_x}-x)\in dy).
		\de
		Then
		\ce
		\partial^{\gamma}_x\mE[\varphi(X_T^{x,\delta_x})]=\int_{\mR^N} \varphi(x+y) \mQ_{\gamma}^x (T,dy).
		\de
		Now choose $\eta\in \cC^{\infty}(\mR^N)$ such that
		\ce
		\eta(z)=\left\{\begin{array}{ll}
			0, & |z| \leq \frac{1}{4};\\
			1, & |z|>\frac{1}{2},
		\end{array}\right.
		\de
		and define
		\ce
		\eta_{k, \epsilon}(y)=\left\{\begin{array}{ll}
			\left(1+|y|^{2}\right)^{k / 2}, & \text { if } \epsilon=0 ; \\
			\left(1+|y|^{2}\right)^{k / 2} \eta(y / \epsilon), & \text { if } \epsilon>0 .
		\end{array}\right.
		\de
		Then using Theorem \ref{IPF variable} again, we have
		\ce
		&&\int_{\mR^N} \partial^{\beta}\varphi(y) \eta_{k, \epsilon}(y)\mQ_{\gamma}^x (T,dy)\\
		&=& \mE[\partial^{\beta}\varphi(X_T^{x,\delta_x}-x) \eta_{k, \epsilon}(X_T^{x,\delta_x}-x) Z_{\gamma}^3(1)(T,x,\delta_x)]\\
		&=& \mE[\varphi(X_T^{x,\delta_x}-x)  Z^3_{\gamma,\beta}(X_T^{x,\delta_x}) ],
		\de
		where
		\ce
		Z^3_{\gamma,\beta}(X_T^{x,\delta_x})=Z_{\beta}^3(\eta_{k, \epsilon}(X_T^{x,\delta_x}-x) Z_{\gamma}^3(1))(T,x,\delta_x).
		\de
		By \cite[Theorem 5.3]{Shi},  for each $(T,x) \in (0,\infty) \times \Xi$,  there is a $q_\gamma^x (T,\cdot) \in \cC_b^{m-1}(\mR^N)$ such that
			\ce
			\eta_{k,\epsilon}(y)\mathbb{Q}_\gamma^x(T,dy)=\eta_{k,\epsilon}(y)q_\gamma^x(T,y)dy.
			\de
			 Moreover,
			\ce
			&&\|\eta_{k,\epsilon}(\cdot)q_\gamma^x(t,\cdot) \|_{\cC_b^{m-1}(\mR^N)}\\
			&\leq &  C\big(\|Z^3_{\gamma}(1)(t,x,\delta_x)\|_{p}+\|Z^3_{\gamma,\beta}(X_t^{x,\delta_x})\|_{p}\big)\|Z_{(i)}^3(1)(t,x,\delta_x)\|_{p}^{N-N/p}\\
			&\leq &
			\begin{cases}
				Ct^{-\frac{2\#\gamma+\#\beta+1}{a}}(1+|x|)^{2\#\gamma+\#\beta+1}, & \epsilon=0; \\ Ct^{-\frac{2\#\gamma+\#\beta+1}{a}}(1+|x|)^{2\#\gamma+\#\beta+1}\|\eta(\cdot/\epsilon)\|_{\mathcal{C}_b^{m+n}(\mathbb{R}^N)}\mathbb{P}(|X_t^{x,\delta_x}-x|\geq\frac{\epsilon}{4})^{\frac{1}{p_1}}, & \epsilon>0.
			\end{cases}
			\de
			The first inequality is the result of \cite[Theorem 5.3]{Shi}, and the second one is due to  Lemma \ref{guocheng}.
		
		Letting $\varepsilon=|y|$, we can obtain
		\ce
		&&\max_{0\leq\#\gamma+\#\beta\leq k} |\partial_y^\beta(1+|y|^2)^{k/2}q_\gamma^x(t,y)| \\
  &\leq&
		\begin{cases}
			Ct^{-\frac{2\#\gamma+\#\beta+1}{a}}(1+|x|)^{2\#\gamma+\#\beta+1}, & y=0; \\
			Ct^{-\frac{2\#\gamma+\#\beta+1}{a}}(1+|x|)^{2\#\gamma+\#\beta+1}\mathbb{P}(|X_t^{x,\delta_x}-x|\geq\frac{|y|}{4})^{\frac{1}{p_1}}, & y\neq 0,
		\end{cases}
		\de
		for all $(t,x,y)\in (0,T]\times U \times \mR^N$. Then it  can be seen from Lemma \ref{est4}  that
		\ce
			\mathbb{P}(|X_t^{x,\delta_x}-x|\geq\frac{|y|}{4})\leq C\frac{t^{\vartheta }}{|y|^{\vartheta}}.
		\de
		Let $\vartheta_1=\vartheta/p_1  $.
		Then
		\ce
		\max_{0\leq\#\gamma+\#\beta\leq k}|\partial_y^\beta\big((1+|y|^2)^{k/2}q_\gamma^x(t,y)\big)|\leq
		\begin{cases}
			Ct^{-\frac{2\#\gamma+\#\beta+1}{a}}(1+|x|)^{2\#\gamma+\#\beta+1}, & y=0; \\
			C\frac{t^{\vartheta_1 -\frac{2\#\gamma+\#\beta+1}{a}}(1+|x|)^{2\#\gamma+\#\beta+1}}{|y|^{\vartheta_1}}, & y\neq 0.
		\end{cases}
		\de
		By calculation, it follows that
		\be\label{est5}
		\max_{0\leq\#\gamma+\#\beta\leq k}|(1+|y|^2)^{k/2}\partial_y^\beta q_\gamma^x(t,y)|\leq\begin{cases}
			Ct^{-\frac{2\#\gamma+\#\beta+1}{a}}(1+|x|)^{2|\gamma|+|\beta|+1}, & y=0; \\
			C\frac{t^{\vartheta_1 -\frac{2\#\gamma+\#\beta+1}{a}}(1+|x|)^{2\#\gamma+\#\beta+1}}{|y|^{\vartheta_1}}, & y\neq 0.
		\end{cases}
		\ee
		
		Since
		\ce
			\partial^\gamma_x\mathbb{E}[\varphi(X_t^{x,[\theta]})] & =&\int_{\mathbb{R}^N}\varphi(x+y)\mathbb{Q}_\gamma^x(t,dy) \\
			& =&\int_{\mathbb{R}^N}\varphi(x+y)q_\gamma^x(t,y)dy \\
			& =&\int_{\mathbb{R}^N}\varphi(y)q_\gamma^x(t,y-x)dy,
		\de
		and
		\ce
		\partial^\gamma_x\mathbb{E}[\varphi(X_t^{x,[\theta]})] & =&\partial^\gamma\int_{\mathbb{R}^N}\varphi(y)p(t,x,y)dy \\
			& =&\int_{\mathbb{R}^N}\varphi(y)\partial^\gamma p(t,x,y)dy,
		\de
		so
		\ce
		q_\gamma^x(t,y-x)=\partial^\gamma p(t,x,y).
		\de
		This combined with the estimate $(\ref{est5})$ gives that, there exists $ \tilde{\vartheta}_1 >\frac{2\#\gamma+\#\beta+1}{a} $ such that
		\ce
		\begin{aligned}
			\max_{0\leq2\#\gamma\leq k}|\partial_x^\gamma p(t,x,x)|\leq Ct^{-\frac{2\#\gamma+1}{a}}(1+|x|)^{2\#\gamma+1}
		\end{aligned}
		\de
		and
		\ce
		\max_{0\leq\#\gamma+\#\beta\leq k}|(1+|y-x|^2)^{k/2}\partial_x^\gamma\partial_y^\beta p(t,x,y)|\leq C\frac{t^{\tilde{\vartheta}_1 -\frac{2\#\gamma+\#\beta+1}{a}}(1+|x|)^{2\#\gamma+\#\beta+1}}{|y|^{\tilde{\vartheta}_1}},
		\de
		for all $(t,x,y) \in\{(t,x,y) \in (0,T] \times \Xi \times \mR^N; y\neq x \}$.
	\end{proof}

	It is easy to generalize the result of Theorem \ref{smoothness} to $X_t^{\theta}$ with general initial distributions $[\theta]$.
	\bc
	Suppose that $b\in\cC_{b,Lip}^{k,k}(\mR^N\times\cP_2(\mR^N);\mR^N)$ and  $c\in\cC_{b,Lip}^{k,k}(\mR^N\times\Xi\times\cP_2(\mR^N);\mR^N)$, and  Assumption \ref{ellipticity} and Assumption \ref{a} are satisfied. Let $\theta\in L^p(\cF_0;\mR^N)$, then for any multi-index $\beta$ on $\{1,\cdots,N\}$ such that $\#\beta+N+2\leq k$, and for all $t\in[0,T]$, $X_t^{\theta}$ has a density $p(t,\cdot)$ such that $\partial_y^{\beta}p(t,y)$ exists and $y\mapsto\partial_y^{\beta}p(t,y)$ is continuous.
	\ec
	
	\section{Connection with PDE}
	We now turn to consider  PDE (\ref{pde1}). It is proven in \cite[Proposition 6.1]{HL} that  the derivatives of
$U(t,x,[\theta]):=\mE[g(X_t^{x,[\theta]},[X_t^{\theta}])]$ exists when $g$ is sufficiently smooth. The following lemma is a straight result from \cite[Proposition 6.1]{HL}.
	
	\bl\label{U}
	Assume that the function $g:\mR^N\times \cP_2(\mR^N)\to \mR^N$ admits continuous derivatives $\partial_x g$ and $\partial_\mu g$ satisfying for some $q>0$ and $0\leq p <2$,
	\ce
	\begin{aligned}
		\left|\partial_{x} g(x,[\theta])\right| & \leq C\left(1+|x|+\|\theta\|_{2}\right)^{q}, \\
		\left|\partial_{\mu} g(x,[\theta], v)\right| & \leq C\left(1+|x|^{q}+\|\theta\|_{2}^{q}+|v|^{p}\right),
	\end{aligned}
	\de
	and  we assume $b\in \cC_{b,Lip}^{1,1}(\mR^N\times \cP_2(\mR^N);\mR^N)$ and $ c\in \cC_{b,Lip}^{1,1}(\mR^N\times \Xi\times \cP_2(\mR^N);\mR^N)$. Then $\partial_\mu U$ exists and takes the following form:
	\ce
	\partial_\mu U(t,x,[\theta],v)&=&\mE[\partial g(X_t^{x,[\theta]},[X_t^{\theta}]) \partial_{\mu}X_t^{x,[\theta]}(v)  ]\\
	&&+\mE\mE^{\prime}[ \partial_{\mu} g(X_t^{x,[\theta]},[X_t^{\theta}],(X_t^{v,[\theta]})^{\prime}) \partial_v (X_t^{v,[\theta]})^{\prime}\\
	&&+\partial_{\mu} g(X_t^{x,[\theta]},[X_t^{\theta}],(X_t^{\theta^{\prime},[\theta]})^{\prime}) \partial_\mu (X_t^{\theta^{\prime},[\theta]})^{\prime}(v)].
	\de
	\el
	Now we introduce the following function class for which we will be able to develop integration by parts formulas. The following definition was first presented in \cite[Definition 5.3]{CM}.
	\bd
	[(IC)$_x$ and (IC)$_v$ ] We say that a function $g:\mR^N\times\cP_2(\mR^N)\to \mR $  is in the class (IC) if the following conditions hold:
	\begin{enumerate}
		\item $g$ is   continuous with polynomial growth: i.e., there exists $q>0$ such that for all $ (x,[\theta])\in\mathbb{R}^N\times\mathcal{P}_2(\mathbb{R}^N):$
		\ce
		|g(x,[\theta])|\leq C(1+|x|+\|\theta\|_2)^q;
		\de
		\item there exists a sequence of functions $(g_l)_{l\geq 1}: \mR^N\times \cP_2(\mR^N)\to \mR$ with  polynomial growth such that $g_l\to g$ uniformly on compacts and $\partial_x g_l$ exists and also has polynomial growth for each $l\geq 1$;
		\item for each $l\geq 1$ there exists a function $G_l:\mathbb{R}^N\times\mathcal{P}_2(\mathbb{R}^N)\times\mathbb{R}^N\to\mathbb{R}$ which is either differentiable in $x$ or $v$ and
		\ce
		\partial_{\mu}g_{l}(x,\mu,v)=\partial_{x}G_{l}(x,\mu,v)\quad\text{or}\quad \partial_{\mu}g_{l}(x,\mu,v)=\partial_{v}G_{l}(x,\mu,v).
		\de
		Moreover, each $G_l$ and its derivatives satisfy the growth condition: there exist $q>0$ and $0\leq r <1$ such that for all $ (x,[\theta],v)\in\mathbb{R}^{N}\times\mathcal{P}_{2}(\mathbb{R}^{N})\times\mathbb{R}^{N}:$
		\ce
		|h(x,[\theta],v)|\leq C\left(1+|x|^q+\|\theta\|_2^q+|v|^r\right),
		\de
		where $h$ is $G_l, \partial_x G_l$ or $ \partial_v G_l$. In addition, we assume that for all $(x,\mu,v)$ the pointwise limit $  \operatorname*{lim}_{l\rightarrow\infty}G_{l}(x,\mu,v)$ exists and the function $G$ defined by $ G(x,\mu,v):=\operatorname*{lim}_{l\rightarrow\infty}G_{l}(x,\mu,v)$ is  continuous and satisfies the same growth condition.
	\end{enumerate}
	If $ \partial_{\mu}g_{l}=\partial_{x}G_{l}$ we say $g$ is in the class (IC)$_x$. If $ \partial_{\mu}g_{l}=\partial_{v}G_{l}$ we say $g$ is in the class (IC)$_v$.
	\ed
	
	In order to prove the existence and uniqueness of a solution to PDE (\ref{pde1}), we need the following hypotheses.

	(H1): Assumption \ref{ellipticity},  Assumption \ref{a} hold, the coefficients $b\in \cC_{b,Lip}^{2,2}(\mR^N\times  \cP_2(\mR^N);\mR^N)$, $c\in \cC_{b,Lip}^{2,2}(\mR^N\times \Xi \times \cP_2(\mR^N);\mR^N)$ and $g: \mR^N\times\cP_2(\mR^N)\to \mR^N$ is in the class (IC)$_x$.
	
	(H2): Assumption \ref{ellipticity}, Assumption \ref{a} hold, the coefficients $b\in \cC_{b,Lip}^{2,2}(\mR^N\times  \cP_2(\mR^N);\mR^N)$, $c\in \cC_{b,Lip}^{2,2}(\mR^N\times \Xi \times \cP_2(\mR^N);\mR^N)$ as well as being uniformly bounded, and \ $g: \mR^N\times\cP_2(\mR^N)\to \mR^N$ is in the class (IC)$_v$.
	
	\bl\label{5.2}
	Under either $(H1)$ or $(H2)$, for the function $ U(t,x,[\theta]):=\mE[ g(X_t^{x,[\theta]},[X_t^{\theta}])]$, the derivative functions
	\ce
	(0,T]\times \mR^N \times \cP_2(\mR^N) &\ni& (t,x,[\theta]) \mapsto  \partial_xU(t,x,[\theta]),  \\
	(0,T]\times \mR^N \times \cP_2(\mR^N)\times \mR^N &\ni& (t,x,[\theta],v) \mapsto  \partial_\mu U(t,x,[\theta],\boldsymbol{v})
	\de
	exist and are  continuous. Moreover, for all compacts $K\subset \cP_2(\mR^N)$
	\ce
	\sup_{[\theta]\in K} \mE[|\partial_\mu U(t,x,[\theta],v)|^2]<\infty.
	\de
	\el
	
	\begin{proof}
		Since (H1) or (H2) holds,  $g$ belongs to the class $(IC)$ and there exists a sequence of functions $(g_l)_{l\geq 1}$ that approximates $g$. Let $U_l(t,x,[\theta])=\mE[ g_l( X_t^{x,[\theta]},[X_t^{\theta}])]$.  From Theorem \ref{IPF variable}, we can obtain that, for $\#\gamma=1$,
		\ce
		\partial_x^{i} U_l(t,x,[\theta])&=&\mE[ g_l(X_t^{x,[\theta]},[X_t^{\theta}])Z_{(i)}^1(1)(t,x,[\theta])].
		\de
		It is obvious that the expectation above is bounded uniformly in $l\geq 1$ by the growth assumption on $g_l$, H\"{o}lder's inequality
and the moment estimates already obtained for the processes $X_t^{x,[\theta]}$, $ X_t^{\theta}$ in Proposition \ref{eau},
and for $Z_{(i)}^1(1)(t,x,[\theta])$ in Lemma \ref{guocheng}.  So we can take the limit in each equation above by dominated convergence thoerem to conclude that
$ \partial_x U(t,x,[\theta])$  exists and
		 \ce
		 \partial_x^i U(t,x,[\theta])&=&\lim_{l\to \infty} \partial_x^\gamma U_l(t,x,[\theta])\\
		 &=&\lim_{l\to \infty}\mE[ g_l(X_t^{x,[\theta]},[X_t^{\theta}])Z_{(i)}^1(1)(t,x,[\theta])]\\
		 &=&\mE[ \lim_{l\to \infty} g_l(X_t^{x,[\theta]},[X_t^{\theta}])Z_{(i)}^1(1)(t,x,[\theta])]\\
		 &=&\mE[  g(X_t^{x,[\theta]},[X_t^{\theta}])Z_{(i)}^1(1)(t,x,[\theta])].
		 \de
Since  $Z_{(i)}^1(1)(t,x,[\theta])$ is jointly  continuous in $(t,x,[\theta])$ under the $L^p(\Omega)$,
by Theorem \ref{KS} and the  continuity of $g$, it follows that $ \partial_x U(t,x,[\theta])$ is  continuous. The existence and continuity of $\partial_x U(t,x,[\theta]) $ has now been established.

		To simplify notation, we consider the case $N=1$ through the rest of this proof. Firstly, we assume (H1) holds, so $g$ is in the class
(IC)$_x$. Note that $g_l$ satisfies the hypotheses of Lemma \ref{U}, which gives
		\be\label{5.2}
		\partial_\mu U_l(t,x,[\theta],v)&=&\mE[\partial g_l(X_t^{x,[\theta]},[X_t^{\theta}]) \partial_{\mu}X_t^{x,[\theta]}(v)  ]\nonumber\\
		&&+\mE\mE^{\prime}[ \partial_{\mu} g_l(X_t^{x,[\theta]},[X_t^{\theta}],(X_t^{v,[\theta]})^{\prime}) \partial_v (X_t^{v,[\theta]})^{\prime}\nonumber\\
		&&+\partial_{\mu} g_l(X_t^{x,[\theta]},[X_t^{\theta}],(X_t^{\theta^{\prime},[\theta]})^{\prime}) \partial_\mu (X_t^{\theta^{\prime},[\theta]})^{\prime}(v)]\nonumber\\
		&=&\mE[\partial g_l(X_t^{x,[\theta]},[X_t^{\theta}]) \partial_{\mu}X_t^{x,[\theta]}(v)  \nonumber]\\
		&&+\mE\mE^{\prime}[ \partial_{x} G_l(X_t^{x,[\theta]},[X_t^{\theta}],(X_t^{v,[\theta]})^{\prime})\partial_v (X_t^{v,[\theta]})^{\prime}\nonumber\\
		&&+\partial_{x} G_l(X_t^{x,[\theta]},[X_t^{\theta}],(X_t^{\theta^{\prime},[\theta]})^{\prime}) \partial_\mu (X_t^{\theta^{\prime},[\theta]})^{\prime}(v)].
		\ee
		Due to part $(1)$ of Theorem \ref{IPF measure}, the first term on the right-hand side of the last equality in (\ref{5.2})  can be written as
		\ce
		\mE[\partial g_l(X_t^{x,[\theta]},[X_t^{\theta}]) \partial_{\mu}X_t^{x,[\theta]}(v)]&=&\mE[\partial_\mu (g_l(X_t^{x,[\theta]},[X_t^{\theta}]))]\\
		&=&\mE[g_l(X_t^{x,[\theta]},[X_t^{\theta}]) Z_{\mu,(i)}^1(1)(t,x,[\theta],v) ].
		\de
		For the second term on the right-hand side of the last equality in (\ref{5.2}), by the identity
		\ce
		I_{N}=\widehat{\mE}[ DX_t^{x,[\theta]} (DX_t^{x,[\theta]})^{*} ] (\Gamma[X_t^{x,[\theta]}])^{-1},
		\de
		we have
		\ce
		&&\partial_{x} G_l(X_t^{x,[\theta]},[X_t^{\theta}],(X_t^{v,[\theta]})^{\prime}) \partial_v (X_t^{v,[\theta]})^{\prime}\\
		&=&\partial_{x} G_l(X_t^{x,[\theta]},[X_t^{\theta}],(X_t^{v,[\theta]})^{\prime}) \widehat{\mE}[ DX_t^{x,[\theta]}  (DX_t^{x,[\theta]})^{*} ] (\Gamma[X_t^{x,[\theta]}])^{-1} \partial_v (X_t^{v,[\theta]})^{\prime}\\
		&=& \widehat{\mE}[\partial_{x} G_l(X_t^{x,[\theta]},[X_t^{\theta}],(X_t^{v,[\theta]})^{\prime})  DX_t^{x,[\theta]}  (DX_t^{x,[\theta]})^{*} (\Gamma[X_t^{x,[\theta]}])^{-1}\partial_v (X_t^{v,[\theta]})^{\prime}]\\
		&=&\widehat{\mE}[ DG_l(X_t^{x,[\theta]},[X_t^{\theta}],(X_t^{v,[\theta]})^{\prime}) (DX_t^{x,[\theta]})^{*}  (\Gamma[X_t^{x,[\theta]}])^{-1} \partial_v (X_t^{v,[\theta]})^{\prime}].
		\de
Applying part $(2)$ of Theorem \ref{IPF variable}  we obtain
		\ce
		&&\mE\mE^{\prime}[ \partial_{x} G_l(X_t^{x,[\theta]},[X_t^{\theta}],(X_t^{v,[\theta]})^{\prime}) \partial_v (X_t^{v,[\theta]})^{\prime}]\\
		&=& \mE \mE^{\prime}\widehat{\mE}[ DG_l(X_t^{x,[\theta]},[X_t^{\theta}],(X_t^{v,[\theta]})^{\prime}) (DX_t^{x,[\theta]})^{*} (\Gamma[X_t^{x,[\theta]}])^{-1} \partial_v (X_t^{v,[\theta]})^{\prime}]\\
		&=& \mE \mE^{\prime}[ \widehat{\mE}[ DG_l(X_t^{x,[\theta]},[X_t^{\theta}],(X_t^{v,[\theta]})^{\prime}) (DX_t^{x,[\theta]})^{*}  (\Gamma[X_t^{x,[\theta]}])^{-1} ] \partial_v (X_t^{v,[\theta]})^{\prime}]\\
		&=& \mE \mE^{\prime}[ G_l(X_t^{x,      [\theta]},[X_t^{\theta}],(X_t^{v,[\theta]})^{\prime}) \delta((DX_t^{x,[\theta]})^{*} (\Gamma[X_t^{x,[\theta]}])^{-1}) \partial_v (X_t^{v,[\theta]})^{\prime}  ]\\
		&=&\mE \mE^{\prime}[ G_l(X_t^{x,[\theta]},[X_t^{\theta}],(X_t^{v,[\theta]})^{\prime}) Z_{(\gamma_i)}^2(1)(t,x,[\theta])  \partial_v (X_t^{v,[\theta]})^{\prime}].
		\de
		Similarly,
		\ce
		&&\mE\mE^{\prime}[\partial_{x} G_l(X_t^{x,[\theta]},[X_t^{\theta}],(X_t^{\theta^{\prime},[\theta]})^{\prime}) \partial_\mu (X_t^{\theta^{\prime},[\theta]})^{\prime}(v) ]\\
		&=& \mE\mE^{\prime}[ G_l(X_t^{x,[\theta]},[X_t^{\theta}],(X_t^{\theta^{\prime},[\theta]})^{\prime}) Z_{(\gamma_i)}^2(1)(t,x,[\theta])  \partial_\mu (X_t^{\theta^{\prime},[\theta]})^{\prime}(v)  ].
		\de
		Therefore, (\ref{5.2}) can be rewritten as
		\be\label{5.3}
		\partial_\mu U_l(t,x,[\theta],v)&=&\mE[g_l(X_t^{x,[\theta]},[X_t^{\theta}]) Z_{\mu,(i)}^1(1)(t,x,[\theta],v) ]\nonumber\\
		&&+\mE\mE^{\prime}[  G_l(X_t^{x,[\theta]},[X_t^{\theta}],(X_t^{v,[\theta]})^{\prime}) Z^2_{(i)}{(1)(t,x,[\theta])}  \partial_v (X_t^{v,[\theta]})^{\prime}\nonumber\\
		&&+G_l(X_t^{x,[\theta]},[X_t^{\theta}],(X_t^{\theta^{\prime},[\theta]})^{\prime}) Z^2_{(i)}(t,x,[\theta])  \partial_\mu (X_t^{\theta^{\prime},[\theta]})^{\prime}(v) ].
		\ee

		Since $g_l$, $ G_l$ satisfy the polynomial growth condition, and $\partial_v(X_t^{v,[\theta]})^{\prime}$ and $ \partial_\mu (X_t^{\theta^{\prime},[\theta]})^{\prime}(v)$ are bounded in $L^p$ by Lemma \ref{linear existence} and $Z_{{(i)}}^2(t,x,[\theta]), Z_{\mu,(i)}^{1}(t,x,[\theta],v)$ are  bounded in $L^p$ by Lemma \ref{guocheng}, it is obvious that all the term in (\ref{5.3}) can be bounded by $\|\theta\|_2$ except $(X_t^{\theta^{\prime},[\theta]})^{\prime} $ in the final term. Now we  calculate the last term of the right hands side of (\ref{5.3}):
		\ce
		&&|\mE\mE^{\prime}[ G_l(X_t^{x,[\theta]},[X_t^{\theta}],(X_t^{\theta^{\prime},[\theta]})^{\prime}) Z^2_{(i)}{(1)}  \partial_\mu (X_t^{\theta^{\prime},[\theta]})^{\prime}(v) ]|^2\\
		&\leq &\|G_l(X_t^{x,[\theta]},[X_t^{\theta}],(X_t^{\theta^{\prime},[\theta]})^{\prime}) \|_{L^{\frac 2 r}(\Omega\times \widetilde{\Omega})}\times \| Z^2_{(i)}{(1)}\|_{L^{\frac {4}{1-r}}(\Omega\times \widetilde{\Omega})}\times \|\partial_\mu (X_t^{\theta^{\prime},[\theta]})^{\prime}(v) \|_{L^{\frac{4}{1-r}}(\Omega\times \widetilde{\Omega})}\\
		&\leq & C\big( \mE\widetilde{\mE}\big[ (1+|X_t^{x,[\theta]}|^q+\| X_t^{\theta}\|_2^q+|(X_t^{\theta^{\prime},[\theta]})^{\prime} |^r)^{\frac 2 r}\big] \big)^r  \times \big( \mE \widetilde{\mE}[Z_{(i)}^2(1)] \big)^2 \times \big( \mE \widetilde{\mE}[\partial_\mu (X_t^{\theta^{\prime},[\theta]})^{\prime}(v) ] \big)^2\\
		&\leq & C\big( \mE\widetilde{\mE}\big[ (1+|X_t^{x,[\theta]}|^{2q/r}+\| X_t^{\theta}\|_2^{2q/r}+|(X_t^{\theta^{\prime},[\theta]})^{\prime} |^{2})\big] \big)^r  \times \big( \mE \widetilde{\mE}[Z_{(i)}^2(1)] \big)^2 \times \big( \mE \widetilde{\mE}[\partial_\mu (X_t^{\theta^{\prime},[\theta]})^{\prime}(v) ] \big)^2\\
		&\leq &C(1+|x|^{2q/r}+\|\theta\|_2^{2q/r}+\|\theta\|_2^{2})^r\times t^{-\frac{1}{a}}(1+|x|+\|\theta\|_2)^3,
		\de
		which gives
		\ce
		\sup_{[\theta]\in K} |\partial_\mu U_l(t,x,[\theta],\theta)|^2<\infty.
		\de
		Notice that  the estimate of (\ref{5.3}) are independent of $l$.  So we can take the limit in each equation above by
the dominated convergence theorem, and this shows that $\partial_\mu U(t,x,[\theta],v)$ exists and is continuous.
		
		Now, assume that (H2) holds, $g$ in the class (IC)$_v$. By Lemma \ref{U}, we can know the  expression for $\partial_\mu U_l$ and using the special form of $\partial_\mu g_l$ for initial conditions in the class (IC)$_v$,  and we get
		\be\label{5.6}
		\partial_\mu U_l(t,x,[\theta],v)&=&\mE[\partial g_l(X_t^{x,[\theta]},[X_t^{\theta}]) \partial_{\mu}X_t^{x,[\theta]}(v)  ]\nonumber\\
		&&+\mE\mE^{\prime}[ \partial_{v} G_l(X_t^{x,[\theta]},[X_t^{\theta}],(X_t^{v,[\theta]})^{\prime}) \partial_v (X_t^{v,[\theta]})^{\prime}\nonumber\\
		&&+\partial_{v} G_l(X_t^{x,[\theta]},[X_t^{\theta}],(X_t^{\theta^{\prime},[\theta]})^{\prime}) \partial_\mu (X_t^{\theta^{\prime},[\theta]})^{\prime}(v)].
		\ee
		The first term in (\ref{5.6}) can be written as
		\ce
		\mE[\partial g_l(X_t^{x,[\theta]},[X_t^{\theta}]) \partial_{\mu}X_t^{x,[\theta]}(v)]=\mE[\partial_\mu (g_l(X_t^{x,[\theta]},[X_t^{\theta}]))]\mE[g_l(X_t^{x,[\theta]},[X_t^{\theta}])Z_{\mu,(i)}^1(1)(t,x,[\theta],v) ].
		\de
		We again  use that
		\ce
		I_N=\widehat{\mE}[ D^{\prime}X_t^{x,[\theta]}  (D^{\prime}X_t^{x,[\theta]})^{*} ] (\Gamma[X_t^{x,[\theta]}])^{-1}
		\de
to calculate the second term in (\ref{5.6}):
		\ce
		&&\partial_{v} G_l(X_t^{x,[\theta]},[X_t^{\theta}],(X_t^{v,[\theta]})^{\prime}) \partial_v (X_t^{v,[\theta]})^{\prime}\\
		&=&\partial_{v} G_l(X_t^{x,[\theta]},[X_t^{\theta}],(X_t^{v,[\theta]})^{\prime}) \widehat{\mE}[ D^{\prime}X_t^{v,[\theta]} (D^{\prime}X_t^{v,[\theta]})^{\prime,*} ] (\Gamma[(X_t^{v,[\theta]})^{\prime}])^{-1} \partial_v (X_t^{v,[\theta]})^{\prime}\\
		&=&\widehat{\mE}[\partial_{v} G_l(X_t^{x,[\theta]},[X_t^{\theta}],(X_t^{v,[\theta]})^{\prime})  D^{\prime}X_t^{v,[\theta]} (D^{\prime}X_t^{v,[\theta]})^{\prime,*} (\Gamma[(X_t^{v,[\theta]})^{\prime}])^{-1} \partial_v (X_t^{v,[\theta]})^{\prime}]\\
		&=&\widehat{\mE}[ D^{\prime}G_l(X_t^{x,[\theta]},[X_t^{\theta}],(X_t^{v,[\theta]})^{\prime}) (D^{\prime}X_t^{v,[\theta]})^{\prime,*} (\Gamma[(X_t^{v,[\theta]})^{\prime}])^{-1} \partial_v (X_t^{v,[\theta]})^{\prime}].
		\de
		By part $(2)$ of Theorem \ref{IPF variable}  we have
		\ce
		&&\mE \mE^{\prime}[\partial_{v} G_l(X_t^{x,[\theta]},[X_t^{\theta}],(X_t^{v,[\theta]})^{\prime}) \partial_v (X_t^{v,[\theta]})^{\prime} ]\\
		&=& \mE \mE^{\prime}\widehat{\mE}[ D^{\prime}G_l(X_t^{x,[\theta]},[X_t^{\theta}],(X_t^{v,[\theta]})^{\prime}) (D^{\prime}X_t^{v,[\theta]})^{\prime,*} (\Gamma[(X_t^{v,[\theta]})^{\prime}])^{-1} \partial_v (X_t^{v,[\theta]})^{\prime}]\\
		&=& \mE \mE^{\prime}[ \partial_v (X_t^{v,[\theta]})^{\prime}\widehat{\mE}[ D^{\prime}G_l(X_t^{x,[\theta]},[X_t^{\theta}],(X_t^{v,[\theta]})^{\prime}) (D^{\prime}X_t^{v,[\theta]})^{\prime,*}  (\Gamma[(X_t^{v,[\theta]})^{\prime}])^{-1}]]\\
		&=& \mE \mE^{\prime}[ \partial_v (X_t^{v,[\theta]})^{\prime} G_l(X_t^{x,[\theta]},[X_t^{\theta}],(X_t^{v,[\theta]})^{\prime}) \delta ((D^{\prime}X_t^{v,[\theta]})^{\prime,*}  (\Gamma[(X_t^{v,[\theta]})^{\prime}])^{-1})]\\
		&=& \mE \mE^{\prime}[  G_l(X_t^{x,[\theta]},[X_t^{\theta}],(X_t^{v,[\theta]})^{\prime})  (Z_{(i)}^2(1)(t,x,[\theta],v))^{\prime} \partial_v (X_t^{v,[\theta]})^{\prime}].
		\de
But
		\ce
		&&\mE \mE^{\prime}[\partial_{v} G_l(X_t^{x,[\theta]},[X_t^{\theta}],(X_t^{\theta^{\prime},[\theta]})^{\prime}) \partial_\mu (X_t^{\theta^{\prime},[\theta]})^{\prime}(v) ]\\
		&=& \mE \mE^{\prime}[  G_l(X_t^{x,[\theta]},[X_t^{\theta}],(X_t^{\theta^{\prime},[\theta]})^{\prime}) (Z_{(i)}^2(1)(t,\theta^{\prime},[\theta],v))^{\prime} \partial_\mu (X_t^{\theta^{\prime},[\theta]})^{\prime}(v)],
		\de
		so (\ref{5.6}) can be rewritten as
		\ce
		\partial_{\mu} U_l(t,x,[\theta],v)&=&\mE[g_l(X_t^{x,[\theta]},[X_t^{\theta}]) Z_{\mu,(i)}^1(1) ]\\
		&&+ \mE \mE^{\prime}[  G_l(X_t^{x,[\theta]},[X_t^{\theta}],(X_t^{v,[\theta]})^{\prime})  (Z_{(i)}^2(1)(t,x,[\theta],v))^{\prime} \partial_v (X_t^{v,[\theta]})^{\prime}]\\
		&&+\mE \mE^{\prime}[  G_l(X_t^{x,[\theta]},[X_t^{\theta}],(X_t^{\theta^{\prime},[\theta]})^{\prime}) (Z_{(i)}^2(1)(t,\theta^{\prime},[\theta],v))^{\prime}\partial_\mu (X_t^{\theta^{\prime},[\theta]})^{\prime}(v)].
		\de
		
In the same way as the last step of the preceding proof under (H1), we can show from this formula that $ \sup_{[\theta]\in K} |\partial_\mu U_l(t,x,[\theta],\theta)|^2<\infty$ and the existence and the continuity of $\partial_\mu U(t,x,[\theta],\theta) $.
	\end{proof}

\begin{remark}
(1)
 Noticing that under hypothesis (H2), the coefficients are required to be uniformly bounded. The reason is that
Lemma \ref{guocheng} only provides a bound that grows like $|x|^3 $ for the $L^p$ norm and  substituting $x=\theta^{\prime}$ would then require controlling moments of $\theta^{\prime}$ up to order three, but $\theta^{\prime}$ is only assumed to be square-integrable. Even if such moments exist, additional growth conditions would complicate the argument. Therefore, assuming boundedness avoids these technical difficulties.

(2)
From the above proof, we can immediately get the following estimates
\ce
|\partial_x^i U(t,x,[\theta])|&\leq& Ct^{-\frac{2}{a}} (1+|x|+\|\theta\|_2)^{q+2},\\
|\partial_{\mu} U(t,x,[\theta])|&\leq& Ct^{-\frac{2}{a}} (1+|x|+\|\theta\|_2)^{q+2}.
\de

(3)
 Since $\partial_{\mu} U(t,x,[\theta])$ exists and is continuous, and for all compacts $K\subset \cP_2(\mR^N)$
	\ce
	\sup_{[\theta]\in K} \mE[|\partial_\mu U(t,x,[\theta],v)|^2]<\infty,
	\de
     we  deduce that $U$ is Lipschitz  continuous.

\end{remark}
		Now, we define what is the classical solution of PDE (\ref{pde1}).
	\bd
	Suppose that $U:[0,T]\times \mR^N \times \cP_2(\mR^N)$ satisfies PDE (\ref{pde1}) and
	\ce
	(0,T]\times \mR^N \times \cP_2(\mR^N) &\ni& (t,x,[\theta]) \mapsto  \partial_xU(t,x,[\theta]),\\
	(0,T]\times \mR^N \times \cP_2(\mR^N)\times \mR^N &\ni& (t,x,[\theta],v) \mapsto  \partial_\mu U(t,x,[\theta],v)
	\de
	exist and are continuous. Moreover, suppose that for all $(x,\theta)\in\mR^N \times L^2(\Omega)$
	\ce
	\lim_{(t,y,[\gamma])\to (0,x,[\theta])} U(t,y,[\gamma])=g(x,[\theta]).
	\de
	Then we say that $U$ is a solution to  PDE (\ref{pde1}).
	\ed
	
	Before presenting the main theorem of this section, we need the following chain rule for $(U(\mu_t))_{t\geq0}$ (cf. \cite[Theorem 7.2]{HL}), where $F$ is a smooth real-valued functional
	defined on the space $\mathcal{P}_2(\mathbb{R}^N)$, and $(\mu_t)_{t\geq0}$ denotes the flow of marginal measures associated with
	an $\mathbb{R}^N$-valued It\^o's process $(X_t)_{t\geq0}.$

	\bl\label{chainrule}
	Assume that $F\in\cC^{1,1}_{b,Lip}(\mR^N\times\cP_2(\mR^N);\mR^N)$ (here $F$ does not depend on the $x$ variable, i.e., $F(x, \mu) \equiv F(\mu)$ for all $x \in \mR^N$), then for any $t\geq 0$, we have
	\be\label{new chain rule}
	F(\mu_t (X_t) )&=&F(\mu_0(X_0))+ \int_0^t \mE[\partial_{\mu} F(\mu_s(X_s),(X_s)^{\prime}) b_s] ds\no\\
	&&+\int_0^t \int_{\Xi}\mE[ \big(F(\mu_s(X_{s-}+c(X_{s-},u)) ) -F(\mu_s(X_{s-}))\no\\
&&-  c(X_{s-},u) \partial_{\mu} [F( \mu_s(X_{s-},(X_s)^{\prime}))] \big)]\lambda(du)ds,
	\ee
	where $\mu_t(X_t)$ denotes the law of $X_t$ and $(X_s)^{\prime}$ is a copy of $X_s$.
	\el

	\bt\label{pdesolution}
	Suppose that either $(H1)$ or $(H2)$ holds. Then
	\ce
	U(t,x,[\theta]):=\mE[ g(X_t^{x,[\theta]},[X_t^\theta])]
	\de
	is a classical solution of  PDE (\ref{pde1}). Moreover, $U$ is unique among all of the classical solution satisfying the polynomial growth
	\ce
	|U(t,x,[\theta])|\leq C(1+|x|+\|\theta\|_2)^q
	\de
	for some $q>0$ and all ~$(t,x,[\theta]) \in [0,T]\times \mR^N \times \cP_2(\mR^N )$.
	\et
	\begin{proof}
		\textbf{Existence.} To prove continuity at the boundary, we use continuity of $g$ and the fact that
		\ce
		\| X_t^\theta-\theta\|_2+\| X_t^{x,[\theta]}-x\|_2 \to 0\quad \text{as}\quad t\to 0,
		\de
		which follows from (\ref{est3}).
		
		Now, note that by the flow property, we have, for ~$h>0$,
		\ce
		( X_{r,t+h}^{x,[\theta]},X_{r,t+h}^{\theta})=\big( X_{r,t}^{ X_{r,h}^{x,[\theta]},[ X_{r,h}^{\theta}]} ,X_{r,t}^{X_{r,h}^{\theta}}\big),\quad r\in[t,T].
		\de
		To lighten the notation, we omit the starting time $r$, so that
		\ce
		U(t+h,x,[\theta])&=&\mE[ g( X_{t+h}^{x,[\theta]},[X_{t+h}^\theta])]=\mE\big[g\big( X_t^{ X_h^{x,[\theta]},[ X_h^{\theta}]} ,[X_t^{X_h^{\theta}}]\big)\big]\\
		&=&\mE[\mE[g\big( X_t^{ X_h^{x,[\theta]},[ X_h^{\theta}]} ,[X_t^{X_h^{\theta}}]\big) ] \big| \cF_h]=\mE[ U(t,X_h^{x,[\theta]},[X_h^\theta])].
		\de
		Hence,
		\be\label{5.11}
		&&U(t+h,x,[\theta])-U(t,x,[\theta])\nonumber\\
		&=&\mE[ U(t,X_h^{x,[\theta]},[X_h^\theta])]-U(t,x,[\theta])\nonumber\\
		&=&\{ U(t,x,[X_h^\theta])-U(t,x,[\theta])\}+\mE[ U(t,X_h^{x,[\theta]},[X_h^\theta])-U(t,x,[X_h^{\theta}])].
		\ee
		Next, we take  different approaches to deal with the above two items and then drive by $h$ and let it to $0$. Together with continuity of the terms appearing in the expansion, we will prove that ~$U$ solves  PDE (\ref{pde1}).
		
		For the first term about measure, we can apply the chain rule introduced in Lemma \ref{chainrule} to the function $U(t,x,\cdot)$ to get
		\ce
		&&U(t,x,[X_h^\theta])-U(t,x,[\theta])\\
        &=& \int_0^h \mE[\partial_{\mu} U(t,x,[X_r^{\theta}],X_r^{\theta}) b(X_r^{x,[\theta]},[X_r^{\theta}])] dr\\
		&&+\int_0^h \int_{\Xi}\mE[ \big(U(t,x,[X_{r-}^{\theta}+c(X_{r-}^{x,[\theta]},u,[X_{r-}^{\theta}])] ) \\ &&-U(t,x,[X_{r-}^{\theta}])-c(X_{r-}^{x,[\theta]},u,[X_{r-}^{\theta}]) \partial_{\mu} [U( t,x,[X_{r-}^{\theta}])] \big)]\lambda(du)dr.
		\de
		
		 For the second term about $X_t^{x,[\theta]}$, It\^{o}'s formula with jumps (see \cite[Theorem  7.2]{HL})  can be applied to $U(t,\cdot,[X_h^\theta])$ to yield
		\be\label{final}
		&&U(t,X_h^{x,[\theta]},[X_h^\theta])-U(t,x,[X_h^{\theta}])\no\\
		&=& \int_0^h \partial_x U(t,X_{r}^{x,[\theta]},[X_r^\theta]) b(X_r^{x,[\theta]},[X_r^{\theta}]) dr\no\\
		&&+\int_0^h \int_{\Xi}( U(t,X_{r-}^{x,[\theta]}+c(X_{r-}^{x,[\theta]},u,[X_{r-}^\theta]),[X_{r}^{\theta}])- U(t,X_{r-}^{x,[\theta]},[X_r^\theta])\no\\
		&&-\partial_x U(t,X_{r-}^{x,[\theta]},[X_{r-}^\theta]) c(X_{r-}^{x,[\theta]},u, [X_{r-}^\theta]))\lambda(du)dr\no\\
		&&+\int_0^h \int_{\Xi}( U(t,X_{r-}^{x,[\theta]}+c(X_{r-}^{x,[\theta]},u,[X_{r-}^\theta]),[X_{r-}^{\theta}])-U(t,X_{r-}^{x,[\theta]},[X_{r-}^\theta]))\widetilde{N}(dr,du).\no\\
		\ee
		We next show that the final term in
		(\ref{final}) is square integrable, so that it is a true martingale with zero expectation. We have
		\ce
		&&U(t,X_{r-}^{x,[\theta]}+c(X_{r-}^{x,[\theta]},u,[X_{r-}^\theta]),[X_{r-}^{\theta}])- U(t,X_{r-}^{x,[\theta]},[X_{r-}^\theta])\\
		&=&\mE[g( X_r^{x,[\theta]}+c(X_r^{x,[\theta]},u,[X_s^{\theta}]),[X_r^{\theta}])]-\mE[g(X_r^{x,[\theta]},[X_r^{\theta}])]\\
		&\leq & |\mE[g( X_r^{x,[\theta]}+c(X_r^{x,[\theta]},u,[X_s^{\theta}]),[X_r^{\theta}])] |+|\mE[g(X_r^{x,[\theta]},[X_r^{\theta}])]|.
		\de
The first term on the RHS is dominated  by
		\ce
		&&\mE[g( X_r^{x,[\theta]}+c(X_r^{x,[\theta]},u,[X_r^{\theta}]),X_r^{\theta})]\\
		&\leq& C\| g( X_r^{x,[\theta]}+c(X_r^{x,[\theta]},u,[X_r^{\theta}]),X_r^{\theta})\|_2\\
		&\leq & C\| (1+|X_r^{x,[\theta]}+c(X_r^{x,[\theta]},u,[X_r^{\theta}]) | +\| X_r^{\theta}\|_2)^q\|_2\\
		&\leq & C\| (1+|X_r^{x,[\theta]}|+|c(X_r^{x,[\theta]},u,[X_r^{\theta}]) | +\| X_r^{\theta}\|_2)^q\|_2\\
		&\leq &C(1+|x|^q+\| \theta\|_2^q),
		\de
where	 the last inequity holds because of (\ref{estX2}), (\ref{estX3}) and Assumption \ref{H_solution} . So,
\ce
U(t,X_{r-}^{x,[\theta]}+c(X_{r-}^{x,[\theta]},u,[X_{r-}^\theta]),[X_{r-}^{\theta}])
\de
	  is square  integrable. Similarly, we can prove $\mE[g(X_r^{x,[\theta]},[X_r^{\theta}])] $ is also square  integrable. Therefore, the final term in (\ref{final}) is square integrable.

	 Putting the expansions back into (\ref{5.11}),
	 \ce
	 &&U(t+h,x,[\theta])-U(t,x,[\theta])\\
	 &=&\int_0^h \mE[\partial_{\mu} U(t,x,[X_r^{\theta}],X_r^{\theta}) b(X_r^{x,[\theta]},[X_r^{\theta}])] dr\\
	 &&+\int_0^h \int_{\Xi}\mE[ \big(U(t,x,[X_{r-}^{\theta}+c(X_{r-}^{x,[\theta]},u,[X_{r-}^{\theta}])] ) -U(t,x,[X_{r-}^{\theta}])\\ &&-c(X_{r-}^{x,[\theta]},u,[X_{r-}^{\theta}]) \partial_{\mu} [U( t,x,[X_{r-}^{\theta}])] \big)]\lambda(du)dr\\
	 &&+\int_0^h \mE[\partial_{x} U(t,X_r^{x,[\theta]},[X_r^{\theta}]) b(X_r^{x,[\theta]},[X_r^{\theta}])] dr\\
	 &&+\int_0^h \int_{\Xi}\mE[( U(t,X_{r-}^{x,[\theta]}+c(X_{r-}^{x,[\theta]},u,[X_{r-}^\theta]),[X_{r-}^\theta])- U(t,X_{r-}^{x,[\theta]},[X_{r-}^\theta])\\
	 &&-\partial_x U(t,X_{r-}^{x,[\theta]},[X_{r-}^\theta]) c(X_{r-}^{x,[\theta]},u, [X_{r-}^\theta]))]\lambda(du)dr,
	 \de
	  and dividing by $h$ and sending it to zero, we see that $U$ solves   PDE (\ref{pde1}).
	
	 \textbf{Uniqueness.}  We suppose $W$ is a solution with polynomial growth of  PDE. Let $\delta>0$, so
	 \ce
	 W(t,x,[\theta])-W(0,X_t^{x,[\theta]},[X_t^{\theta}])&=&W(t,x,[\theta])-W(\delta,X_{t-\delta}^{x,[\theta]},[X_{t-\delta}^{\theta} ])\\
&&+W(\delta,X_{t-\delta}^{x,[\theta]},[X_{t-\delta}^{\theta} ])-W(0,X_t^{x,[\theta]},[X_t^{\theta}]).
	 \de
	 Since $W$ is of polynomial growth, this is square integrable. Now applying It\^o's formula to
$(W(t-s,X_s^{x,[\theta]}, [X_s^{\theta}]))_{s\in [\delta,t]}$  and using the fact that $W$ is a solution of  PDE (\ref{pde1}),
we obtain
	 \ce
	 &&W(t,x,[\theta])-W(0,X_t^{x,\theta},[X_t^{\theta}])\\
	 &=& \int_\delta^t\int_{\Xi} [ W( r,X_{r-}^{x,[\theta]}+c(X_{r-}^{x,[\theta],u,[X_{r-}^{\theta}]}),[X_{r-}^{\theta}])-W( r,X_{r-}^{x,[\theta]},[X_{r-}^{\theta}])] \widetilde{N}(ds,du)\\
	 &&+W(\delta,X_{t-\delta}^{x,[\theta]},[X_{t-\delta}^{\theta} ])-W(0,X_t^{x,[\theta]},[X_t^{\theta}]).
	 \de
	
	 Since this is square integrable,  the stochastic integral is a true martingale. So taking expectation in the above expansion, we get
	 \ce
	 W(t,x,[\theta])-\mE[ W(0,X_t^{x,[\theta]},[X_t^{\theta}])]=\mE[W(\delta,X_{t-\delta}^{x,[\theta]},[X_{t-\delta}^{\theta} ]) ]-W(0,X_t^{x,[\theta]},[X_t^{\theta}]).
	 \de
	 Now, Let $\delta$ tend to $0$  and use continuity of $W$ at boundary, we can get
	 \ce
	 W(t,x,[\theta])=\mE[ W(0,X_t^{x,[\theta]},[X_t^{\theta}])]=\mE[g(X_t^{x,[\theta]},[X_t^{\theta}])],
	 \de
	 which finishes the proof.
	\end{proof}

	\section{Appendix}

	\subsection{Proofs of Theorem \ref{KS}}\label{appA}
	In order to prove Theorem \ref{KS}, the following two lemmas need to be introduced.

	\bl\label{linear existence}
	Let $Y_t^{x,[\theta]}(\boldsymbol{v})$ solve the following SDE with jumps
	\be\label{first}
	Y_t^{x,[\theta]}(\boldsymbol{v}) &=& a_0 + \int_0^t \{ a_1^1(s, x, [\theta]) Y_{s}^{x,[\theta]}(\boldsymbol{v}) + a_2^1(s, x, [\theta], \boldsymbol{v}) \no\\
	&&+  {\mathbb{E}^{\prime}}[ a_3^1(s, x, [\theta], {\theta^{\prime}}) ({Y}_{s}^{{\theta^{\prime}},[\theta]}(\boldsymbol{v}))^{\prime}
	+ \sum_{r=1}^{\#\boldsymbol{v}} a_3^1(s, x, [\theta], v_r) ({Y}_{s}^{v_r,[\theta]}(\boldsymbol{v}) )^{\prime}]\} ds\no\\
	&&+ \int_0^t\int_{\Xi} \{ a_1^2(s, x,u, [\theta]) Y_{s-}^{x,[\theta]}(\boldsymbol{v}) + a_2^2(s, x,u, [\theta], \boldsymbol{v}) \no\\
	&&+  {\mathbb{E}^{\prime}}[ a_3^2(s, x,u, [\theta], {\theta^{\prime}}) ({Y}_{s-}^{{\theta^{\prime}},[\theta]}(\boldsymbol{v}))^{\prime}
	+ \sum_{r=1}^{\#\boldsymbol{v}} a_3^2(s, x, u,[\theta], v_r) ({Y}_{s-}^{v_r,[\theta]}(\boldsymbol{v}) )^{\prime}]\} \widetilde{N}(ds,du),\no\\
	\ee
	where $v_r$ denotes an element in the $\boldsymbol{v}=(v_1,\cdots,v_{\#\boldsymbol{v}})$ and
\ce
	\begin{aligned}
		& a_{0}\in\mathbb{R}^{N}, \\
		& a_{1}^1:\Omega\times[0,T]\times\mathbb{R}^{N}\times\mathcal{P}_{2}(\mathbb{R}^{N})\rightarrow\mathbb{R}^{N\times N}, \\
		&
		 a_{1}^2:\Omega\times[0,T]\times\mathbb{R}^{N}\times\Xi\times\mathcal{P}_{2}(\mathbb{R}^{N})\rightarrow\mathbb{R}^{N\times N}, \\
		& a^1_{2}:\Omega\times[0,T]\times\mathbb{R}^{N}\times\mathcal{P}_{2}(\mathbb{R}^{N})\times(\mathbb{R}^{N})^{\#\boldsymbol{v}}\to\mathbb{R}^{N}, \\
		& a^2_{2}:\Omega\times[0,T]\times\mathbb{R}^{N}\times \Xi \times\mathcal{P}_{2}(\mathbb{R}^{N})\times(\mathbb{R}^{N})^{\#\boldsymbol{v}}\to\mathbb{R}^{N}, \\
		& a_{3}^1:{\Omega}^{\prime}\times\Omega\times[0,T]\times\mathbb{R}^{N}\times\mathcal{P}_{2}(\mathbb{R}^{N})\times\mathbb{R}^{N}\to\mathbb{R}^{N\times N},\\
		& a_{3}^2:{\Omega}^{\prime}\times\Omega\times[0,T]\times\mathbb{R}^{N}\times\Xi\times\mathcal{P}_{2}(\mathbb{R}^{N})\times\mathbb{R}^{N}\to\mathbb{R}^{N\times N}.
	\end{aligned}
	\de
  Suppose that for all $p\geq 2$ and $k=1,2,3$, the coefficients  $(t,x,[\theta],v)\mapsto a^1_k(t,x,[\theta],\boldsymbol{v})$ and $(t,x,u,[\theta],\boldsymbol{v})\mapsto a^2_k(t,x,u,[\theta],\boldsymbol{v})$  are   continuous in $L^p(\Omega)$.
	
	In (\ref{first}), $(Y^{\theta^{\prime},[\theta]})^{\prime}$ is a copy of ~$ Y^{x,[\theta]}$ on the probability space~$( \Omega^{\prime},\cF^{\prime},\mP^{\prime})$ driven by the Poisson process ~$N^{\prime}$ and with initial ~$x=\theta^{\prime}$.  Similarly, $(Y^{v_r,[\theta]})^{\prime}$ is a copy of ~
	$ Y^{x,[\theta]}$ on the probability space~$( \Omega^{\prime},\cF^{\prime},\mP^{\prime})$ driven by the Poisson process ~$N^{\prime}$ and with initial ~$x=v_r$. Furthermore, we make the following  assumptions:
	\begin{enumerate}
		\item $ a_k^2(t,x,u,[\theta],\boldsymbol{v})\leq CH(u) h_k(t,x,[\theta],\boldsymbol{v})$;
        \item $a_k(t,x,[\theta],\boldsymbol{v})=2\max(a_k^1(t,x,[\theta],\boldsymbol{v}),h_k(t,x,[\theta],\boldsymbol{v}))$;
		\item $\sup_{x\in\mathbb{R}^N,[\theta]\in\mathcal{P}_2(\mathbb{R}^N),\boldsymbol{v}\in(\mathbb{R}^N)^{\# v}}\|a_2(\cdot,x,[\theta],\boldsymbol{v})\|_{\mathcal{S}_T^p}<\infty$;
		\item $a_1$ and ~$a_3$ are uniformly bounded;
		\item $ \sup_{x\in\mathbb{R}^N,[\theta]\in\mathcal{P}_2(\mathbb{R}^N),\boldsymbol{v}\in(\mathbb{R}^N)^{\#\boldsymbol{v}}}\|{a}_2(\cdot,\theta,[\theta],\boldsymbol{v})\|_{\mathcal{S}_T^2}<\infty,$
	\end{enumerate}
    where $(t,x,[\theta],\boldsymbol{v})\mapsto h_k(t,x,[\theta],\boldsymbol{v})$ is continuous in $L^p(\Omega)$ and $H(u)\in L^2(\mP)$.
	Then, there is a constant $C=C(p,T,a_1,a_3)>0$ such that
	\be\label{guji}
	\| Y^{x,[\theta]}\|_{S_T^p} \leq C\big(|a_0|+  \| a_2(\cdot,x,[\theta],\boldsymbol{v})\|_{S_T^2}+ \| a_2(\cdot,x,[\theta],\boldsymbol{v})\|_{S_T^p}\big).
	\ee
	Moreover, the mapping
	\ce
	[0,T]\ni t \mapsto Y_t^{x,[\theta]}(\boldsymbol{v})\in L^p(\Omega)
	\de
	is c\`adl\`ag and the mapping
    \ce
    \mathbb{R}^N\times\mathcal{P}_2(\mathbb{R}^N)\times(\mathbb{R}^N)^{\#\boldsymbol{v}} \ni (x,[\theta],\upsilon) \mapsto Y_t^{x,[\theta]}(\boldsymbol{v})\in L^p(\Omega)
    \de
    is continuous.
	\el
	
	\begin{proof}
		Let ~$ \iota,\kappa:[0,T]\mapsto [0,\infty)$ be defined as
		\ce
		\begin{aligned}
			& \iota\left(t\right)=\left\|Y^{{\theta}^{\prime},[\theta]}(\boldsymbol{v})\right\|_{\mathcal{S}_{t}^{2}}^{2}+\sum_{r=1}^{\#\boldsymbol{v}}\left\|Y^{v_{r},[\theta]}(\boldsymbol{v})\right\|_{\mathcal{S}_{t}^{2}}^{2},t\in[0,T], \\
			& \kappa\left(t\right)=\left\|Y^{x,[\theta]}(\boldsymbol{v})\right\|_{S_{t}^{p}}^{p},t\in[0,T].
		\end{aligned}
		\de
        By Kunita's first inequality, H\"{o}lder's inequality and the assumptions on coefficients, we conclude that there exists a constant $C>0$ such that
        \ce
        \iota(t) \leq C\bigl(|a_0|^2 + \|a_2\|_{S_t^2}^2 + (\|a_1\|_{\infty}^2 + \|a_3\|_{\infty}^2) \int_0^t \iota(h) \, dh\bigr).
        \de
        By Gr\"{o}nwall's inequality,
        \ce
        \iota(t) &\leq& C\Bigl(|a_0|^2 + \|a_2\|_{S_t^2}^2 + \int_0^t (\|a_1\|_{\infty}^2 + \|a_3\|_{\infty}^2) \iota(h) \, dh\Bigr)\\
        &\leq& C\bigl(|a_0|^2 + \|a_2\|_{S_t^2}^2\bigr) e^{\int_0^t (\|a_1\|_{\infty}^2 + \|a_3\|_{\infty}^2) \, dh}\\
        &\leq& C\bigl(|a_0|^2 + \|a_2\|_{S_t^2}^2\bigr) e^{T \cdot (\|a_1\|_{\infty}^2 + \|a_3\|_{\infty}^2)}.
        \de
		For $\kappa(t)$, we deduce
		\ce
		\kappa(t) &\leq& C(|a_0|^p+  \|a_3\|_{\infty}^p \iota(T)^{ \frac p  2} + \| a_2(\cdot,x,[\theta],v)\|_{S_T^p}^p +\int_0^t  \|a_{1}(\cdot,x,[\theta])\|_{\infty}^p \kappa(s)dh ),
		\de
which gives by Gr\"onwall's inequality
		\ce
		\kappa(t) &\leq & C(|a_0|^p +   \|a_3\|_{\infty}^p \iota(T)^{ \frac p  2} +\| a_2(\cdot,x,[\theta],v)\|_{S_T^p}^p )  e^{ \int_0^t  \|a_{1}(\cdot,x,[\theta])\|_{\infty}^p dh )}.
		\de
Since ~$a_1$ and~$a_3$ are informally bounded, we have
		\ce
		\kappa(t) &\leq & C( |a_0|^p +\iota(T)^{ \frac p  2} + \| a_2(\cdot,x,[\theta],\boldsymbol{v})\|_{S_T^p}^p ).
		\de
Using our estimate on ~$\iota(t)$ we get
		\ce
		\kappa(t) &\leq & C( |a_0|^p +\| a_2(\cdot,x,[\theta],\boldsymbol{v})\|_{S_T^2}^p + \| a_2(\cdot,x,[\theta],\boldsymbol{v})\|_{S_T^p}^p ) .
		\de
		So we have
		\ce
		\| Y^{x,[\theta]}\|_{S_T^p} \leq C(|a_0|+  \| a_2(\cdot,x,[\theta],\boldsymbol{v})\|_{S_T^2}+ \| a_2(\cdot,x,[\theta],\boldsymbol{v})\|_{S_T^p}).
		\de
		Now we turn to prove the continuity of $Y_t^{x,[\theta]}$. We use the following decomposition
		\ce
			&&Y_t^{x,[\theta]}(\boldsymbol{v})-Y_{t^{\prime}}^{x^{\prime},[\theta^{\prime}]}(\boldsymbol{v}^{\prime}) \\
			&=& Y_t^{x,[\theta]}(\boldsymbol{v})-Y_{t^{\prime}}^{x,[\theta]}(\boldsymbol{v})\\
&&+Y_{t^{\prime}}^{x,[\theta]}(\boldsymbol{v})-Y_{t^{\prime}}^{x^{\prime},[\theta]}(\boldsymbol{v})+Y_{t^{\prime}}^{x^{\prime},[\theta]}(\boldsymbol{v})-Y_{t^{\prime}}^{x^{\prime},[\theta^{\prime}]}(\boldsymbol{v})+Y_{t^{\prime}}^{x^{\prime},[\theta^{\prime}]}(\boldsymbol{v})-Y_{t^{\prime}}^{x^{\prime},[\theta^{\prime}]}(\boldsymbol{v}^{\prime})\\
			&=:&\Delta_tY^{x,[\theta]}(\boldsymbol{v})+\Delta_xY_{t^{\prime}}^{\theta}(\boldsymbol{v})+\Delta_\theta Y_{t^{\prime}}^{x^{\prime}}(\boldsymbol{v})+\Delta_{\boldsymbol{v}}Y_{t^{\prime}}^{x^{\prime},[\theta^{\prime}]},
		\de
		and deal with the RHS term by term. Firstly,
		\ce
		\Delta_tY^{x,[\theta]}(\boldsymbol{v})&=&\int_{t^{\prime}}^t\{a_1^1Y_s^{x,[\theta]}(\boldsymbol{v})+a_2^1+{\mathbb{E}^{\prime}}[a_3^1|_{v={\theta^{\prime}}}(Y_s^{{\theta^{\prime}},[\theta]})^{\prime}(\boldsymbol{v})+\sum_{r=1}^{\#\boldsymbol{v}}a_3^1|_{v=v_r}({Y}_s^{v_r,[\theta]})^{\prime}(\boldsymbol{v})]\}ds\\
		&&+\int_{t^{\prime}}^t\int_{\Xi}\{a_1^2Y_s^{x,[\theta]}(\boldsymbol{v})+a_2^2+{\mathbb{E}^{\prime}}[a_3^2|_{v={\theta^{\prime}}}(Y_s^{{\theta^{\prime}},[\theta]})^{\prime}(\boldsymbol{v})\\
&&+\sum_{r=1}^{\#\boldsymbol{v}}a_3^2|_{v=v_r}({Y}_s^{v_r,[\theta]})^{\prime}(\boldsymbol{v})]\}\widetilde{N}(ds,du).
		\de
		Since the integrand is uniformly bounded in $L^p (\Omega)$ in time,  we have by  Kunita's second inequality,
		\ce
		\left\|\Delta_tY^{x,[\theta]}(\boldsymbol{v})\right\|_p\leq C|t-t^{\prime}|^{\frac{1}{2}} ,
		\de
		 which goes to ~$0$ as ~$t\to t^{'}$.
		
		Secondly,
		\ce
			\Delta_{x}Y_{t}^{\theta}(\boldsymbol{v}) & =&\Delta_{x}a_{0}+\int_{0}^{t}\{a_{1}^1\Delta_{x}Y_{s}^{[\theta]}(\boldsymbol{v})+Y_{s}^{x,[\theta]}(\boldsymbol{v})\Delta_{x}a_{1}^1+\Delta_{x}a_{2}^1 \\
			&&+{\mathbb{E}^{'}}[\Delta_{x}a_{3}^1|_{v={\theta^{'}}}({Y}_{s}^{{\theta^{'}},[\theta]}(\boldsymbol{v}))^{\prime}+\sum_{r=1}^{\#\boldsymbol{v}}\Delta_{x}a_{3}^1|_{v=v_{r}}({Y}_{s}^{v_{r},[\theta]}(\boldsymbol{v}))^{\prime}]\}ds\\
			&&+\int_{0}^{t}\{a_{1}^2\Delta_{x}Y_{s}^{[\theta]}(\boldsymbol{v})+Y_{s}^{x,[\theta]}(\boldsymbol{v})\Delta_{x}a_{2}^1+\Delta_{x}a_{2}^2 \\
			&&+{\mathbb{E}^{'}}[\Delta_{x}a_{3}^2|_{v={\theta^{'}}}({Y}_{s}^{{\theta^{'}},[\theta]}(\boldsymbol{v}))^{\prime}+\sum_{r=1}^{\#\boldsymbol{v}}\Delta_{x}a_{3}^2|_{v=v_{r}}({Y}_{s}^{v_{r},[\theta]}(\boldsymbol{v}))^{\prime}]\}\widetilde{N}(ds,du) .
		\de
		An identical line of reasoning to that which yielded the result (\ref{guji}) shows that
		\ce
		\begin{aligned}
			\left\|\Delta_{x}Y_{t}^{\theta}(\boldsymbol{v})\right\|_{p}^{p}\leq &  C\operatorname*{sup}_{s\in[0,t]}\mathbb{E}(Y_{s}^{x,[\theta]}(\boldsymbol{v})\Delta_{x}a_{1}+\Delta_{x}a_{2}+{\mathbb{E}^{\prime}}[({Y}_{s}^{{\theta^{\prime}},[\theta]}(\boldsymbol{v}))^{\prime}\Delta_{x}a_{3}|_{v={\theta^{\prime}}}\\ &+\sum_{r=1}^{\#\boldsymbol{v}}\Delta_{x}a_{3}|_{v=v_{r}}({Y}_{s}^{v_{r},[\theta]}(\boldsymbol{v}))^{\prime}])^{p}.
		\end{aligned}
		\de
		By H\"{o}lder's inequality, the boundedness of ~$Y_s^{x,[\theta]}(\boldsymbol{v})$ in ~$L^p(\Omega)$, and the continuity assumptions on ~$a_1, a_2, a_3$, we can get the above term goes to ~$0$ when $x\to x^{\prime}$.
		
		 The arguments for ~$ \partial_\theta Y_{t^{'}}^{x^{'}}(\boldsymbol{v} )$ and ~$\partial_{\boldsymbol{v}}Y_{t^{'}}^{x^{'},\theta^{'}}(\boldsymbol{v} ) $ are same.
	\end{proof}
	
	We  now consider the differentiability of $Y_t^{x,[\theta]}(v)$ in Lemma \ref{linear existence}. For simplification, we omit the ~$(t,x,u,[\theta])$ in ~$a_k$, and write ~$a_k|_{v=\theta^{\prime}}$ to denote ~$ a_k(s,x,u,[\theta],\theta^{\prime})$.

	\bl\label{linear second order}
	Suppose that the process~$Y^{x,[\theta]}$ is as in Lemma \ref{linear existence} and satisfies the assumptions in Lemma \ref{linear existence}. Here, we introduce the following additional  assumptions:
	\begin{enumerate}[label=(\alph*)]
		\item For $k=1,2,3$, all $(s,u,[\theta],\boldsymbol{v}) \in [0,T]\times\Xi\times \cP_2(\mR^N)\times(\mR^N)^{\#\boldsymbol{v}} $ and each~$ p\geq 1$, $ \mR^N \ni x \mapsto a_k^1(s, x,[\theta],\boldsymbol{v}) \in L^p(\Omega)$ and $ \mR^N \ni x \mapsto a_k^2(s, x,u,[\theta],\boldsymbol{v}) \in L^p(\Omega)$ are differentiable.
		\item For $k=1,2,3$, all $(s,x,u,[\theta]) \in [0,T]\times\mR^N\times\Xi\times \cP_2(\mR^N) $ and each~$ p\geq 1$, $ (\mR^N)^{\# \boldsymbol{v}} \ni v \mapsto a_k^1(s, x,[\theta],\boldsymbol{v}) \in L^p(\Omega)$ and $ (\mR^N)^{\# v} \ni \boldsymbol{v} \mapsto a_k^2(s, x,u,[\theta],\boldsymbol{v}) \in L^p(\Omega)$ are differentiable;
		\item for all $(s,x,u,\boldsymbol{v}) \in [0,T]\times \mR^N\times\Xi\times(\mR^N)^{\#\boldsymbol{v}} $ the mapping ~$ L^2(\cF_0;\mR^N) \ni \theta \mapsto a_k^1(s, x,[\theta],\boldsymbol{v}) \in L^p(\Omega)$ and ~$ L^2(\cF_0;\mR^N) \ni \theta \mapsto a_k^2(s, x,u,[\theta],\boldsymbol{v}) \in L^p(\Omega)$ are Fr\'{e}chet differentiable;
		\item for ~$k=1,2,3$ and all ~$(s,x,[\theta],\boldsymbol{v}) \in [0,T]\times \mR^N\times \cP_2(\mR^N)\times(\mR^N)^{\#\boldsymbol{v}} ,$ $\Xi \ni u\mapsto a_k^2(s,x,u,[\theta],\boldsymbol{v}) \in \cH_{\overline{\mD}^{1,\infty}}^N $.  Moreover, we assume the following estimates on the Malliavin  derivatives holds
		\ce
		\mathbb{E} [\sup_{s \in [0,T]}|D a_k^1(s,x,[\theta],\boldsymbol{v})|^p]+\mathbb{E} [\sup_{s \in [0,T]}|D a_k^2(s,x,u,[\theta],\boldsymbol{v})|^p]<\infty, \quad k=0,1,2,3.
		\de
	\end{enumerate}
	Then, for all~$t\in [0,T] $ the following statements hold true:
	\begin{enumerate}
		\item Under Assumption $(a)$, $ x\mapsto Y_t^{x,[\theta]}(\boldsymbol{v})$ is differentiable in ~$L^p(\Omega)$ for all ~$ p\geq 1$, i.e., there exists a process $\partial_{x}Y_{t}^{x,[\theta]}(\boldsymbol{v}) $ such that
		\ce
		\mE[\sup_{s\in [t,T]}|Y_s^{x+h,[\theta]}(\boldsymbol{v})-Y_s^{x,[\theta]}(\boldsymbol{v})-\partial_{x}Y_{t}^{x,[\theta]}(\boldsymbol{v})|^p]=o(|h|^p),\quad \mR^N \ni h \to 0.
		\de
		Moreover, $\partial_{x}Y_{t}^{x,[\theta]}(\boldsymbol{v}) $ is the unique solution of the following SDE with jumps:
		\ce
			\partial_{x}Y_{t}^{x,[\theta]}(\boldsymbol{v}) &=&\int_{0}^{t}\big\{\partial_{x}a^1_{1}Y_{s}^{x,[\theta]}(\boldsymbol{v})+a_{1}^1\partial_{x}Y_{s}^{x,[\theta]}(\boldsymbol{v})+\partial_{x}a_{2}^1 \\
			&&+{\mathbb{E}^{\prime}}[\partial_{x}a_{3}^1|_{v=\theta^{'}}(Y_{s}^{\theta^{'},[\theta]}(\boldsymbol{v}))^{'}+\sum_{r=1}^{\#\boldsymbol{v}}\partial_{x}a_{3}^1|_{\upsilon=v_{r}}(Y_{s}^{v_r,[\theta]}(\boldsymbol{v}))^{\prime}] \big\}ds\\
			&&+\int_{0}^{t}\int_{\Xi}\big\{\partial_{x}a_{1}^2Y_{s}^{x,[\theta]}(\boldsymbol{v})+a_{1}^2\partial_{x}Y_{s}^{x,[\theta]}(\boldsymbol{v})+\partial_{x}a_{2}^2 \\
			&&+{\mathbb{E}^{\prime}}[\partial_{x}a_{3}^2|_{v=\theta^{'}}(Y_{s}^{\theta^{'},[\theta]}(\boldsymbol{v}))^{'}+\sum_{r=1}^{\#\boldsymbol{v}}\partial_{x}a_{3}^2|_{\upsilon=v_{r}}(Y_{s}^{v_r,[\theta]}(\boldsymbol{v}))^{\prime}] \big\}\widetilde{N}(ds,du).
		\de
		\item Under Assumption $(b)$, $ v\mapsto Y_t^{x,[\theta]}(\boldsymbol{v})$ is differentiable in $L^p(\Omega)$ for all ~$p\geq 1$, i.e., there exists a process $\partial_{\boldsymbol{v}}Y_{t}^{x,[\theta]}(\boldsymbol{v}) $ such that
		\ce
		\mE[\sup_{s\in [t,T]}|Y_s^{x,[\theta]}(\boldsymbol{v}+h)-Y_s^{x,[\theta]}(\boldsymbol{v})-\partial_{v}Y_{t}^{x,[\theta]}(\boldsymbol{v})|^p]=o(|h|^p),\quad \mR^N \ni h \to 0.
		\de
		Moreover, $\partial_{v_j}Y_{t}^{x,[\theta]}(\boldsymbol{v}) $ is the unique solution of the following SDE with jumps:
		\ce
			\begin{aligned}
				\partial_{v_j} Y_t^{x,[\theta]}(\boldsymbol{v})
			&	=  \int_0^t \Big\{ \,
				a_1^1 \partial_{v_j} Y_s^{x,[\theta]}(\boldsymbol{v})
				+ \partial_{v_j} a_2^1+\mathbb{E}^\prime \big[ a_3^1 \big|_{v=\theta^\prime} \partial_{v_j} (Y_s^{\theta^\prime,[\theta]}(\boldsymbol{v}))^\prime \big]  \\
				& + \mathbb{E}^\prime \big[ a_3^1 \big|_{v=v_j} \partial_x (Y_s^{v_j,[\theta]}(\boldsymbol{v}))^\prime \big] +\mathbb{E}^\prime \big[ \partial_{v_j} a_3^1 \big|_{v=v_j} (Y_s^{v_j,[\theta]}(\boldsymbol{v}))^\prime \big]  \\
				& + \sum_{r=1,r\neq j}^{\# \boldsymbol{v}} \mathbb{E}^\prime \big[ a_3^1 \big|_{v=v_r} \partial_{v_j} (Y_s^{v_r,[\theta]}(\boldsymbol{v}))^\prime \big] \Big\}
				\, ds\\
				 & +\int_0^t \int_{\Xi}\Big\{ \,
				a_1^2 \partial_{v_j} Y_s^{x,[\theta]}(\boldsymbol{v}) +\partial_{v_j} a_2^2
				+ \mathbb{E}^\prime \big[ a_3^2 \big|_{v=\theta^\prime} \partial_{v_j} (Y_s^{\theta^\prime,[\theta]}(\boldsymbol{v}))^\prime \big] \\
				& + \mathbb{E}^\prime \big[ a_3^2 \big|_{v=v_j} \partial_x (Y_s^{v_j,[\theta]}(\boldsymbol{v}))^\prime \big] + \mathbb{E}^\prime \big[ \partial_{v_j} a_3^2 \big|_{v=v_j} (Y_s^{v_j,[\theta]}(\boldsymbol{v}))^\prime \big]  \\
				& + \sum_{r=1,r\neq j}^{\# \boldsymbol{v}} \mathbb{E}^\prime \big[ a_3^2 \big|_{v=v_r} \partial_{v_j} (Y_s^{v_r,[\theta]}(\boldsymbol{v}))^\prime \big] \Big\}
				\, \widetilde{N}(ds, du).
			\end{aligned}	
		\de
		\item Under Assumptions $(a),(b)$ and $(c)$, the mapping ~$L^2(\cF_0;\mR^N) \ni\theta\mapsto Y_{t}^{x,[\theta]}(\boldsymbol{v})$ is Fr\'{e}chet differentiable for all ~$(x,\boldsymbol{v})\in\mathbb{R}^N\times(\mathbb{R}^N)^{\#\boldsymbol{v}}, $ so ~$\partial_\mu Y_t^{x,[\theta]}(\boldsymbol{v},{v}^{\prime}) $ exists and it satisfies
		\ce
		\partial_\mu Y_t^{x,[\theta]}(\boldsymbol{v},{v}^{\prime}) &
	    =&\int_0^t\{\partial_\mu a_1^1 Y_s^{x,[\theta]}(\boldsymbol{v})+a_1^1 \partial_\mu Y_s^{x,[\theta]}(\boldsymbol{v},v^{\prime})+\partial_\mu a_2^1 \\
		&&+\mE^{\prime}\bigg[\partial_\mu a_3^1 ({Y}_s^{\theta^{\prime},[\theta]}(\boldsymbol{v}))^{\prime}+\partial_{\boldsymbol{v}} a_3^1 ({Y}_s^{\upsilon^{\prime},[\theta]}(\boldsymbol{v}) )^{\prime}\\ &&+a_3^1|_{v=\theta^{\prime}}\partial_\mu({Y}_s^{\theta^{\prime},[\theta]}(\boldsymbol{v},v^{\prime}))^{\prime}]+a_3^1|_{v=\upsilon^{\prime}}\partial_{x}({Y}_{s}^{v^{\prime},[\theta]}(\boldsymbol{v}) )^{\prime}\\
		&&+\mE^{\prime}[\sum_{r=1}^{\# \boldsymbol{v}}a_3^1|_{v=v_r}\partial_\mu({Y}_s^{\upsilon_r,[\theta]}(\boldsymbol{v},v^{\prime}))^{\prime}+ \sum_{r=1}^{\# \boldsymbol{v}} \partial_{\mu} a_3^1 ({Y}_s^{v_r,[\theta]}(\boldsymbol{v},v^{\prime}))^{\prime}]\}ds\\
		&&+\int_0^t\int_{\Xi}\{\partial_\mu a_1^2 Y_s^{x,[\theta]}(\boldsymbol{v})+a_1^2 \partial_\mu Y_s^{x,[\theta]}(\boldsymbol{v},v^{\prime})+\partial_\mu a_2^2\\
		&&+\mE^{\prime}\bigg[\partial_\mu a_3^2 ({Y}_s^{\theta^{\prime},[\theta]}(\boldsymbol{v}))^{\prime}+\partial_{\boldsymbol{v}} a_3^2 ({Y}_s^{\upsilon^{\prime},[\theta]}(\boldsymbol{v}) )^{\prime}\\ &&+a_3^2|_{v=\theta^{\prime}}\partial_\mu({Y}_s^{\theta^{\prime},[\theta]}(\boldsymbol{v},v^{\prime}))^{\prime}+a_3^2|_{v=\upsilon^{\prime}}\partial_{x}({Y}_{s}^{\upsilon^{\prime},[\theta]}(\boldsymbol{v}) )^{\prime}] \\
		&&+\mE^{\prime}[ \sum_{r=1}^{\#\boldsymbol{v}}a_3^2|_{v=v_r}\partial_\mu({Y}_s^{\upsilon_r,[\theta]}(\boldsymbol{v},v^{\prime}))^{\prime}+ \sum_{r=1}^{\# \boldsymbol{v}} \partial_{\mu} a_3^2 ({Y}_s^{\upsilon_r,[\theta]}(\boldsymbol{v},v^{\prime}))^{\prime}]\}\widetilde{N}(ds,du).
		\de
		Moreover, the Fr\'echet derivative $\cD(Y_t^{\theta,[\theta]})(\gamma)$  has the following representation, for all ~$\gamma \in L^2(\cF_0;\mR^N)$,
		\ce
		\cD\left(Y_{t}^{\theta,[\theta]}(\boldsymbol{v})\right)(\gamma)=\left.\left(\partial_{x}Y_{t}^{x,[\theta]}(\boldsymbol{v})\gamma+{\mathbb{E}^{\prime}}\left[\partial_{\mu}Y_{t}^{x,[\theta]}(\boldsymbol{v},{\theta^{\prime}}){\gamma}^{\prime}\right]\right)\right|_{x=\theta}.
		\de
		\item Under Assumption \ref{R}, $ Y_t^{x,[\theta]} \in \cH_{\overline{\mD}^{1,\infty}}^N$ and its Malliavin derivative  satisfies the following SDE with jumps:
		\ce
		DY^{x,[\theta]}_t
		&=&   \int_0^t Da_1^1(s, x, [\theta]) Y_{s-}^{x,[\theta]}(\boldsymbol{v}) +a_1^1(s, x, [\theta])DY_{s-}^{x,[\theta]}(\boldsymbol{v}) + Da_2^1(s, x, [\theta], \boldsymbol{v}) \\
		&&+  {\mathbb{E}^{\prime}}[ Da_3^1(s, x, [\theta], {\theta^{\prime}}) ({Y}_{s-}^{{\theta^{\prime}},[\theta]}(\boldsymbol{v}))^{\prime}
		+ \sum_{r=1}^{\#v} Da_3^1(s, x,[\theta], v_r) ({Y}_{s-}^{v_r,[\theta]}(\boldsymbol{v}) )^{\prime}] ds\\
		&&+   \int_0^t\int_{\Xi} Da_1^2(s, x,u, [\theta]) Y_{s-}^{x,[\theta]}(\boldsymbol{v}) +a_1^2(s, x,u, [\theta]) DY_{s-}^{x,[\theta]}(\boldsymbol{v}) + Da_2^2(s, x,u, [\theta], \boldsymbol{v}) \\
		&&+  {\mathbb{E}^{\prime}}[ Da_3^2(s, x,u, [\theta], {\theta^{\prime}}) ({Y}_{s-}^{{\theta^{\prime}},[\theta]}(\boldsymbol{v}))^{\prime}
		+ \sum_{r=1}^{\#\boldsymbol{v}} Da_3^2(s, x, u,[\theta], v_r) ({Y}_{s-}^{v_r,[\theta]}(\boldsymbol{v}) )^{\prime}] \widetilde{N}(ds,du)\\
		&&+\int_0^t\int_{\Xi\times R} \{ a_1^2(s, x,u, [\theta]) Y_{s-}^{x,[\theta]}(\boldsymbol{v},r) + a_2^2(s, x,u, [\theta], \boldsymbol{v},r) \\
		&&+  {\mathbb{E}^{\prime}}[ a_3^2(s, x,u, [\theta], {\theta^{\prime}}) ({Y}_{s-}^{{\theta^{\prime}},[\theta]}(\boldsymbol{v},r))^{\prime}\\
		&&+ \sum_{r=1}^{\#\boldsymbol{v}} a_3^2(s, x, u,[\theta], v_r) ({Y}_{s-}^{v_r,[\theta]}(\boldsymbol{v},r) )^{\prime}]\}^{\flat} N\odot \rho (ds,du,dr).
		\de
		Moreover, the following bound holds:
		\be\label{6.11}
		\mathbb{E}[\sup_{s\in[0,T]}|DY_s^{x,[\theta]}(\boldsymbol{v})|^p]\leq C\mathbb{E}[\sup_{s\in[0,T]}|D a_1|^p].
		\ee
	\end{enumerate}
	\el
	\begin{proof}
		Part $(1)$ and $(2)$ are results on the differentiability of SDEs with jumps with respect to a real parameter, so we omit them here.
			Part $(3)$ is a straight result from \cite[Proposition 5.1]{HL}.
			Finally, since $(Y^{\theta^{\prime},[\theta]}_t)^{\prime} $ and $(Y^{v_r,[\theta]}_t)^{\prime} $ do not depend on
$\omega\in\Omega$,  the Malliavin differentiability of the solution to Equation (\ref{first}) follows by a standard application of the framework in \cite[Proposition 8.25, Proposition 8.26]{BL}, with the corresponding estimate derived similarly to the proof of (\ref{guji}).
	\end{proof}
	 We now proceed to establish the proof of Theorem \ref{KS}. For the sake of notational simplicity, we will demonstrate this result in the one-dimensional setting.
	 \begin{proof}
	 	We will show via induction that the following statements hold for  ~$I=1,\cdots,k$:
	 	
	 	$(S1)$ For all ~$\alpha,\boldsymbol{\beta},\gamma$ satisfying ~$ \#\alpha+\# \boldsymbol{\beta}+\# \gamma=I$, $ \partial_x^\gamma\partial_{\boldsymbol{v}}^{\boldsymbol{\beta}}\partial_{\mu}^{\alpha}X_{t}^{x,[\theta]}(\boldsymbol{v})$ exists and is c\`adl\`ag. In the meanwhile, it solves a linear equation of the form (\ref{first}). Moreover, $ \|\partial_x^\gamma\partial_{\boldsymbol{v}}^{\boldsymbol{\beta}}\partial_{\mu}^{\alpha}X_{t}^{x,[\theta]}(\boldsymbol{v})\|_{S_T^p}$ is bounded independently of ~$(x,[\theta],\boldsymbol{v}) $ for all ~$p\geq 2$.
	 		
$(S2)$ $\partial_x^\gamma\partial_{\boldsymbol{v}}^{\boldsymbol{\beta}}\partial_{\mu}^{\alpha}X_{t}^{x,[\theta]}(\boldsymbol{v})\in
	 		\mD^{M-1,\infty} $ and moreover,
	 		\ce
	 		\sup_{t\in[0,T],\boldsymbol{v}\in (\mR^N)^{\# \boldsymbol{v}}}\mathbb{E}\left[\left|{D}^{(M-I-1)}\partial_{x}^{\gamma}\partial_{\boldsymbol{v}}^{\boldsymbol{\beta}}\partial_{\mu}^{\alpha}X_{t}^{x,[\theta]}(\boldsymbol{v})\right|^{p}\right]\leq C(1+|x|+\|\theta\|_2)^{mp},
	 		\de
	 		for all ~$p\geq 2$, $m=1$; when the coefficients ~$c$ and $b$ are uniformly bounded, then ~$m=0$.
	 	
	 	$I=1$:
	 	
	 		$(S1)$ By Theorem \ref{dfo}$, \partial_x X_t^{x,[\theta]}$ and $ \partial_\mu X_t^{x,[\theta]}(v_1)$ exist.  Since they have the same form as $X_t^{x,[\theta]}$, so $ \partial_x X_t^{x,[\theta]}$ and $ \partial_\mu X_t^{x,[\theta]}(v_1)$ are c\`ald\`ag.  At this stage, no differentiation with respect to $v$ is required. We can denote
	 		\ce
	 		Y_t^{x,[\theta]}(v_1):=
	 		\begin{pmatrix}
	 			\partial_xX_t^{x,[\theta]} \\
	 			\partial_\mu X_t^{x,[\theta]}(v_1)
	 		\end{pmatrix}
	 		\de
	 		in the form of Eq. (\ref{first}) and identify the coefficients,
	 		{
	 		\ce
	 		a_0&=&
	 		\begin{pmatrix}
	 			1 \\
	 			0
	 		\end{pmatrix}, \\
	 		a_1^1(s,x,[\theta])&=&
	 		\begin{pmatrix}
	 			\partial b\left(X_s^{x,[\theta]},\left[X_s^{\theta}\right]\right) & 0 \\
	 			0 & \partial b\left(X_s^{x,[\theta]},\left[X_s^{\theta}\right]\right)
	 		\end{pmatrix}, \\
	 		a_1^2(s,x,u,[\theta])&=&
	 		\begin{pmatrix}
	 			\partial c\left(X_{s-}^{x,[\theta]},u,\left[X_{s-}^{\theta}\right]\right) & 0 \\
	 			0 & \partial c\left(X_{s-}^{x,[\theta]},u,\left[X_{s-}^{\theta}\right]\right)
	 		\end{pmatrix}, \\
	 		a_2^1(s,x,[\theta],v_1)&=&a_2^2(s,x,[\theta],v_1,u)=
	 		\begin{pmatrix}
	 			0  \\
	 			0
	 		\end{pmatrix}, \\
	 		a_3^1(s,x,[\theta],v_1) &=&
	 		\begin{pmatrix}
	 			0 & 0 \\
	 			\begin{aligned}
	 				&\partial_\mu b\left(X_s^{x,[\theta]},[X_s^{\theta}],(X_s^{v,[\theta]})^{\prime}\right)\mathbf{1}_{v = v_1}
	 			\end{aligned}
	 			&
	 			\begin{aligned}
	 				&\partial_\mu b\left(X_s^{x,[\theta]},[X_s^{\theta}],(X_s^{v,[\theta]})^{\prime}\right)\mathbf{1}_{v \ne v_1}
	 			\end{aligned}\\
	 		\end{pmatrix},\\
	 		a_3^2(s,x,u,[\theta],v_1) &=&
	 		\begin{pmatrix}
	 			0 & 0 \\
	 			\begin{aligned}
	 				& \partial_\mu c\left(X_{s-}^{x,[\theta]}, u, [X_{s-}^{\theta}], (X_{s-}^{v,[\theta]})^{\prime}\right)\mathbf{1}_{v = v_1}
	 			\end{aligned}
	 			&
	 			\begin{aligned}
	 				&\partial_\mu c\left(X_{s-}^{x,[\theta]}, u, [X_{s-}^{\theta}], (X_{s-}^{v,[\theta]})^{\prime}\right)\mathbf{1}_{v \ne v_1}
	 			\end{aligned}
	 			\end{pmatrix},\\
	 			a_3^1(s,x,[\theta],{\theta}^\prime) &=&
	 		\begin{pmatrix}
	 			0 & 0 \\
	 			0&
	 			\begin{aligned}
	 				&\partial_\mu b\left(X_s^{x,[\theta]},[X_s^{\theta}],(X_s^{\theta^{\prime},[\theta]})^{\prime}\right)
	 			\end{aligned}\\
	 		\end{pmatrix},\\
	 		a_3^2(s,x,u,[\theta],{\theta}^\prime) &=&
	 		\begin{pmatrix}
	 			0 & 0 \\
	 			0&
	 			\begin{aligned}
	 				& \partial_\mu c\left(X_{s-}^{x,[\theta]}, u, [X_{s-}^{\theta}], (X_{s-}^{\theta^{\prime},[\theta]})^{\prime}\right)
	 			\end{aligned}\\
	 		\end{pmatrix}.
	 		\de}
             It is obvious that the coefficients
	 	satisfy the assumptions of the  Lemma \ref{linear existence}, and we can obtain a bound on ~$\| Y^{x,[\theta]}(v_1)\|_{S_T^p}$.
	 		
Recall to the equations satisfied by ~$\partial_x X_t^{x,[\theta]}$ and ~$\partial_\mu X_t^{x,[\theta]}(v) $, we see that the coefficients are ~$(k-1)$ times differentiable with bounded Lipschitz derivatives, since the coefficients $b$ and ~$c$ are in ~$\cC_{b,Lip}^{k,k}$. \cite[Proposition8.2]{BL}  immediately tells us that  ~$\partial_x X_t^{x,[\theta]}, \partial_\mu X_t^{x,[\theta]}(v) \in \mD^{k-1,\infty}$,
  	 		\ce
	 		\mE[\sup_{s\in [0,T]} |D( Y_t^{x,[\theta]}(v))|^p] \leq C \mE[\sup_{s\in [0,T]} |Da_1|^p]
	 		\leq  C\mE[\sup_{s\in [0,T]} |\partial^2 c(X_s^{x,[\theta]},u,[X_s^{\theta}]) DX_{s}^{x,\theta}|^p].
	 		\de
Since $\partial^2 c $ is bounded  it suffices to  prove that
	 		\ce
	 		\mE[\sup_{s\in [0,T]} | DX_{s}^{x,\theta}|^p]\leq C(1+x+\|\theta\|_2)^{p}.
	 		\de
	 		
First, recall the Eq. (\ref{Md}) satisfied by the Malliavin derivative of $X_t^{x,[\theta]}$ we we can easily get by Kunita's second inequality,
	 		\ce
	 		\mE[\sup_{s\in [0,t]} | DX_{s}^{x,\theta}|^p]&\leq& C(1+|x|+\|\theta\|_2)^p+C\mE[\int_0^t  | DX_{h}^{x,\theta}|^p dh ]\\
	 		&\leq& C(1+|x|+\|\theta\|_2)^p+C\mE[\sup_{s\in [0,t]}\int_0^s  | DX_{h}^{x,\theta}|^p dh ].
	 		\de
By Gr\"onwall's inequality,
	 		\ce
	 		\mE[\sup_{s\in [0,t]} | DX_{s}^{x,\theta}|^p]&\leq&C(1+|x|+\|\theta\|_2)^p.
	 		\de
	 		We have thus established a bound for the first Malliavin derivative of $Y^{x,[\theta]}(v)$. For higher-order derivatives, the result in \cite[Proposition 8.19]{BL} indicates that each satisfies a linear equation. A key observation is that in these equations for higher-order derivatives, the coefficient $a_1^1$
	 		is  ~$ \partial b(X_s^{x,[\theta]},[X_s^\theta]) $ and the coefficient $a_1^2$
	 		is  ~$ \partial c(X_s^{x,[\theta]},[X_s^\theta]) $.
	 		From the bound on the Malliavin derivative for the general linear equation derived in (\ref{6.11}), it is evident that only a single term governs the estimate. Consequently, under the boundedness of all coefficients, the bound remains uniform with respect to the parameters ~$(x,[\theta],\boldsymbol{v})$.
	 		
	 		$2\leq I\leq k$:
	 		
	 			$(S1)$ By the induction hypothesis, for any ~$\alpha,\boldsymbol{\beta},\gamma$ satisfying ~$\#\alpha+\# \boldsymbol{\beta}+\# \gamma=I$, let $ Y_{t}^{x,[\theta]}(\boldsymbol{v}):=\partial_{x}^{\gamma} \partial_{\boldsymbol{v}}^{\boldsymbol{\beta}} \partial_{\mu}^{\alpha} X_{t}^{x,[\theta]}(\boldsymbol{v})$ and $Z_t^{x,[\theta]}(\boldsymbol{v},{v}^{\prime})$ be its derivatives which have the form
	 			\ce
	 			Z_t^{x,[\theta]}(\boldsymbol{v},{v}^{\prime}) &=& + \int_0^t \{ b_1^1(s, x, [\theta]) Z_{s}^{x,[\theta]}(\boldsymbol{v},{v}^{\prime}) + b_2^1(s, x, [\theta], \boldsymbol{v}) \no\\
	 		&&+  {\mathbb{E}^{\prime}}[ b_3^1(s, x, [\theta], {\theta^{\prime}}) ({Y}_{s}^{{\theta^{\prime}},[\theta]}(\boldsymbol{v},{v}^{\prime}))^{\prime}
	 		+ \sum_{r=1}^{\#\boldsymbol{v}} b_3^1(s, x, [\theta], v_r) ({Z}_{s}^{v_r,[\theta]}(\boldsymbol{v},{v}^{\prime}) )^{\prime}]\} ds\no\\
	 		&&+ \int_0^t\int_{\Xi} \{ b_1^2(s, x,u, [\theta]) Z_{s-}^{x,[\theta]}(\boldsymbol{v},{v}^{\prime}) + b_2^2(s, x,u, [\theta], \boldsymbol{v}) \no\\
	 		&&+  {\mathbb{E}^{\prime}}[ b_3^2(s, x,u, [\theta], {\theta^{\prime}}) ({Z}_{s-}^{{\theta^{\prime}},[\theta]}(\boldsymbol{v},{v}^{\prime}))^{\prime}\\
	 		&&+ \sum_{r=1}^{\#\boldsymbol{v}} b_3^2(s, x, u,[\theta], v_r) ({Z}_{s-}^{v_r,[\theta]}(\boldsymbol{v},{v}^{\prime}) )^{\prime}]\} \widetilde{N}(ds,du).
	 			\de
	 	        Then denote
	 			\ce
	 			Z_t^{x,[\theta]}(\boldsymbol{v},v^{\prime}):=
	 			\begin{pmatrix}
	 				\partial_xY_t^{x,[\theta]}(\boldsymbol{v}) \\
	 				\partial_\mu Y_t^{x,[\theta]}(\boldsymbol{v},{v}^{\prime}) \\
	 				\partial_{v_j}Y_t^{x,[\theta]}(\boldsymbol{v})
	 			\end{pmatrix}.
	 			\de
	 			By Lemma \ref{linear second order}, we can identify these new coefficients ~$b_1^i,b_2^i,b_3^i,i=1,2$ as
	 			\begin{align*}
b_{1}^1(s, x,[\theta]) &= \partial b\left(X_{s}^{x,[\theta]},\left[X_{s}^{\theta}\right]\right) I_{3},\quad b_1^2(s,x,u,[\theta]) = \partial c\left(X_{s-}^{x,[\theta]},u,\left[X_{s-}^{\theta}\right]\right)I_{3}, \\
b_{2}^i(s, x,u,[\theta], v) &= \begin{pmatrix}
b_{21}^i \\
b_{22}^i \\
b_{23}^i
\end{pmatrix}, \\
b_{3}^1(s, x,[\theta], v_r) &= \begin{pmatrix}
0 & 0 & 0 \\
a_{3}^1(s, x,[\theta], v_r) \mathbf{1}_{v_r=v^{\prime}} & a_{3}^1(s, x,[\theta], v_r) & 0 \\
a_{3}^1(s, x,[\theta], v_r) \mathbf{1}_{v_r=v_{j}} & 0 & a_{3}^1(s, x,[\theta], v_r)
\end{pmatrix}, \\
b_{3}^2(s, x,u,[\theta], v_r) &= \begin{pmatrix}
0 & 0 & 0 \\
a_{3}^2(s, x,u,[\theta], v_r) \mathbf{1}_{v_r=v^{\prime}} & a_{3}^2(s, x,u,[\theta], v_r) & 0 \\
a_{3}^2(s, x,u,[\theta], v_r) \mathbf{1}_{v_r=v_{j}} & 0 & a_{3}^2(s, x,u,[\theta], v_r)
\end{pmatrix},\\
b_{3}^1(s, x,[\theta], \theta^{\prime}) &= \begin{pmatrix}
0 & 0 & 0 \\
0 & a_{3}^1(s, x,[\theta], \theta^{\prime}) & 0 \\
0 & 0 & a_{3}^1(s, x,[\theta], \theta^{\prime})
\end{pmatrix}, \\
b_{3}^2(s, x,u,[\theta], \theta^{\prime}) &= \begin{pmatrix}
0 & 0 & 0 \\
0 & a_{3}^2(s, x,u,[\theta], \theta^{\prime}) & 0 \\
0 & 0 & a_{3}^2(s, x,u,[\theta], \theta^{\prime})
\end{pmatrix},
\end{align*}
	 			where
	 			\ce
	 			b_{21}^i&=&\partial_{x} a_{1}^i Y_{s}^{x,[\theta]}(\boldsymbol{v})+\partial_{x} a_{2}^i+{\mathbb{E}^{\prime}}[\partial_{x} a_{3}^i|_{v=\theta^{\prime}} ({Y}_{s}^{\theta^{\prime},[\theta]}(\boldsymbol{v}))^{\prime}+\sum_{r=1}^{\# \boldsymbol{v}} \partial_{x} a_{3}^i|_{v=v_{r}} ({Y}_{s}^{v_{r},[\theta]}(v))^{\prime}],\\
	 			b_{22}^i&=&\partial_{\mu} a_{1}^i Y_{s}^{x,[\theta]}(\boldsymbol{v})+\partial_{\mu} a_{2}^i+{\mathbb{E}^{\prime}}[\partial_{\boldsymbol{v}} a_{3}^i|_{v=v^{\prime}} Y_{s}^{v^{\prime},[\theta]}(\boldsymbol{v})+\partial_{\mu} a_{3}^i|_{v=\theta^{\prime}} Y_{s}^{\theta^{\prime},[\theta]}(\boldsymbol{v})+\sum_{r=1}^{\# \boldsymbol{v}} \partial_\mu a_3^i|_{v=v_r} (Y_s^{v_r,[\theta]})^{\prime}],\\
	 			b_{23}^i&=&\partial_{v_{j}} a_{2}^i+{\mathbb{E}^{\prime}}[\partial_{v_{j}} a_{2}^i ({Y}_{s}^{v_{j},[\theta]}(\boldsymbol{v}))^{\prime}].
	 			\de
	 			Now in order to obtain  boundedness on the ~$S_T^p$-norm of ~$Z^{x,[\theta]}(v,v^{\prime})$, we just need to check that the coefficients ~$b_1,b_2,b_3$ satisfy the assumptions of Lemma \ref{linear existence}, which is straightforward.

	 			$(S2)$  This is the same as the case ~$I=1$.
	 \end{proof}
	
    \subsection{Proofs of Lemma \ref{guocheng}}\label{appB}
	 In order to prove Lemma \ref{guocheng}, we need the following lemmas. The function belonging to the set~$\mK_r^q(E,M)$ satisfy the following properties, which we make use of when developing integration by parts formulas in Section 4.
	
	 \bl\label{pro} The following statements hold:
	 \begin{enumerate}
	 	\item If ~$\Psi_i\in \mK_{r_i}^{q_i}(\mR,M_i)$ for ~$i=1,\cdots,n,$ then
	 	\ce
	 	\begin{aligned}
	 		\prod_{i=1}^n\Psi_i\in\mathbb{K}_{r_1+...+r_n}^{q_1+\cdots+q_n}(\mR,\min_iM_i), & & & & \sum_{i=1}^n\Psi_i\in\mathbb{K}_{\min_ir_i}^{\max_iq_i}(\mR,\min_iM_i).
	 	\end{aligned}
	 	\de
	 	\item If~$\Psi\in \mK_r^q(\mR,M)$, then
	 	\ce
	 	\partial_x\Psi(x,\mu)\in\mathbb{K}_r^q(\mathbb{R},M-1)
	 	\de
	 	and
	 	\ce
	 	\partial_\mu\Psi(x,\mu,\boldsymbol{v}) \in\mK_r^q(\mR,M-1).
	 	\de
	 \end{enumerate}
	 \el
	 \begin{proof}
	 	The proof is similar to \cite{CM}, so we omit it here.
	 \end{proof}
	
	 \bl\label{A}
	 Under Assumption \ref{R}, for all $t\in[0,T]$, $p\geq 2$ and all $ i\in\{ 1,\cdots,N\}$, we have
	 	\ce
	 &&A[X_{i,t}^{x,[\theta]}]\\
	 &=&\int_0^t \int_{\Xi} (\sum_{j=1}^{N} \frac{\partial c_{i}}{\partial x_{j}}( X_{s-}^{x,[\theta]}, u,[X_{s-}^{\theta}]) A[X_{j, s-}^{x,[\theta]}] \widetilde{N}(d s, d u)\\
	 &&+\int_0^t \int_{\Xi}\frac{1}{2} \sum_{j, k=1}^{N} \frac{\partial^{2} c_{i}}{\partial x_{j} \partial x_{k}}( X_{s-}^{x,[\theta]}, u,[X_{s-}^{\theta}]) \Gamma[X_{j, s-}^{x,[\theta]}, X_{k, s-}^{x,[\theta]}]) \widetilde{N}(d s, d u)\\
	 &&+\int_{0}^{t} \int_{\Xi}\sum_{j=1}^{N} a[c_{i}(X_{s-}^{x,[\theta]}, \cdot,[X_{s-}^{\theta}])](u) \widetilde{N}(d s, d u) \\
	 &&+\int_{0}^{t}(\sum_{j=1}^{N} \frac{\partial b_{i}}{\partial x_{j}}(X_{s}^{x,[\theta]}, [X_{s}^{\theta}]) A[X_{j,s}^{x,[\theta]}]+\frac{1}{2} \sum_{j, k=1}^{N} \frac{\partial^{2} b_{i}}{\partial x_{j} \partial x_{k}}(X_{s}^{x,[\theta]},[X_{s}^{\theta}]) \Gamma[X_{j,s}^{x,[\theta]}, X_{k,s}^{x,[\theta]}]) ds
	 \de
	 and
	 \ce
	 \mE[A[X_{i,t}^{x,[\theta]}]^p] \leq Ct(1+|x|+\|\theta\|_2)^p.
	 \de
	 \el
	 \begin{proof}
	 The explicit expression for $A[ X_t^{x,[\theta]}]$  has already been computed in \cite[Proposition 8.26]{BL}.
	 By Kunita's second inequality and the assumptions on the cofficients, we have
	 \ce
	 &&\mE[A[X_{i,t}^{x,[\theta]}]^p]\\
	 &\leq & \mE[|\int_0^t \int_{\Xi} (\sum_{j=1}^{N} \frac{\partial c_{i}}{\partial x_{j}}( X_{s-}^{x,[\theta]}, u,[X_{s-}^{\theta}]) A[X_{j, s-}^{x,[\theta]}]\widetilde{N}(ds,du)|^p] \\
	 &&+\mE[|\frac{1}{2} \sum_{j, k=1}^{N} \int_0^t \int_{\Xi}\frac{\partial^{2} c_{i}}{\partial x_{j} \partial x_{k}}( X_{s-}^{x,[\theta]}, u,[X_{s-}^{\theta}]) \Gamma[X_{j, s-}^{x,[\theta]}, X_{k, s-}^{x,[\theta]}]) \widetilde{N}(d s, d u)|^p]\\
	 &&+\mE[|\int_{0}^{t} \int_{\Xi} a[c_{i}(X_{s-}^{x,[\theta]}, \cdot, [X_{s-}^{\theta}])](u) \widetilde{N}(d s, d u)|^p] \\
	 &&+\mE[|\int_{0}^{t}(\sum_{j=1}^{N} \frac{\partial b_{i}}{\partial x_{j}}(X_{s}^{x,[\theta]}, [X_{s}^{\theta}]) A[X_{j,s}^{x,[\theta]}]ds|^p]\\
	 &&+\mE[ |\frac{1}{2}\sum_{j, k=1}^{N} \int_0^t\frac{\partial^{2} b_{i}}{\partial x_{j} \partial x_{k}}(X_{s}^{x,[\theta]},[X_{s}^{\theta}]) \Gamma[X_{j,s}^{x,[\theta]}, X_{k,s}^{x,[\theta]}])ds|^p]\\
	 &\leq & C\sum_{j=1}^{N}(\mE[\int_0^t |A[X_{j, s-}^{x,[\theta]}]|^p ds] + \mE[\int_0^t |\Gamma[X_{j, s-}^{x,[\theta]}, X_{k, s-}^{x,[\theta]}])|^p ds]\\
&&+ \mE[ \int_0^t\int_{\Xi} |a[c_{i}(X_{s-}^{x,[\theta]}, \cdot, [X_{s-}^{\theta}])](u)|^p \lambda(du)ds]).
	 \de
	 Under Assumption \ref{R}, we have
$\partial_x^{\gamma}\partial_{\boldsymbol{v}}^{\boldsymbol{\beta}} \partial_\mu^{\alpha} c(X_s^{x,[\theta]},u,[X_s^{\theta}]) \in \bar{\mathbf{d}}^{m,p} $ and
\ce
\|\partial_x^\gamma \partial_{\boldsymbol{v}}^{\boldsymbol{\beta}} \partial_{\mu}^{\alpha} c(X_s^{x,[\theta]},u,[X_s^{\theta}])) \|_{\bar{\mathbf{d}}^{m,q}}\leq C.
\de
 Then there exists $\eta(u)\in L^2(\mP)$ such that
    \ce
    |a[c_i( X_s^{x,[\theta]}, \cdot,[X_s^{\theta}])](u)-a[c_i(0, \cdot,\delta_0)](u)|\leq C \eta(u)(| X_s^{x,[\theta] }|+W_2(X_s^{\theta},\delta_0)).
    \de
    Therefore
	 \ce
	 &&\mE[ \int_0^t\int_{\Xi} |a[c_{i}(X_{s-}^{x,[\theta]}, \cdot, [X_{s-}^{\theta}])](u)|^p \lambda(du)ds]\\
&=&\int_0^t \int_{\Xi}\mE[ |a[c_{i}(X_{s-}^{x,[\theta]}, \cdot, [X_{s-}^{\theta}])](u)|^p ]\lambda(du)ds\\
&\leq& \int_0^t C(1+|x|+\|\theta\|_2)^p ds=Ct(1+|x|+\|\theta\|_2)^p.
	 \de
Since $\Gamma[X_t^{x,[\theta]}]=\widehat{\mE}[X_t^{x,[\theta]}]$ and $\|DX_t^{x,\theta}\|_p^p \leq C(1+|x|+\|\theta\|_2)^p$. Then
	 \ce
	 &&\mE[\int_0^t |\Gamma[X_{j, s-}^{x,[\theta]}, X_{k, s-}^{x,[\theta]}])|^p ds]=\int_0^t \mE[|\Gamma[X_{j, s-}^{x,[\theta]}, X_{k, s-}^{x,[\theta]}])|^p ]ds\\
	 &=&\int_0^t \mE[|\widehat{\mE}[DX_{j, s-}^{x,[\theta]} DX_{k, s-}^{x,[\theta]}]|^p ]ds\\
	 &\leq & \int_0^t \mE[(\widehat{\mE}[|DX_{j, s-}^{x,[\theta]}|^{2p}])^{\frac{1}{2}} (\widehat{\mE}[|DX_{k, s-}^{x,[\theta]}|^{2p})^{\frac{1}{2}} ]ds\\
	 &\leq & \int_0^t (\mE[\widehat{\mE}[|DX_{j, s-}^{x,[\theta]}|^{2p}])^{\frac{1}{2}} (\mE[\widehat{\mE}[|DX_{k, s-}^{x,[\theta]}|^{2p}]) ]^{\frac{1}{2}}ds\\
	 &\leq & C \int_0^t  (1+|x|+\|\theta\|_2)^p ds\leq  Ct(1+|x|+\|\theta\|_2)^p.
	 \de
	 Therefore, we have
	 \ce
	 \mE[A[X_{i,t}^{x,[\theta]}]^p]
	 &\leq& C\mE[\int_0^t |A[X_{j, s-}^{x,[\theta]}]|^p ds] \\
&&+C \mE[\int_0^t |\Gamma[X_{j, s-}^{x,[\theta]}, X_{k, s-}^{x,[\theta]}])|^p ds]+C \mE[ \int_0^t |a[c_{i}(X_{s-}^{x,[\theta]}, u,[X_{s-}^{\theta}], \cdot)](u)|^p ds]\\
	 &\leq & C\int_0^t \mE[|A[X_{j, s-}^{x,[\theta]}]|^p] ds +Ct+Ct(1+|x|+\|\theta\|_2)^p\\
	 &\leq & C(t(1+|x|+\|\theta\|_2)^p+ \int_0^t \mE[|A[X_{j, s-}^{x,[\theta]}]|^p] ds).
	 \de
	 By Gr\"onwall's ineqaulity,
	 \ce
	 \mE[A[X_{i,t}^{x,[\theta]}]^p] \leq Ct(1+|x|+\|\theta\|_2)^pe^{CT} \leq Ct(1+|x|+\|\theta\|_2)^p.
	 \de	
	 \end{proof}

	  Now we can give a detailed proof of Lemma \ref{guocheng}.
	  \begin{proof}
	(1) Since
$$\langle Du,\Psi \rangle_H=-\delta(u\cdot \Psi)+u\delta(\Psi),$$
using \cite[Lemma 7.20]{BL}, we have
	  		\ce
	  		Z(\Psi)(t,x,u,[\theta])&:=&\delta( \big( (DX_t^{x,[\theta]})^* (\Gamma[X_t^{x,[\theta]}])^{-1}\cdot \partial_x X_t^{x,[\theta]}\big)_i \Psi(t,x,u,[\theta]))\\
	  		&=&-\langle D\Psi,\big((DX_t^{x,[\theta]})^* (\Gamma[X_t^{x,[\theta]}])^{-1} \partial_x X_t^{x,[\theta]}\big)_i \rangle_{H}\\
&&+\Psi \delta(\big((DX_t^{x,[\theta]})^* (\Gamma[X_t^{x,[\theta]}])^{-1} \partial_x X_t^{x,[\theta]}\big)_i)\\
	  		&=&-\langle D\Psi,\big((DX_t^{x,[\theta]})^* (\Gamma[X_t^{x,[\theta]}])^{-1} \partial_x X_t^{x,[\theta]}\big)_i\rangle_{H}\\
&&+2\Psi \big(A[(X_t^{x,[\theta]})^* ] (\Gamma[X_t^{x,[\theta]}])^{-1} \partial_x X_t^{x,[\theta]}\big)_i.
	  		\de
	  		Combining the facts $ DX_t^{x,[\theta]}\in \mK_0^1(\mR^N;k-1) $, $\partial_x X_t^{x,[\theta]}\in \mK_0^1(\mR^{N\times N};k-1)$, $D\Psi\in \mK_r^q(\mR;n-1)$ and the Lemma \ref{A} and \ref{G},  we have  by H\"older's inequality the following estimate
	  		\ce
	  		&&\mE[|Z(t,x,u,[\theta])|^p ]\\
	  		&=&\mE[|-\langle D\Psi,[(DX_t^{x,[\theta]})^* (\Gamma[X_t^{x,[\theta]}])^{-1} \partial_x X_t^{x,[\theta]}]_1 \rangle_{H}\\
	  		&&+\Psi\cdot \delta(\big((DX_t^{x,[\theta]})^*] (\Gamma[X_t^{x,[\theta]}])^{-1} \partial_x X_t^{x,[\theta]}\big)_i )|^p ]\\
	  		&\leq & C\mE[ |\langle D\Psi,\big((DX_t^{x,[\theta]})^* (\Gamma[X_t^{x,[\theta]}])^{-1} \partial_x X_t^{x,[\theta]}\big)_i) \rangle_{H} |^p]\\
	  		&&+\mE[| \Psi\cdot 2\big(A[(X_t^{x,[\theta]})^*](\Gamma[X_t^{x,[\theta]}])^{-1} \partial_x X_t^{x,[\theta]}\big)_i|^p]\\
	  		&\leq & C(\mE[ |D\Psi|^{4p} ])^{\frac{1}{4}} (\mE[| D(X_t^{x,[\theta]})^*|^{4p} ])^{\frac{1}{4}}(\mE[| (\Gamma[X_t^{x,[\theta]}])^{-1}|^{4p} ])^{\frac{1}{4}} (\mE[| \partial_x X_t^{x,[\theta]}|^{4p} ])^{\frac{1}{4}}\\
	  		&+&C(\mE[|\Psi|^{4p}])^{\frac{1}{4}}(\mE[|A[(X_t^{x,[\theta]})^*]|^{4p}])^{\frac{1}{4}}(\mE[(\Gamma[X_t^{x,[\theta]}])^{-4p}])^{\frac{1}{4}}(\mE[|\partial_x X_t^{x,[\theta]}|^{4p}])^{\frac{1}{4}}\\
	  		&\leq & Ct^{\frac{rp}{2}-\frac{p}{a}}(1+|x|+\|\theta\|_2)^{(2+q)p}+Ct^{\frac{rp+2}{2}-\frac{p}{a}}(1+|x|+\|\theta\|_2)^{(2+q)p}\\
	  		&\leq & 2Ct^{\frac{rp}{2}-\frac{p}{a}}(1+|x|+\|\theta\|_2)^{(2+q)p}.
	  		\de
	  			It is easy to check that the order of differentiability of $Z(\Psi)(t,x,[\theta])$ is $(k-2)\wedge (n-1)-1$. And by induction, the order of differentiability of $Z_{(i)}^1(\Psi)(t,x,[\theta])$ is $(k-2)\wedge (n-1)-1$.
	  		
	  		Then, based on the form of the derivative computed in the Theorem \ref{dfo}, the proof of higher derivatives in  Theorem \ref{KS} and by applying the same calculation method as above, we can obtain, for all $m+\#\alpha+\# \boldsymbol{\beta}+\# \gamma\leq (k-2)\wedge (n-1)$,
	  		\ce
	  		\mE[|D^m \partial_x^{\gamma} \partial_{\boldsymbol{v}}^{\boldsymbol{\beta}} \partial_\mu^{\alpha}Z_{(i)}^1(\Psi)|^p ]\leq Ct^{\frac{rp}{2}-\frac{p}{a}}(1+|x|+\|\theta\|_2)^{(2+q)p}.
	  		\de
	  		Therefore, $Z_{(i)}^1(\Psi)(t,x,[\theta])\in \mK_{r-\frac{2}{a}}^{q+2}(\mR;(k-2)\wedge (n-1)-1) $. Then repeating the above steps, we can obtain $Z_{\gamma}^1(\Psi)(t,x,[\theta])\in \mK_{r-\frac{2}{a}\# \gamma}^{q+2\# \gamma}(\mR;[(k-2)\wedge (n-1)]-\# \gamma) $.
	
(2) Similar to the proof above,  $Z_{\gamma}^2(\Psi)(t,x,[\theta])\in \mK_{r-\frac{2}{a}\# \gamma}^{q+\# \gamma}(\mR;[(k-2)\wedge (n-1)]-\# \gamma) $.

(3) From Lemma \ref{pro} (2), $\partial_x \Psi(t,x,[\theta]) \in \mK_r^q(\mR;n-1)$, we have by recurrence  $\partial_x^{\gamma} \Psi (t,x,[\theta])\in \mK_r^q(\mR;n-\# \gamma)$. Then by Lemma \ref{pro} (1),
	  		\ce
	  		Z_{\gamma}^3(\Psi)(t,x,[\theta])\in \mK_{r-\frac{2}{a}\# \gamma}^{q+2\# \gamma}(\mR;[(k-2)\wedge( n-1)]-\# \gamma).
	  		\de

(4) From the definition of~$Z_{\alpha,\beta}^4(\Psi)$,
	  		\ce
	  		Z_{\gamma,\alpha}^4(\Psi)(t,x,[\theta])&=&(Z_{\alpha}^2( Z_{\gamma}^1(\Psi))(t,x,[\theta]))+ Z_{\alpha}^2( \partial_x^\gamma \Psi)(t,x,[\theta]).
	  		\de
	  		
	  		Here, it suffices to substitute ~$ \Psi$ with ~$ Z_\gamma^1(\Psi)$ and ~$ \partial_x^{\gamma} \Psi$ respectively,
	  		\ce
	  		Z_{\alpha}^2(Z_{\gamma}^1(\Psi) )(t,x,[\theta]) \in \mK_{r-\frac{2}{a}(\# \gamma+\# \alpha)}^{ q+2\# \gamma+\# \alpha}(\mR, [(k-2) \wedge n]-\# \gamma+\# \alpha ),
	  		\de
	  		and
	  		\ce
	  		Z_{\alpha}^2(\partial_x^{\gamma} \Psi )(t,x,[\theta])\in\mK_{r-\frac{2}{a}\#\alpha}^{ q+\#\alpha}(\mR, [(k-2)
\wedge (n-\#\gamma])-\#\alpha ).
	  		\de
So by Lemma \ref{pro} (1),
	  		\ce
	  		Z^4_{\gamma,\alpha}(\Psi)(t,x,[\theta]) \in \mK_{ r-\frac{2}{a}(\# \gamma+\# \alpha)}^{ q+2\#\gamma+\# \alpha}(\mR, [(k-2) \wedge n]-\# \gamma-\# \alpha).
	  		\de

(5) The proof follows the same line of argument as in  Part ~$(1)$.
	
 (6) By Lemma \ref{pro} (2),
	  		\ce
	  		\partial_\mu^\alpha \Psi(t,x,[\theta],\boldsymbol{v})\in \mK_r^q(\mR,n-\#\alpha),
	  		\de
	  		and
	  		\ce
	  		Z_{\mu,\alpha}^1(\Psi)(t,x,[\theta],\boldsymbol{v})\in
\mK_{r-\frac{2}{a}\#\alpha}^{ q+2\#\alpha}(\mR, [(k-2) \wedge (n-1)]-\#\alpha ).
	  		\de
	  		By Lemma \ref{pro} (1),
	  		\ce
	  		Z_{\mu,\alpha}^2(\Psi)(t,x,[\theta],\boldsymbol{v})&\in&\mK_{r-\frac{2}{a}\#\alpha}^{ q+2\#\alpha}(\mR, [(k-2) \wedge (n-1)]-\#\alpha ).
	  		\de
(7) Since
	  		\ce
	  		Z_{\mu,\gamma,\alpha}^3(\Psi)(t,x,[\theta],\boldsymbol{v})&=&Z_{\gamma}^2(Z_{\mu,\alpha}^1(\Psi) ) (t,x,[\theta],\boldsymbol{v} )+Z_{\gamma}^2(\partial_{\mu}^{\gamma} \Psi)(t,x,[\theta],\boldsymbol{v}),
	  		\de
we have from the above result,
	  		\ce
	  		Z_{\mu,\alpha}^1(\Psi)(t,x,[\theta],\boldsymbol{v}) &\in& \mK_{r-\frac{2}{a}\#\alpha}^{ q+2\#\alpha}(\mR, [(k-2) \wedge (n-1)]-\#\alpha),\\
	  		Z_{\gamma}^2(\Psi)(t,x,[\theta])&\in&\mK_{r-\frac{2}{a}\#\gamma}^{ q+\#\gamma}(\mR, [(k-2) \wedge (n-1)]-\#\gamma ),\\
	  		\partial_\mu^{\alpha} \Psi(t,x,[\theta],\boldsymbol{v})&\in& \mK_r^q(\mR,n-\#\alpha).
	  		\de
Consequently,
	  		\ce
	  		Z_{\gamma}^2(Z_{\mu,\alpha}^1(\Psi) ) (t,x,[\theta],\boldsymbol{v} )&\in& \mK_{r-\frac{2}{a}(\#\gamma+\#\alpha)}^{q+\#\gamma+2\#\alpha}(\mR,[k-2-\#\gamma-\#\alpha] \wedge [n-1-\#\gamma-\#\alpha] ),\\
	  		Z_{\gamma}^2(\partial_{\mu}^{\alpha} \Psi)(t,x,[\theta],\boldsymbol{v})&\in&\mK_{r-\frac{2}{a}\#\gamma}^{ q+\#\gamma}(\mR, [k-2-\#\gamma] \wedge  [n-1-\#\gamma-\#\alpha] ).
	  		\de
	  		That is
	  		\ce
	  		Z_{\mu,\gamma,\alpha}^3(\Psi)(t,x,[\theta],\boldsymbol{v}) &\in& \mK_{r-\frac{2}{a}(\#\gamma+\#\alpha)}^{q+\#\gamma+2\#\alpha}(\mR, [(k-2) \wedge (n-1)]-\#\gamma-\#\alpha).
	  		\de
	 	  \end{proof}

\end{document}